\UseRawInputEncoding

\documentclass[11pt]{article}
\usepackage{}
\usepackage{amssymb}
\usepackage{amsfonts}
\usepackage{mathrsfs}
\usepackage{graphicx}
\usepackage{amsmath}
\usepackage{color}
\usepackage{amsthm}
\usepackage{epstopdf}
\usepackage{subfigure}
\usepackage{pifont,bm}
\usepackage{tikz}
\usepackage{enumerate}
\usepackage{mathrsfs}
\usepackage{cite}
\usepackage{hyperref}
\usepackage{booktabs}
\usepackage{diagbox}
\usepackage{makecell}

\theoremstyle{definition}

\makeatletter
\def\@biblabel#1{[#1]}
\makeatother

\makeatletter \@addtoreset{equation}{section}

\begin{document}

\begin{titlepage}
\title{\bf{Data-driven forward-inverse problems for Yajima-Oikawa system using deep learning with parameter regularization
\footnote{Corresponding authors.\protect\\
\hspace*{3ex} E-mail addresses: ychen@sei.ecnu.edu.cn (Y. Chen)}
}}
\author{Juncai Pu$^{a}$, Yong Chen$^{a,b,*}$\\
\small \emph{$^{a}$School of Mathematical Sciences, Shanghai Key Laboratory of Pure Mathematics and} \\
\small \emph{Mathematical Practice, East China Normal University, Shanghai, 200241, China} \\
\small \emph{$^{b}$College of Mathematics and Systems Science, Shandong University }\\
\small \emph{of Science and Technology, Qingdao, 266590, China} \\
\date{}}
\thispagestyle{empty}
\end{titlepage}
\maketitle

\vspace{-0.5cm}
\begin{center}
\rule{15cm}{1pt}\vspace{0.3cm}

\parbox{15cm}{\small
{\bf Abstract}\\
\hspace{0.5cm}
We investigate data-driven forward-inverse problems for Yajima-Oikawa (YO) system by employing two technologies which improve the performance of neural network in deep physics-informed neural network (PINN), namely neuron-wise locally adaptive activation functions and $L^2$ norm parameter regularization. Indeed, we not only recover three different forms of vector rogue waves (RWs) by means of three distinct initial-boundary value conditions in the forward problem of YO system, including bright-bright RWs, intermediate-bright RWs and dark-bright RWs, but also study the inverse problem of YO system by using training data with different noise intensity. In order to deal with the problem that the capacity of learning unknown parameters is not ideal when the PINN with only locally adaptive activation functions utilizes training data with noise interference in the inverse problem of YO system, thus we introduce $L^2$ norm regularization, which can drive the weights closer to origin, into PINN with locally adaptive activation functions, then find that the PINN model with two strategies shows amazing training effect by using training data with noise interference to investigate the inverse problem of YO system.

}

\vspace{0.5cm}
\parbox{15cm}{\small{

\vspace{0.3cm} \emph{Key words: $L^2$ norm parameter regularization, PINN, rogue waves, Yajima-Oikawa system}  \\

} }
\end{center}
\vspace{0.3cm} \rule{15cm}{1pt} \vspace{0.2cm}

\section{Introduction}
With the revolution of computer hardware equipment and software technology again and again, the increasing amount of data, model scale, accuracy, complexity and impact on the real world promote the continuous and successful application of deep learning in more and more practical problems \cite{Hinton2006,LeCun2015}. Currently, deep learning has achieved remarkable success in practical problems in various fields. In the field of object recognition, modern object recognition network can not only recognize at least 1000 different categories of objects, but also process rich high-resolution photographs without cropping photos near the objects to be recognized \cite{Krizhevsky2012,Cao2020}. The introduction of deep learning has a great impact on speech recognition, which makes the error rate of speech recognition drop sharply \cite{Hinton2012}. Furthermore, deep networks have also achieved spectacular successes for pedestrian detection and image segmentation \cite{LiGF2020,ZengNY2021}, as well as yielded superhuman performance in traffic sign classification \cite{Zhang2020}. Moreover, deep learning detonated a wide range of landing applications, such as language understanding \cite{Collobert2011}, medical imaging \cite{Alipanahi2015}, face recognition \cite{Sun2018}, video surveillance \cite{Shao2020} and forward and inverse problems of nonlinear partial differential equations \cite{Raissi2019}.

Deep feedforward network, also often called feedforward neural network (NN) or multilayer perceptron, is the quintessential deep learning model \cite{Haykin2008}. The universal approximation theorem points out a feedforward NN with a linear output layer and at least one hidden layer with any ``squashing'' activation function can approximate any Borel measurable function from one finite-dimensional space to another with any desired non-zero amount of error, provided that the network is given enough hidden units \cite{Hornik1989}. With the successful use of back-propagation in training deep NNs with internal representation and the popularity of back-propagation algorithms \cite{Rumelhart1986}, many optimization methods based on the idea of calculating gradient in back-propagation came into being, such as stochastic gradient descent \cite{Ruder2017}, Adam \cite{Kingma2014} and L-BFGS \cite{Liu1989}. Furthermore, automatic differentiation (AD), also called algorithmic differentiation or simply ``autodiff'', is a family of techniques similar to but more general than back-propagation for efficiently and accurately evaluating derivatives of numeric functions expressed as computer programs \cite{Baydin2018}. Recently, due to the general approximation ability of NN architecture \cite{Lagaris1998} and wide application of AD technology, after taking the NN space as an ansatz space for the solution of governing equation, a physics-informed neural network (PINN) has been successfully constructed to accurately solve both the forward problems, where the approximate solutions of governing partial differential equations are obtained, as well as the inverse problems, where parameters involved in the governing equation are discovered from the training data \cite{Raissi2019}. Moreover, Karniadakis research team has recently successfully applied the PINN framework to many physical problems \cite{KarniadakisGE2021}, such as discovering turbulence models from scattered/noisy measurements \cite{RaissiM2019}, fractional differential equations \cite{PangG2019}, high speed aerodynamic flows \cite{Mao2020}, heat transfer problems \cite{CaiS2021} and stochastic differential equation by generative adversarial networks \cite{YangL2020}.

The type selection of hidden units is extremely significant and difficult in the design process of NN, and the activation function plays an important role during the selection of hidden units due to the derivative of the loss function depends on the optimization parameters, which depend on the derivative of the activation function \cite{Haykin2008}. The activation function of common hidden units in NNs usually acts on affine transformation, and popular activation functions include rectified linear units, maxout units, logistic sigmoid and hyperbolic tangent function, as well as some unpublished activation functions perform just as well as the popular ones, such as sine and cosine functions \cite{Goodfellow-book2016}. The choice of activation function completely depends on the problem at hand in practical application. However, these activation functions are fixed in the training process of NN, which will greatly limit the performance of NN and the convergence speed of objective function. Based on the basic framework of PINN, Jagtap et al. proposed two different adaptive activation functions, that is global adaptive activation function and locally adaptive activation functions, to approximate smooth and discontinuous functions as well as solutions of linear and nonlinear partial differential equations by introducing scalable parameters into the activation function and adding a slope recovery term based on activation slope to the loss function of locally adaptive activation functions, it proved that the locally adaptive activation functions further improves the performance of the NN and speeds up the training process of the NN \cite{Jagtap2020,JagtapA2020}. Recently, we found that PINN with neuron-wise locally adaptive activation functions shows better training efficiency, robustness and accuracy when studying the forward problem of derivative nonlinear Schr\"{o}dinger equation \cite{PuND2021,PuWM2021}.

As is known to all, a central problem in machine learning is how to make an algorithm that will perform well not just on the training data, but also on new inputs. Many strategies, which are known collectively as regularization, are explicitly designed to reduce the test error in machine learning, possibly at the expense of increased training error \cite{Goodfellow-book2016}. Indeed, regularization has been used for decades prior to the advent of deep learning \cite{Poggio1990}. Many regularization approaches are based on limiting the capacity of models, such as NNs, linear regression, or logistic regression, by adding a parameter norm penalty $\Omega(\bm{\theta})$ to the objective function, of which the most common and simplest is $L^2$ parameter norm regularization. $L^2$ parameter norm penalty is usually called weight decay, ridge regression \cite{ZhangR2021} or Tikhonov regularization \cite{Chada2020}, which drives the weight closer to the origin by adding a regularization term to the objective function. Furthermore, Japtap et al. demonstrated that a gradient descent algorithm minimizing objective function with global adaptive activation function and $L^2$ regularization does not converge to a suboptimal critical point or local minimum with an appropriate initialization and learning rates by providing corresponding theoretical result \cite{Jagtap2020}. Therefore, a natural inspiration is to further enhance the performance of the novel PINN by introducing the parameter regularization strategy into the PINN with neuron-wise locally adaptive activation functions. As far as we know, the forward-inverse problems of nonlinear integrable systems have not been studied by employing PINN with neuron-wise locally adaptive activation functions and $L^2$ norm regularization at the same time.

With the rapid development of deep learning, predicting the generation and evolution of rogue wave (RW) by utilizing obtained initial boundary value data plays an important role. At present, the significant researches have been done to study the RWs of some classical single equations by means of PINN method \cite{PuND2021,PuWM2021,PuCPB2021,Peng2021}. However, compared with single dynamical equation, coupled systems usually allow energy transfer between their additional degrees of freedom, which potentially generates families of intricate vector RWs, that is the coupled systems can describe RWs more accurately than the single model \cite{Manakov1974}. Among coupled field dynamics systems, the coupled long wave-short wave resonance equation (LSWR) \cite{Benney1977}, which also be called one-dimensional Yajima-Oikawa (YO) system \cite{Yajima1976}, is a fascinating nonlinear physical system, it describes a resonant interaction between long wave in complex envelope of rapidly varying field and short wave in real low-frequency field. Once the resonance condition is satisfied, that is, the group velocity of a short wave (high-frequency wave) exactly or almost matches the phase velocity of a long wave (low-frequency wave), this coupled system can be derived from the Davey-Stewartson system \cite{Djordjevic1977}. In 1972, Zakharov made a theoretical investigation for LSWR for the first time when analyzing the Langmuir wave in plasma \cite{Zakharov1972}. In the case of long wave unidirectional propagation, the general Zakharov system was reduced to YO system \cite{Yajima1976}. Surprisingly, despite its simple form, this system can describe various nonlinear wave phenomena, such as capillary-gravity wave in fluid \cite{Djordjevic1977}, optical-terahertz waves in second-order nonlinear negative refractive index medium \cite{Chowdhury2008}, while the resonance interaction of this system can occur between long and short internal waves \cite{Grimshaw1977}, and between a long internal wave and a short surface wave in a two layer fluid \cite{Funakoshi1983}. Different from other classical coupled systems, such as Manakov system, YO system is coupled by two completely different types of equations, which plays an irreplaceable role in many physical application scenarios. This paper is committed to studying the data-driven RWs and learning unknown parameter for this unique and important physical system.

Recently, in order to build appropriate PINN with neuron-wise locally adaptive activation functions suitable for coupled systems, we have proposed a PINN algorithm to obtain the data-driven vector localized waves including vector solitons, breathers and RWs for Manakov system in complex space \cite{PuJC2021}. In this paper, we will employ the regularization strategy into the PINN with neuron-wise locally adaptive activation functions to construct an improved PINN with three outputs and three physical constraints for studying the initial boundary value problem of YO system as follow
\begin{align}\label{E1}
\begin{split}
\begin{cases}
\mathrm{i}S_{t}+\lambda_1S_{xx}+SL=0,\,x\in[X_0,X_1],\, t\in[T_0,T_1],\\
L_{t}=\lambda_2(|S|^2)_x,\,x\in[X_0,X_1],\, t\in[T_0,T_1],\\
S(x,T_0)=S^0(x),\,L(x,T_0)=L^0(x),\,x\in[X_0,X_1],\\
S(X_0,t)=S^{\mathrm{lb}}(t),\,S(X_1,t)=S^{\mathrm{ub}}(t),\,t\in[T_0,T_1],\\
L(X_0,t)=L^{\mathrm{lb}}(t),\,L(X_1,t)=L^{\mathrm{ub}}(t),\,t\in[T_0,T_1],\\
\end{cases}
\end{split}
\end{align}
here the short wave component $S(x,t)$ stands for the complex envelope of the rapidly varying field and the long wave component $L(x,t)$ represents the real low-frequency field, in which $x$ and $t$ are two independent evolution variables. While $\lambda_1$ and $\lambda_2$ are real valued parameters, which can be known parameters or unknown parameters to be learned. The $|S|$ represents the module of complex valued short wave $S$, which also means $|S|^2=SS^*$ with $S^*$ indicates the conjugate of $S$. For physics discussions, The first equation of YO system \eqref{E1} is arranged in a form similar to the standard nonlinear Schr\"{o}dinger equation, which clearly indicates that its nonlinearity is driven by long wave field $L$ rather than Kerr term $|S|^2$. Whereas, the second equation of YO system is the KdV equation without nonlinear term, in which the dispersion term is driven by the module $|S|$ of short wave. We note that the RWs of this system have been effectively obtained with the aid of Hirota bilinear method \cite{Chow2013}, KP hierarchy reduction method \cite{ChenJC2018} and Darboux transformation \cite{ChenSH2014}.

This paper is organized as follows: after the introduction in section 1, section 2 gives a brief discussion of the improved PINN methodology with locally adaptive activation functions and $L^2$ parameter norm regularization for the coupled YO systems, where we also discuss about training data, loss function, parameter regularization, optimization methods and the operating environment. The algorithm flow schematic and algorithm steps for the YO system are also exhibited in detail. Section 3 provides the results and detailed discussions for forward problems on improved PINN approximations of data driven vector RWs in three different states. Section 4 presents experimental results with different trade-off norm penalty term coefficients in inverse problems. Finally, we summarize the conclusions of our work are given out in last section.

\section{Methodology}

Although the classical PINN is a general and efficient deep learning framework for solving forward and inverse problems involving nonlinear partial differential equations. However, during the training process of dealing with complex RWs of complex nonlinear systems, the classical PINN algorithm may show slow training speed, obvious fluctuation of loss function and unsatisfactory approach effect. Indeed, the PINN model with neuron-wise locally adaptive activation functions and slope recovery term can improve the convergence speed, stability of the loss function and approximation ability in the training process from Ref. \cite{JagtapA2020,PuND2021}. Furthermore, during the process of studying the inverse problem of YO system, we find that the PINN method with only neuron-wise locally adaptive activation functions has good training effect by means of clean data, but the network displays poor parameters learning capacity in the case of data training with noise interference. Therefore, considering that parameter regularization can modify the weight and reduce the influence of large weight on the network, it may enhance the performance of NN when learning unknown parameters with noisy training data, so we further improve the deep learning algorithm by introducing $L^2$ norm regularization into the PINN method with neuron-wise locally adaptive activation functions.

\subsection{NN and adaptive activation function}

We establish a NN of depth $D$ with an input layer, $D-1$ hidden-layers and an output layer, in which the $d$th hidden-layer contain $N_d$ number of neurons. Each NN hidden-layer accepts an output $\textbf{x}^{d-1}\in\mathbb{R}^{N_{d-1}}$ from the previous layer, where an affine transformation can be yielded as follows form
\begin{align}\label{E2}
\mathcal{L}_d(\textbf{x}^{d-1})\triangleq\textbf{W}^d\textbf{x}^{d-1}+\textbf{b}^d,
\end{align}
in which the network weights term $\textbf{W}^{d}\in\mathbb{R}^{N_d\times N_{d-1}}$ and bias term $\textbf{b}^d\in\mathbb{R}^{N_d}$ associated with the $d$th layer. Eventually, in order to introduce adaptive activation function, the $d$th layer neuron-wise locally adaptive activation functions are defined as bellow
\begin{align}\nonumber
\sigma\Big(na^d_i\big(\mathcal{L}_d(\textbf{x}^{d-1})\big)_i\Big),\,i=1,2,\cdots,N_d,
\end{align}
where $\sigma$ is the activation function, $n>1$ are scaling factors and $\{a^d_i\}$ are additional $\sum\limits_{d=1}^{D-1}N_d$ parameters to be optimized. Whereas, $n$ exist a critical scaling factor $n_{c}$, the NN optimization algorithm becomes sensitive as $n\geqslant n_c$ in each problem set \cite{JagtapA2020}. The introduction of neuron-wise locally adaptive activation makes all neurons in each hidden layer have own activation function slope. The improved NN with neuron-wise locally adaptive activation functions can be represented as
\begin{align}\label{E3}
\textbf{q}(\textbf{x};\bar{\Theta})=\Big(\left(\mathcal{L}_D\right)_{i'}\circ\sigma\circ na^{D-1}_{i}\left(\mathcal{L}_{D-1}\right)_{i}\circ\cdots\circ\sigma\circ na^1_i\left(\mathcal{L}_1\right)_i\Big)(\textbf{x}),\,i'=1,2,3,
\end{align}
where $\textbf{x}$ and $\textbf{q}(\textbf{x};\bar{\Theta})$ represent the two inputs and three outputs in the NN, respectively. The set of trainable parameters $\bar{\Theta}\in\bar{\mathcal{P}}$ consists of $\big\{\textbf{W}^d,\textbf{b}^d\big\}_{d=1}^{D}$ and $\big\{a_i^d\big\}_{d=1}^{D-1},\forall i=1,2,\cdots,N_d$ (in fact, $\bar{\mathcal{P}}$ also includes parameters to be predicted in the inverse problem), that is $\bar{\mathcal{P}}$ is the parameter space.

\subsection{YO system constraint}

Especially, we consider the (1 + 1)-dimensional coupled YO system as the physical constraint of aforementioned NN to construct PINN, its specific operator representations reduced from YO system \eqref{E1} are as shown below
\begin{align}\label{E4}
\begin{split}
&S_{t}+\mathcal{N}[S,L;\lambda_1]=0,\\
&L_{t}+\mathcal{N}'[S;\lambda_2]=0,
\end{split}
\end{align}
where $S$ and $L$ are complex valued solution and real valued solution of $x$ and $t$ to be determined later respectively, $\mathcal{N}[\cdot,\cdot;\lambda_1]$ and $\mathcal{N}'[\cdot;\lambda_2]$ are nonlinear differential operators with unknown parameters $\lambda_1$ and $\lambda_2$ in space. Due to the complexity of the structure of the complex-valued solution $S(x,t)$ in Eq. \eqref{E3}, we decompose $S(x,t)$ into the real part $u(x,t)$ and the imaginary part $v(x,t)$ by employing real-valued functions $u(x,t)$ and $v(x,t)$, that is $S(x,t)=u(x,t)+\mathrm{i}v(x,t)$. Then substituting it into Eq. \eqref{E4}, we have
\begin{align}\label{E5}
\begin{split}
&u_t+\mathcal{N}_u[u,v,L;\lambda_1]=0,\quad v_t+\mathcal{N}_v[u,v,L;\lambda_1]=0,\quad L_t+\mathcal{N}'_L[u,v;\lambda_2]=0.\\
\end{split}
\end{align}
Accordingly, the $\mathcal{N}_u$, $\mathcal{N}_v$ and $\mathcal{N}'_L$ are nonlinear differential operators in space. Then $f_u(x,t)$, $f_v(x,t)$ and $f_L(x,t)$ constitute the physics-informed parts of the NN, which can be defined as
\begin{align}\label{E6}
\begin{split}
&f_u:=u_t+\mathcal{N}_u[u,v,L;\lambda_1],\quad f_v:=v_t+\mathcal{N}_v[u,v,L;\lambda_1],\quad f_L:=L_t+\mathcal{N}'_L[u,v;\lambda_2].\\
\end{split}
\end{align}

\subsection{Loss function}

We attempt to find the optimized parameters, including the weights, biases and additional coefficients in the activation, to minimize two loss functions $\mathscr{L}(\bar{\Theta})$ without parameter regularization as well as $\widetilde{\mathscr{L}}(\bar{\Theta})$ with weights parameter regularization, which are defined as the following forms respectively
\begin{align}\label{E7}
\begin{split}
&\mathscr{L}(\bar{\Theta})=Loss=Loss_{S}+Loss_{L}+Loss_{f_S}+Loss_{f_L}+Loss_a,\\
&\widetilde{\mathscr{L}}(\bar{\Theta})=Loss_{PR}=\mathscr{L}(\bar{\Theta})+\alpha \Omega(\bar{\Theta}),
\end{split}
\end{align}
in which $Loss_{S}, Loss_{L}, Loss_{f_S}$ and $Loss_{f_L}$ are defined as following
\begin{align}\label{E8}
\begin{split}
&Loss_{S}=\frac{1}{N_q}\left[\sum^{N_q}_{j=1}\big|\hat{u}(x^j,t^j)-u^j\big|^2+\sum^{N_q}_{j=1}\big|\hat{v}(x^j,t^j)-v^j\big|^2\right],\\
&Loss_{L}=\frac{1}{N_q}\sum^{N_q}_{j=1}\big|\hat{L}(x^j,t^j)-L^j\big|^2,
\end{split}
\end{align}
and
\begin{align}\label{E9}
\begin{split}
&Loss_{f_S}=\frac{1}{N_f}\left[\sum^{N_f}_{l=1}\big|f_u(x_f^l,t_f^l)\big|^2+\sum^{N_f}_{l=1}\big|f_v(x_f^l,t_f^l)\big|^2\right],\\
&Loss_{f_L}=\frac{1}{N_f}\sum^{N_f}_{l=1}\big|f_L(x_f^l,t_f^l)\big|^2,
\end{split}
\end{align}
where $\{x^j,t^j,u^j,v^j,L^j\}^{N_q}_{j=1}$ represent the initial and boundary value input data on \eqref{E6}. Here $\hat{u}(x^j,t^j), \hat{v}(x^j,t^j)$ and $\hat{L}(x^j,t^j)$ indicate the optimal training output data through the NN. Likewise, the collocation points on networks $f_u(x,t)$, $f_v(x,t)$ and $f_L(x,t)$ are denoted via $\{x_f^l,t_f^l\}^{N_{f}}_{l=1}$. Then the slope recovery term $Loss_a$ in the $\mathscr{L}(\bar{\Theta})$ is defined as
\begin{align}\label{E10}
Loss_a=\frac{1}{\frac{N_a}{D-1}\sum\limits_{d=1}^{D-1}\mathrm{exp}\Bigg(\frac{\sum\limits_{i=1}^{N_d}a_i^d}{N_d}\Bigg)},
\end{align}
where $1/N_a$ is the hyperparameter for slope recovery term $Loss_a$, thus we uniformly take $N_a=100$ to dominate the loss function $\widetilde{\mathscr{L}}(\bar{\Theta})$ and ensure that the final loss value is not too large in this paper. Here, the rapidly increasing activation slope value is imposed on the NN through the $Loss_a$ term, it ensures the non-vanishing gradient of the loss function and speeds up NN's traning process. Therefore, $Loss_{S}$ and $Loss_{L}$ correspond to the loss functions on the initial and boundary data, the $Loss_{f_S}$ and $Loss_{f_L}$ penalize the YO system not being satisfied on the collocation points, as well as the $Loss_a$ improves the convergence speed and promotes network optimization ability by changing the topology of $\mathscr{L}(\bar{\Theta})$.

\subsection{Parameter regularization}

In fact, $\alpha$ is a hyperparameter that weights the relative contribution of the norm penalty term $\Omega$ and loss function $\mathscr{L}(\bar{\Theta})$ in Eq. \eqref{E7}, and the $L^2$ parameter norm penalty $\Omega(\bar{\Theta})$ can be defined as following
\begin{align}\label{E11}
\Omega(\bar{\Theta})=\frac12\|\textbf{W}\|^2_2,
\end{align}
which drives the weights closer to the origin. Due to the biases typically require less data to fit accurately than the weights, we note that for NNs, we typically choose to use a parameter norm penalty $\Omega$ that penalizes only the weights of the affine transformation at each layer and leaves the biases unregularized. It is worth mentioning that the aforementioned slope recovery term $Loss_a$ can be regarded as a self-defined parameter regularization strategy for additional parameters $\{a^d_i\}$.

Moreover, one can gain some insight into the behavior of weight decay regularization by studying the gradient of the regularized objective function. Such a NN model has the
following total objective function:
\begin{align}\label{E-PR-1}
&\widetilde{\mathscr{L}}(\bar{\Theta})=\mathscr{L}(\bar{\Theta})+\frac{\alpha}{2}\textbf{W}^{\rm T}\textbf{W},
\end{align}
with the corresponding parameter gradient
\begin{align}\label{E-PR-2}
\nabla_{\textbf{W}}\widetilde{\mathscr{L}}(\bar{\Theta})=\nabla_{\textbf{W}}\mathscr{L}(\bar{\Theta})+\alpha\textbf{W}.
\end{align}
In order to take a single gradient step to update the weights, we perform following update:
\begin{align}\nonumber
\textbf{W}\leftarrow\textbf{W}-\eta(\nabla_{\textbf{W}}\mathscr{L}(\bar{\Theta})+\alpha\textbf{W}),
\end{align}
that is
\begin{align}\label{E-PR-3}
\textbf{W}\leftarrow(1-\eta\alpha)\textbf{W}-\eta\nabla_{\textbf{W}}\mathscr{L}(\bar{\Theta}),
\end{align}
where $\eta$ is learning rate, a positive scalar determining the size of the step. In the case that happens in a single step, one can see that addding the weight decay term has modified the learning rule to multiplicatively shrink the weight vector by a constant factor on each step, just before performing the usual gradient update. Next, we further analyze the influence of weight adecay in the whole training process.

In order to further study the impact of weight decay, we assume that $\textbf{W}^*$ is the weight vector which the objective function $\mathscr{L}(\bar{\Theta})$ without regularization obtain the minimal training cost, namely $\textbf{W}^*=\mathrm{arg min}_\textbf{W}\mathscr{L}(\bar{\Theta})$, and make a quadratic approximation of the objective function in the field of $\textbf{W}^*$. The approximation $\widehat{\mathscr{L}}$ is as follows
\begin{align}\label{E-PR-4}
\widehat{\mathscr{L}}(\bar{\Theta})=\mathscr{L}(\textbf{W}^*)+\frac12(\textbf{W}-\textbf{W}^*)^{\rm T}H(\textbf{W}-\textbf{W}^*),
\end{align}
where $H$ is the Hessian matrix of $\mathscr{L}$ with respect to $\textbf{W}$ evaluated at $\textbf{W}^*$. Since $\textbf{W}^*$ is defined to be a minimum of $\mathscr{L}$, the first-order term (the gradient) vanishes in this quadratic approximation, and one can also indicate that $H$ is positive semidefinite. When $\mathscr{L}$ gets the minimum, its gradient
\begin{align}\label{E-PR-5}
\nabla_{\textbf{W}}\widehat{\mathscr{L}}(\bar{\Theta})=H(\textbf{W}-\textbf{W}^*)=0.
\end{align}

Now we consider the effect of weight decay, and add the weight decay gradient in Eq. \eqref{E-PR-5}. Then we utilize the variable $\widehat{\textbf{W}}$ to represent the location of the minimum, and solve for the minimum of the regularized version of $\widehat{\mathscr{L}}$, we have
\begin{align}\label{E-PR-6}
\begin{split}
&H(\widehat{\textbf{W}}-\textbf{W}^*)+\alpha\widehat{\textbf{W}}=0,\\
&(H+\alpha I)\widehat{\textbf{W}}=H\textbf{W}^*,\\
&\widehat{\textbf{W}}=(H+\alpha I)^{-1}H\textbf{W}^*.
\end{split}
\end{align}
Obviously, once $\alpha$ approaches 0, the regularized solution $\widehat{\textbf{W}}$ approaches $\textbf{W}^*$. To further understand what happens when $\alpha$ grows, we decompose the real symmetric matrix $H$ into a diagonal matrix $\Lambda$ and a set of standard orthogonal basis $Q$ of eigenvectors, such that $H=Q\Lambda Q^{\rm T}$, then substitute the decomposition into \eqref{E-PR-6}, and obtain
\begin{align}\label{E-PR-7}
\widehat{\textbf{W}}=Q(\Lambda+\alpha I)^{-1}\Lambda Q^{\rm T}\textbf{W}^*.
\end{align}
From matrix $(\Lambda+\alpha I)^{-1}\Lambda$, one can see that the effect of weight decay is to rescale $\textbf{W}^*$ along the axes defined by the eigenvectors of $H$. Specifically, the component of $\textbf{W}^*$ that is aligned with the $i$th eigenvector of $H$ is rescaled by a factor of $\frac{\lambda_i}{\lambda_i+\alpha}$. The effect of regularization along the direction with large $H$ eigenvalues (such as $\lambda_i\gg\alpha$) is relatively small. Whereas, the components with $\lambda_i\ll\alpha$ will be shrunk to almost zero magnitude.

In other words, only the parameters in the direction of significantly reducing the objective function are preserved relatively intact. In directions that do not contribute to reducing the objective function, a small eigenvalue of the Hessian matrix tells us that movement in this direction will not significantly increase the gradient. Components of the weight vector corresponding to such unimportant directions are decayed away through the use of the regularization throughout training. This also requires that when we add regularization to NN, the weight coefficient $\alpha$ should not be too small or too large in practical application, which will also be reflected in the later training results.

\subsection{Optimization algorithm and improved PINN}

The resulting optimization problem leads to finding the minimum value of the loss function by optimizing the parameters $\bar{\Theta}$, that is, we seek
\begin{align}\nonumber
\bar{\Theta}^*=\mathop{\mathrm{arg\,min}}\limits_{\bar{\Theta}\in\bar{\mathcal{P}}}\widetilde{\mathscr{L}}(\bar{\Theta}).
\end{align}
Generally, one can approximate the solutions to this minimization problem iteratively by one of the forms of gradient descent algorithm. The stochastic gradient descent (SGD) and its variants are probably the most used optimization algorithms for machine learning in general and for deep learning in particular \cite{Ruder2017}. In this work, we introduce Adam optimizer and L-BFGS optimizer to optimize the loss function. Specifically, we employ the Adam optimizer, which is a variant of the SGD algorithm, as well as the L-BFGS optimizer, which is a full-batch gradient descent optimization algorithm based on a quasi-Newton method to optimize the loss function \cite{Kingma2014,Liu1989}. Moreover, in order to better measure the error from training results, we introduce $\mathbb{L}_2$ norm error, which is defined as follows
\begin{align}\nonumber
\mathrm{Error}=\frac{\sqrt{\sum\limits_{k=1}^{N}\big|q^{\mathrm{exact}}(\textbf{x}_{k})-q^{\mathrm{predict}}(\textbf{x}_{k};\bar{\Theta})\big|^2}}{\sqrt{\sum\limits_{k=1}^{N}\big|q^{\mathrm{exact}}(\textbf{x}_{k})\big|^2}},
\end{align}
where $q^{\mathrm{predict}}(\textbf{x}_{k};\bar{\Theta})$ and $q^{\mathrm{exact}}(\textbf{x}_{k})$ represent the model training prediction solution and exact analytical solution at point $\textbf{x}_{k}=(x_k,t_k)$, respectively.

In order to master the improved PINN approach more systematically, the improved PINN algorithm flow chart tailored for the YO system is exhibited in following Fig. \ref{F1}, in which vividly shows the NN along with the supplementary physics-informed part, and the loss function arising from the NN part as well as the residual part from the governing equation given via the physics-informed part. Thus, one can minimizing the loss function below certain tolerance $\varepsilon$ until a prescribed maximum number of iterations by seeking the optimal values of weights $\textbf{W}$, biases $\textbf{b}$ and scalable parameter $a^d_i$. Since the YO system contains two components $S(x,t)$ and $L(x,t)$, one can see that the ``NN" part has three output functions $\{u,v,L\}$, and there are three nonlinear equation constraints in the ``PDE" part in Fig. \ref{F1}. Once increasing the number of coupling system components, the number of output functions and nonlinear equation constraints of the improved PINN will multiply increase. Moreover, the corresponding procedure steps for the improved PINN with adaptive activation function and $L^2$ norm parametric regularization is displayed in the following Tab. \ref{Tab:IPINN}.

\begin{figure}[htbp]
\centering
\begin{minipage}[t]{0.99\textwidth}
\centering
\includegraphics[height=9cm,width=15cm]{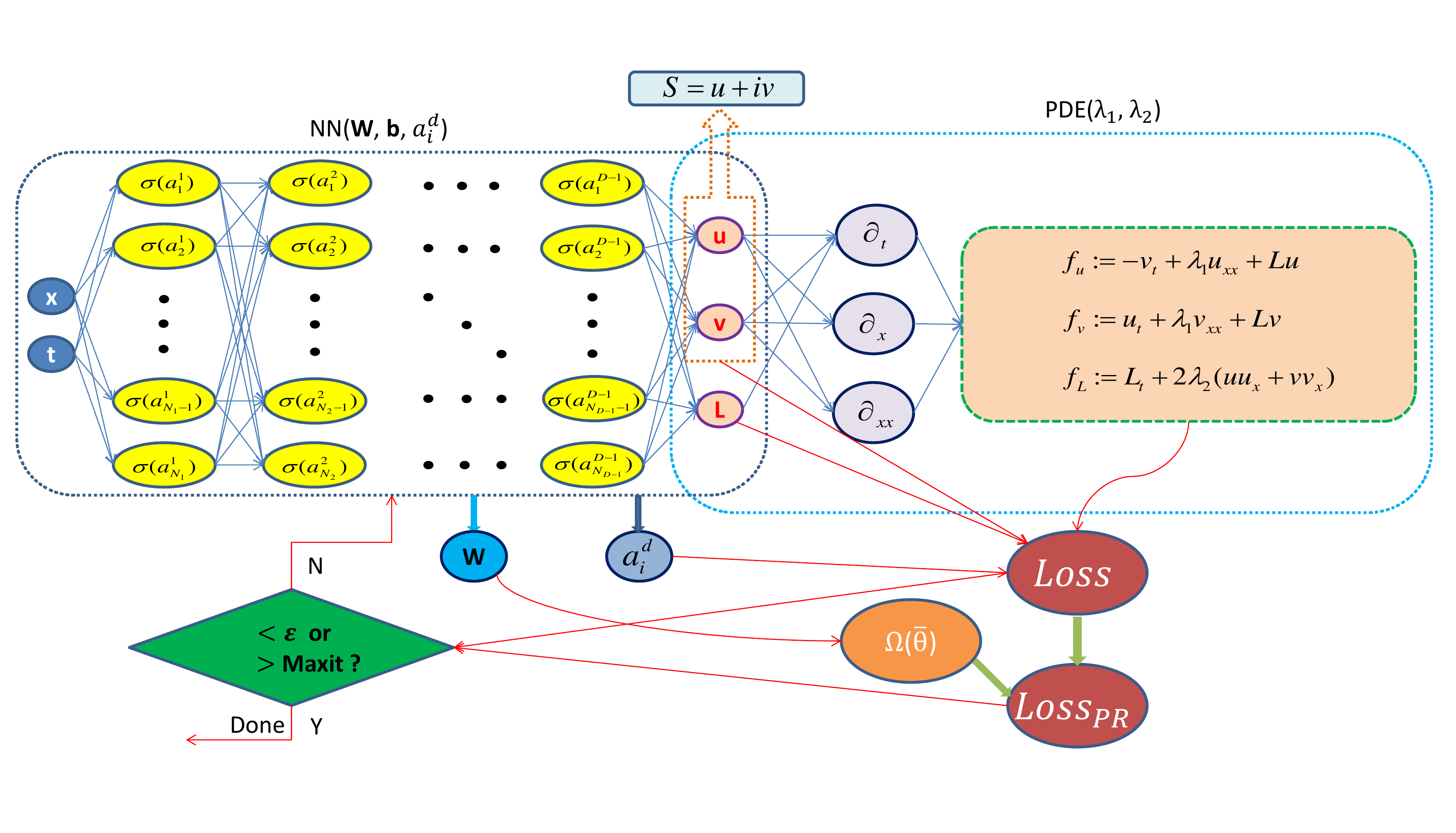}
\end{minipage}
\centering
\caption{(Color online) Schematic of improved PINN with adaptive activation function and $L^2$ norm parameter regularization for the YO system. The left NN is the universal approximation network without informed while the right NN is physics-informed network dominated by the governing equation, both NNs share hyper-parameters that contribute to the loss function.}
\label{F1}
\end{figure}

\begin{table}[htbp]
  \caption{Improved PINN algorithm with adaptive activation function and $L^2$ norm parameter regularization for the YO system.}
  \label{Tab:IPINN}
  \centering
  \begin{tabular}{p{15cm}}
  \toprule[2pt]
  \quad \textbf{Step 1}: Specification of training set in computational domain:\\
  \quad \emph{Training data}: $\{x^j,t^j,u^j,v^j,L^j\}^{N_q}_{j=1}$, \emph{Residual training points}: $\{x_f^l,t_f^l\}^{N_{f}}_{l=1}.$\\
  \quad \textbf{Step 2}: Construct NN $\textbf{q}(\textbf{x};\bar{\Theta})$ with random initialization of parameters $\bar{\Theta}$.\\
  \quad \textbf{Step 3}: Construct the residual NN $\{f_u,f_v,f_L\}$ by substituting surrogate $\textbf{q}(\textbf{x};\bar{\Theta})$ into the governing equations using automatic differentiation and other arithmetic operations.\\
  \quad \textbf{Step 4}: Specification of the loss function $\widetilde{\mathscr{L}}(\bar{\Theta})$ that includes the slope recovery term and parameter regularization term.\\
  \quad \textbf{Step 5}: Find the best parameters $\bar{\Theta}^*$ using a suitable optimization method for minimizing the loss function $\widetilde{\mathscr{L}}(\bar{\Theta})$ as\\
  \qquad\qquad\qquad\qquad\qquad\qquad\qquad\quad $\bar{\Theta}^*=\mathop{\mathrm{arg\,min}}\limits_{\bar{\Theta}\in\bar{\mathcal{P}}}\mathscr{L}(\bar{\Theta})$.\\
  \bottomrule[2pt]
  \end{tabular}
\end{table}

\subsection{Training data and network environment}

In supervised learning, training data is important to train the NN, which can be obtained from the exact solution (if available) or from high-resolution numerical solution using numerical methods like spectral method, finite element method, Chebfun numerical method, discontinuous Galerkin method etc, as per the problem at hand. Furthermore, training data can also be obtained from carefully performed physical experiments that may produce high- and low-fidelity data sets. Fortunately, the nonlinear integrable systems like YO system has very good integrability, and there are some effective methods to solve a variety of accurate localized wave solutions, which provide rich sample space for the extraction of training data. Here, we use the traditional finite difference scheme on the even grid to discretize these exact solutions to obtain the training data in Matlab.

In the adaptive activation function, the initialization of scalable parameters are carried out in the case of $n=10,a_i^d=0.1$, namely $na_i^d=1$. In addition, we select relatively simple multi-layer perceptrons (i.e., feedforward NNs) with the Xavier initialization and the hyperbolic tangent ($\tanh$) as activation function. The residual training points $\{x_f^l,t_f^l\}^{N_{f}}_{l=1}$ are generated by the Latin Hypercube Sampling method (LHS) \cite{Stein1987}. All the codes in this article is based on Python 3.7 and Tensorflow 1.15, and all numerical experiments reported here are run on a DELL Precision 7920 Tower computer with 2.10 GHz 8-core Xeon Silver 4110 processor, 64 GB memory and 11 GB Nvidia GeForce GTX 1080 Ti video card.

\section{The forward problem of the YO system}

In this section, we focus on the forward problem of the YO system, that is reveal the data-driven RWs for the YO system by means of small data set and 9 hidden layers deep PINN with 40 neurons per layer. Specifically, after the unknown parameters of YO system \eqref{E1} are determined, that is $\lambda_1=0.5$ and $\lambda_2=1$, once some initial boundary value data points are given, one can successfully approximate the various RWs through improved PINN with parameter regularization of hyper-parameter $\alpha=0.0001$. Here the physics-informed parts Eq. \eqref{E6} of the improved PINN for YO sysem \eqref{E1} become the following formula
\begin{align}\label{TE1}
\begin{split}
&f_u:=-v_t+0.5u_{xx}+uL,\quad f_v:=u_t+0.5v_{xx}+vL,\quad f_L:=L_t-(2uu_x+2vv_x).
\end{split}
\end{align}

The exact form of these RWs for the YO system \eqref{E1} with $\lambda_1=0.5$ and $\lambda_2=1$ have been derived by the Darboux transformation \cite{ChenSH2014}, the general vector form of RW can be expressed as
\begin{align}\label{E12}
\begin{split}
&S(x,t)=a\mathrm{e}^{\mathrm{i}kx-\mathrm{i}\big(\frac12k^2-b\big)t}\bigg[1-\frac{\mathrm{i}t+(\mathrm{i}x)/(2m-k)+1/\big(2(2m-k)(m-k)\big)}{(x-mt)^2+n^2t^2+1/(4n^2)}\bigg],\\
&L(x,t)=b+2\frac{n^2t^2-(x-mt)^2+1/(4n^2)}{[(x-mt)^2+n^2t^2+1/(4n^2)]^2},
\end{split}
\end{align}
and $m$, $n$, $a$, $b$ and $k$ satisfy the following relationship
\begin{align}\label{E13}
\begin{split}
&m=\frac16\big[5k-\sqrt{3(k^2+\eta+\sigma/\eta)}\big],\,n=\pm\sqrt{(3m-k)(m-k)},\\
&\sigma=\frac19k^4+6a^2k,\,\rho=\frac12k^6-\frac{1}{54}(27a^2+5k^3)^2,\\
&\eta=\begin{cases}
-\big(\rho-\sqrt{\rho^2-\sigma^3}\big)^{1/3}& k\leqslant-3k_n,\\
\big(-\rho+\sqrt{\rho^2-\sigma^3}\big)^{1/3}& -3k_n<k\leqslant\frac32k_n,
\end{cases}
\end{split}
\end{align}
where $k_n=(2a^2)^{1/3}$, $a>0$, $b\geqslant0$ and $k\in\mathbb{R}$. It is clear that the complex short-wave RW $|S|$ is characterized by the second-order polynomial of $x$ and $t$, while the real long-wave RW $L$ involves the fourth-order polynomial. One can observe that although having a form similar to that of Peregrine soliton \cite{Peregrine1983}, the RWs in \eqref{E12} admit more complex dynamics than the latter.

According to the central amplitude and the relative position of two zero-amplitude points, one can divide the regime of the short-wave RW $S(x,t)$ into three regions, that is bright RW region, intermediate RW region and dark RW region, which are corresponding to parametric conditions $k\leqslant0$, $0<k<(4/3)^{1/3}k_n$ (here $(4/3)^{1/3}\approx1.1$) and $(4/3)^{1/3}k_n\leqslant k<1.5k_n$, respectively \cite{ChenSH2014}. Specifically, as follows:

$\bullet$ \textbf{Bright-Bright RWs $S_{\rm brw}$ and $L_{\rm brw1}$}

In order to obtain the bright-bright RWs $S_{\rm brw}$ and $L_{\rm brw1}$ of YO system, one can lead to $m=-\frac12,\,n=\pm\frac{\sqrt{3}}{2}$ by taking $a=1,\, b=0,\, k=0$ in Eq. \eqref{E13} and combining the value range of $k$ from the aforementioned results. Substituting the above parameters into the formula Eq. \eqref{E12}, the specific form of bright-bright RWs $S_{\rm brw}$ and $L_{\rm brw1}$ for YO system are shown as follows
\begin{align}\label{E14}
\begin{split}
&S_{\rm brw}(x,t)=\frac{-3\mathrm{i}t+3\mathrm{i}x+3t^2+3tx+3x^2-2}{3t^2+3tx+3x^2+1},\\
&L_{\rm brw1}(x,t)=\frac{3(3t^2-6tx-6x^2+2)}{(3t^2+3tx+3x^2+1)^2}.
\end{split}
\end{align}
Apparently, the complex short-wave RW $S_{\rm brw}(x,t)$ takes plane $|S_{\rm brw}|=1$ as the background wave, while real long-wave RW $L_{\rm brw1}(x,t)$ takes plane $L_{\rm brw1}=0$ as the background wave. Furthermore, $|S_{\rm brw}(x,t)|$ obtains the maximum amplitude 2 at $(x,t)=(0,0)$, while $L_{\rm brw1}(x,t)$ obtains the maximum amplitude 6 at $(x,t)=(0,0)$.

$\bullet$ \textbf{Intermediate-Bright RWs $S_{\rm irw}$ and $L_{\rm brw2}$}

Similarly, in order to obtain the intermediate-bright RWs $S_{\rm irw}$ and $L_{\rm brw2}$ of YO system, one can also take $a=1,\,b=0$, but parameter $k$ have to meet the condition $0<k<(4/3)^{1/3}k_n$, here $k_n=2^{1/3}$. In particular, we derive the values of parameters $m$ and $n$ by taking $k=\frac122^{1/3}$, then substitute them into the Eq. \eqref{E12}, and then one can obtain the intermediate-bright RWs $S_{\rm irw}$ and $L_{\rm brw2}$ of YO system. Since the forms of expressions for $m$, $n$, $S_{\rm irw}$ and $L_{\rm brw2}$ are complex, we omit them here.

$\bullet$ \textbf{Dark-Bright RWs $S_{\rm drw}$ and $L_{\rm brw3}$}

In the same way, the dark-bright RWs $S_{\rm drw}$ and $L_{\rm brw3}$ of YO system are obtained by taking $a=1,\,b=0$ and making the parameter $k$ satisfy the condition $(4/3)^{1/3}k_n\leqslant k<1.5k_n$, here $k_n=2^{1/3}$. After taking $k=\frac652^{1/3}$, then we obtain the dark-bright RWs $S_{\rm drw}$ and $L_{\rm brw3}$ of YO system by means of known and derived parameters. Similarly, due to the complex form, we omit the specific expression form here.

Next, we employ aforementioned three exact RWs to obtain the initial boundary value condition data, so as to construct the training data set. Under different initial boundary value conditions, three different RWs dynamic structures are recovered by utilizing improved PINN with parameter regularization. In order to more intuitively exhibit the training effect of improved PINN with parameter regularization, we also compare it with the training results from PINN without parameter regularization.

\subsection{The data-driven bright-bright RWs}
In this section, in order to recover the data-driven bright-bright RWs of the YO system, we commit to introducing the initial boundary value conditions of the YO system to the 9-layer improved PINN with 40 neurons per layer. Selecting $[X_0,X_1]$ and $[T_0,T_1]$ in Eq. \eqref{E1} as $[-5.0,5.0]$ and $[-2.0,2.0]$ respectively, we can write the corresponding initial value conditions as follows
\begin{align}\label{E15}
\begin{split}
&S^0(x)=S_{\rm brw}(x,-2.0),\,L^0(x)=L_{\rm brw1}(x,-2.0),\,x\in[-5.0,5.0],
\end{split}
\end{align}
and the Dirichlet boundary conditions become
\begin{align}\label{E16}
\begin{split}
&S^{\mathrm{lb}}(t)=S_{\rm brw}(-5.0,t),\,L^{\mathrm{lb}}(t)=L_{\rm brw1}(-5.0,t),\,t\in[-2.0,2.0],\\
&S^{\mathrm{ub}}(t)=S_{\rm brw}(5.0,t),\,L^{\mathrm{ub}}(t)=L_{\rm brw1}(5.0,t),\,t\in[-2.0,2.0].
\end{split}
\end{align}

In order to obtain the original training data set of the above initial boundary value conditions \eqref{E15} and \eqref{E16}, we discretize the exact bright-bright RWs $S_{\rm brw}$ and $L_{\rm brw1}$ \eqref{E14} based on the finite difference method by dividing the spatial region $[-5.0,5.0]$ into 2000 points and the temporal region $[-2.0,2.0]$ into 1000 points in Matlab. Furthermore, in addition to the data set composed of the aforementioned initial boundary value conditions, the residual data set is used to calculate the $\mathbb{L}_2$ norm error by comparing with the data-driven bright-bright RWs. After that, a smaller training dataset that containing initial-boundary data will be generated by randomly extracting $N_q=1000$ initial boundary value data points from original dataset and $N_f=20000$ collocation points which are produced by the LHS. According to 20000 Adam iterations and 50000 L-BFGS iterations, the latent bright-bright RWs $S(x,t)$ and $L(x,t)$ have been successfully learned by employing the improved PINN with parameter regularization, and the network achieved relative $\mathbb{L}_2$ error of 4.968430$\rm e$-04 for the bright RW $S(x,t)$ and relative $\mathbb{L}_2$ error of 1.763312$\rm e$-03 for the bright RW $L(x,t)$, and the total number of iterations is 70000.

\begin{figure}[htbp]
\centering
\subfigure[]{
\begin{minipage}[t]{0.48\textwidth}
\centering
\includegraphics[height=3.5cm,width=7cm]{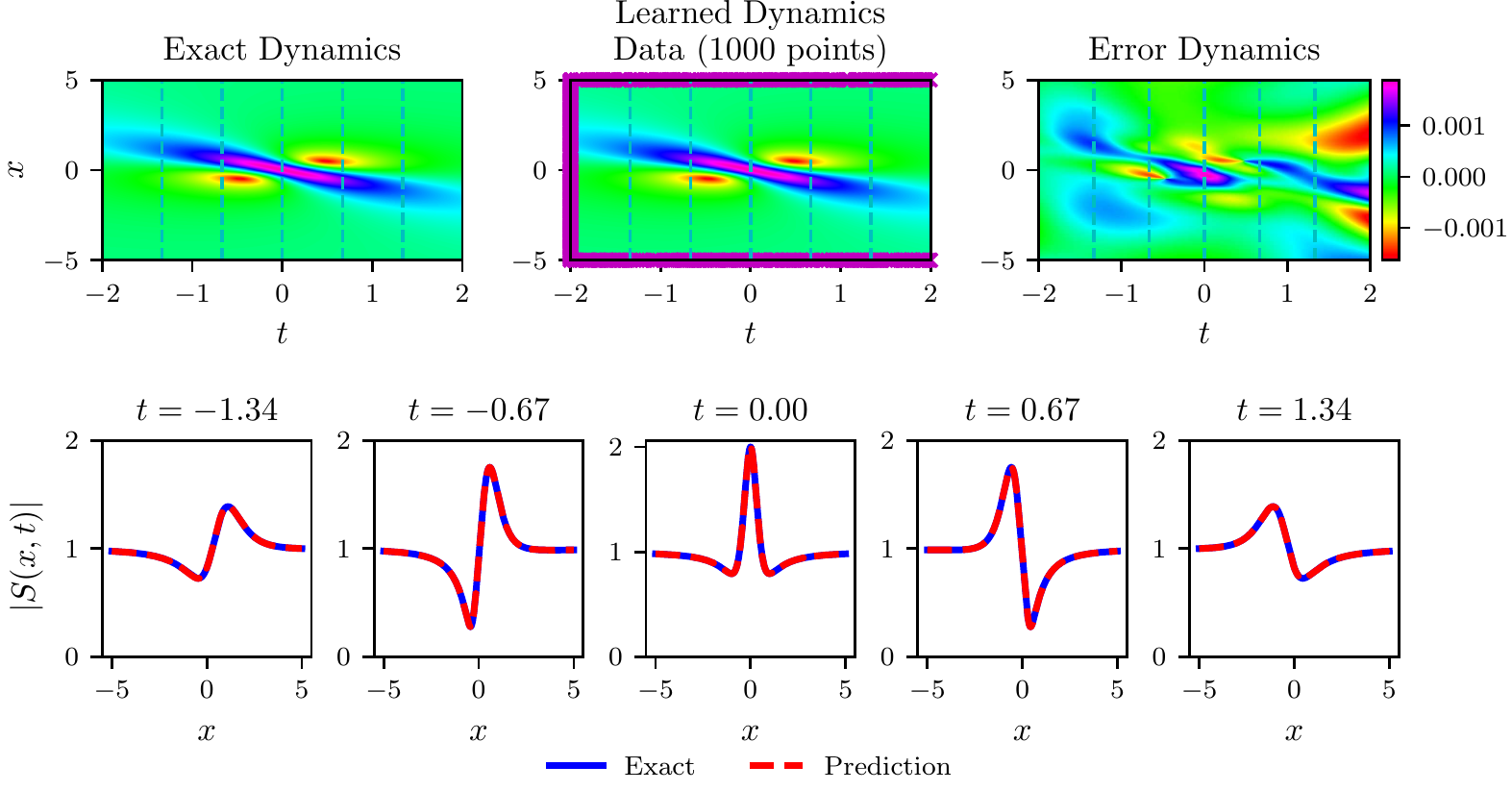}
\end{minipage}
}%
\subfigure[]{
\begin{minipage}[t]{0.48\textwidth}
\centering
\includegraphics[height=3.5cm,width=6.5cm]{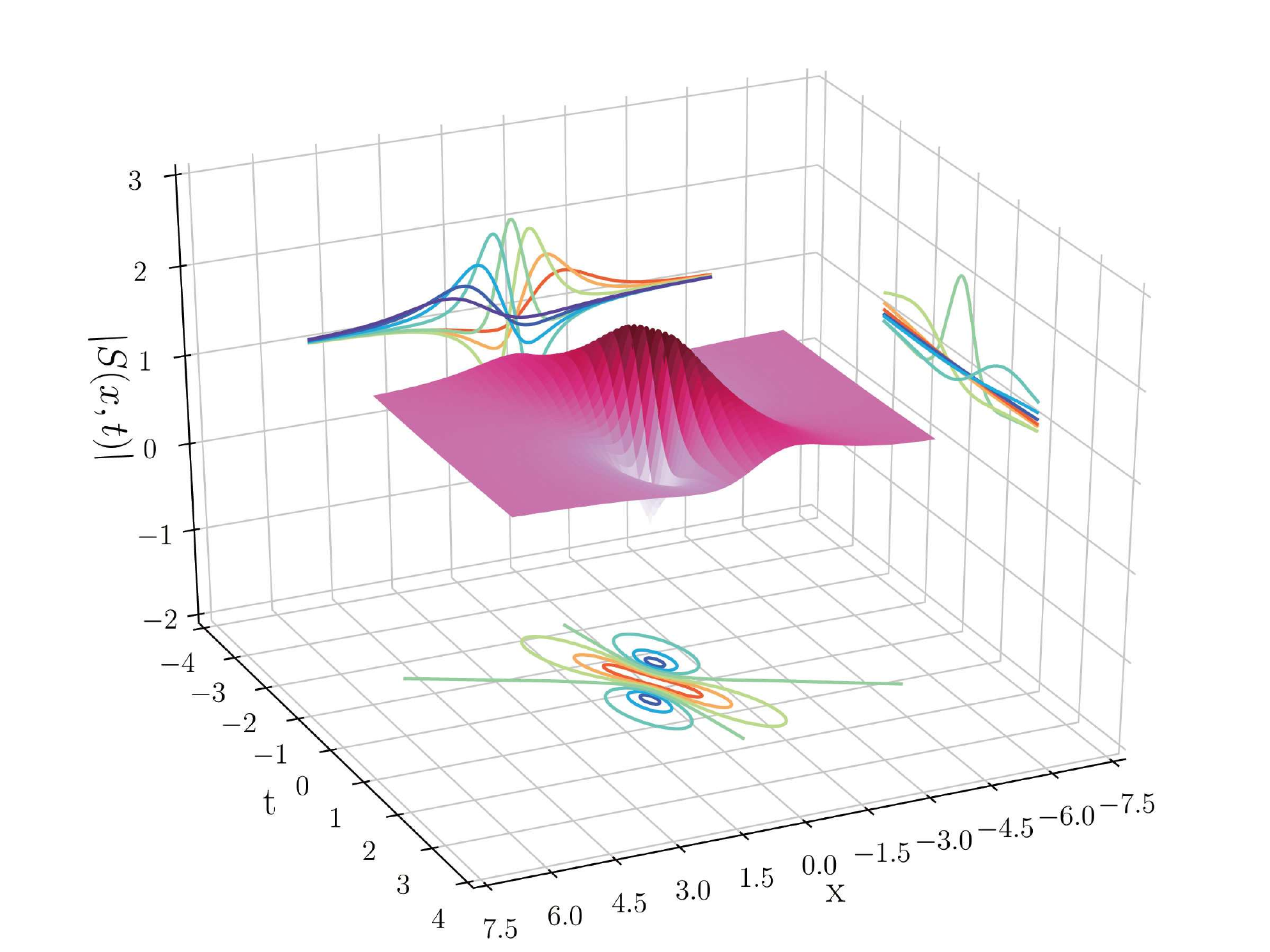}
\end{minipage}%
}\\%
\subfigure[]{
\begin{minipage}[t]{0.48\textwidth}
\centering
\includegraphics[height=3.5cm,width=7cm]{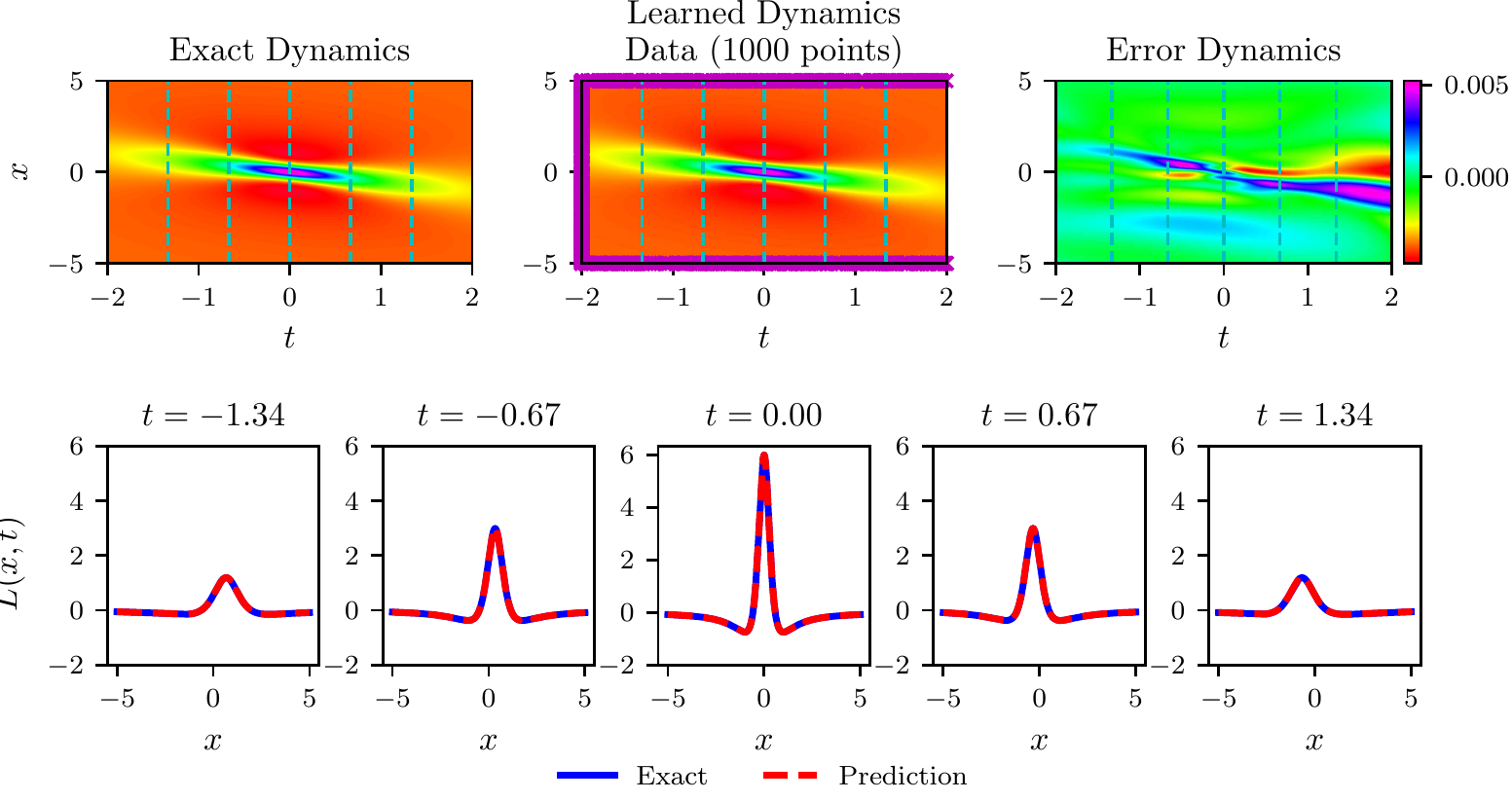}
\end{minipage}
}%
\subfigure[]{
\begin{minipage}[t]{0.48\textwidth}
\centering
\includegraphics[height=3.5cm,width=6.5cm]{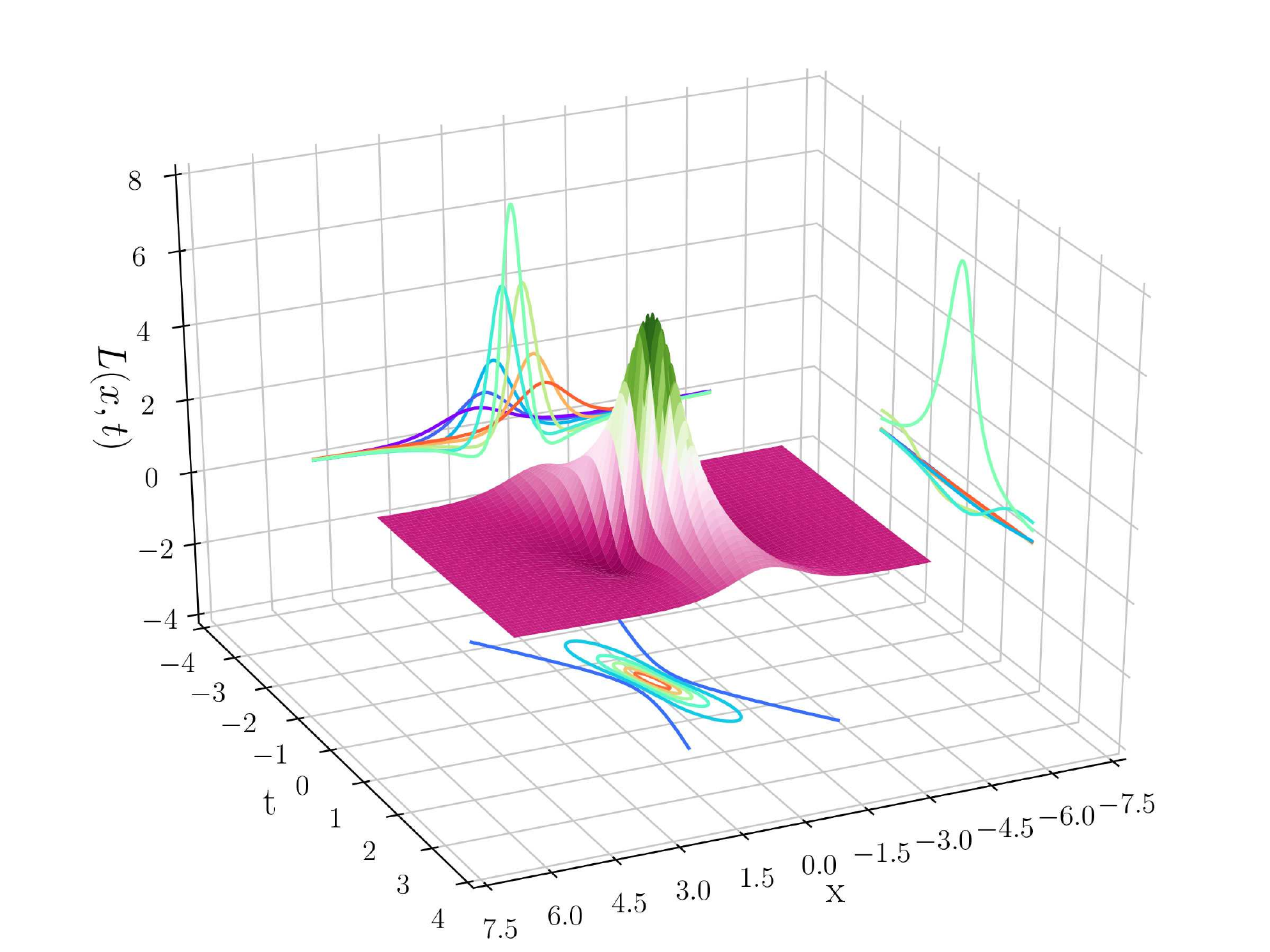}
\end{minipage}%
}%
\centering
\caption{(Color online) The data-driven bright-bright RWs $S(x,t)$ and $L(x,t)$ arising from the improved PINN with parameter regularization of hyper-parameter $\alpha=0.0001$ by randomly choosing $N_q=1000$ initial and boundary points which have been shown by means of mediumorchid $``\times"$ in learned dynamics , as well as $N_f = 20000$ collocation points in the corresponding spatiotemporal region: (a) and (c) The exact, learned and error dynamics density plots with five distinct tested times $t=-1.34, -0.67, 0.00, 0.67$ and 1.34 (darkturquoise dashed lines), as well as sectional drawings which contain the learned and explicit bright-bright RWs at the aforementioned five distinct times; (b) and (d) The 3D plot with contour map for the data-driven bright-bright RWs.}
\label{F2}
\end{figure}

\begin{figure}[htbp]
\centering
\subfigure[]{
\begin{minipage}[t]{0.48\textwidth}
\centering
\includegraphics[height=4cm,width=6cm]{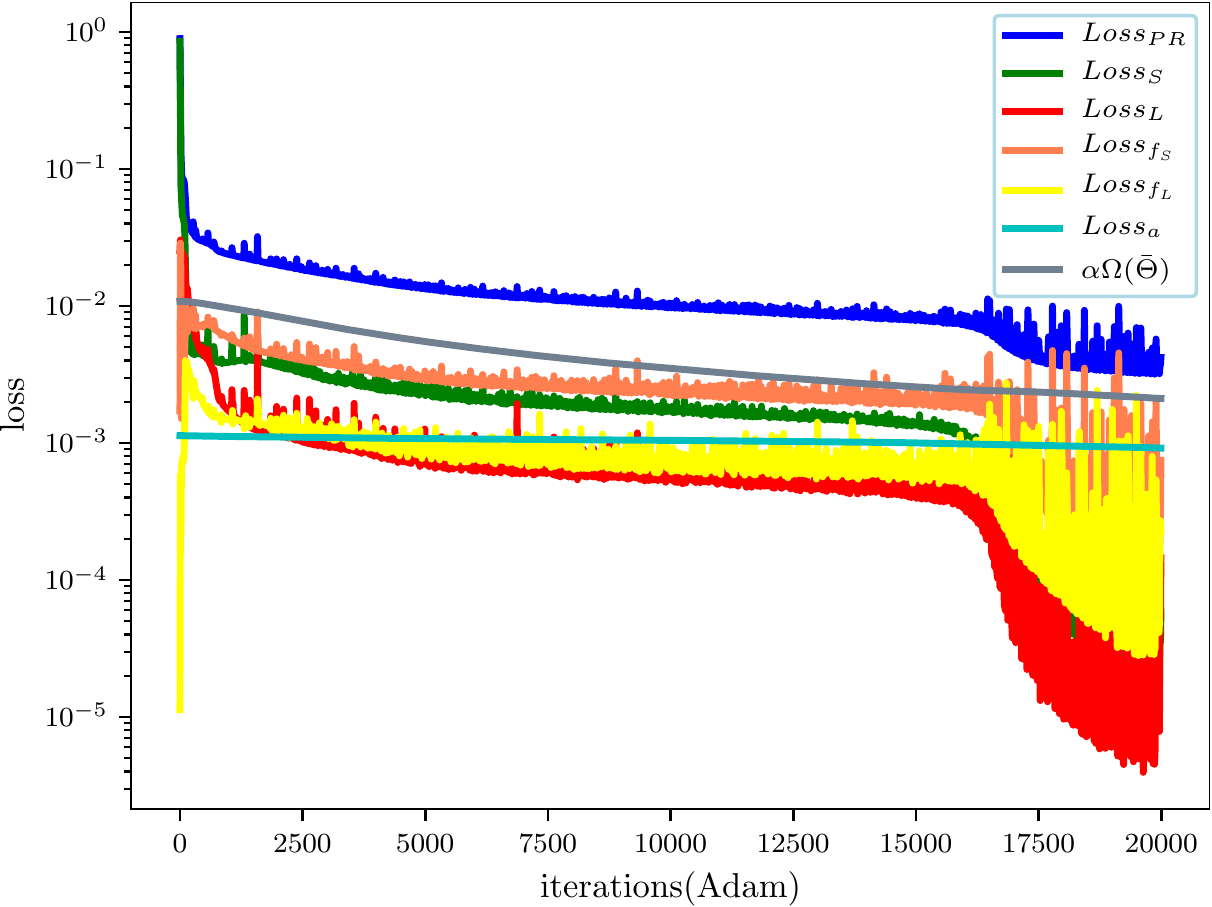}
\end{minipage}
}%
\subfigure[]{
\begin{minipage}[t]{0.48\textwidth}
\centering
\includegraphics[height=4cm,width=6cm]{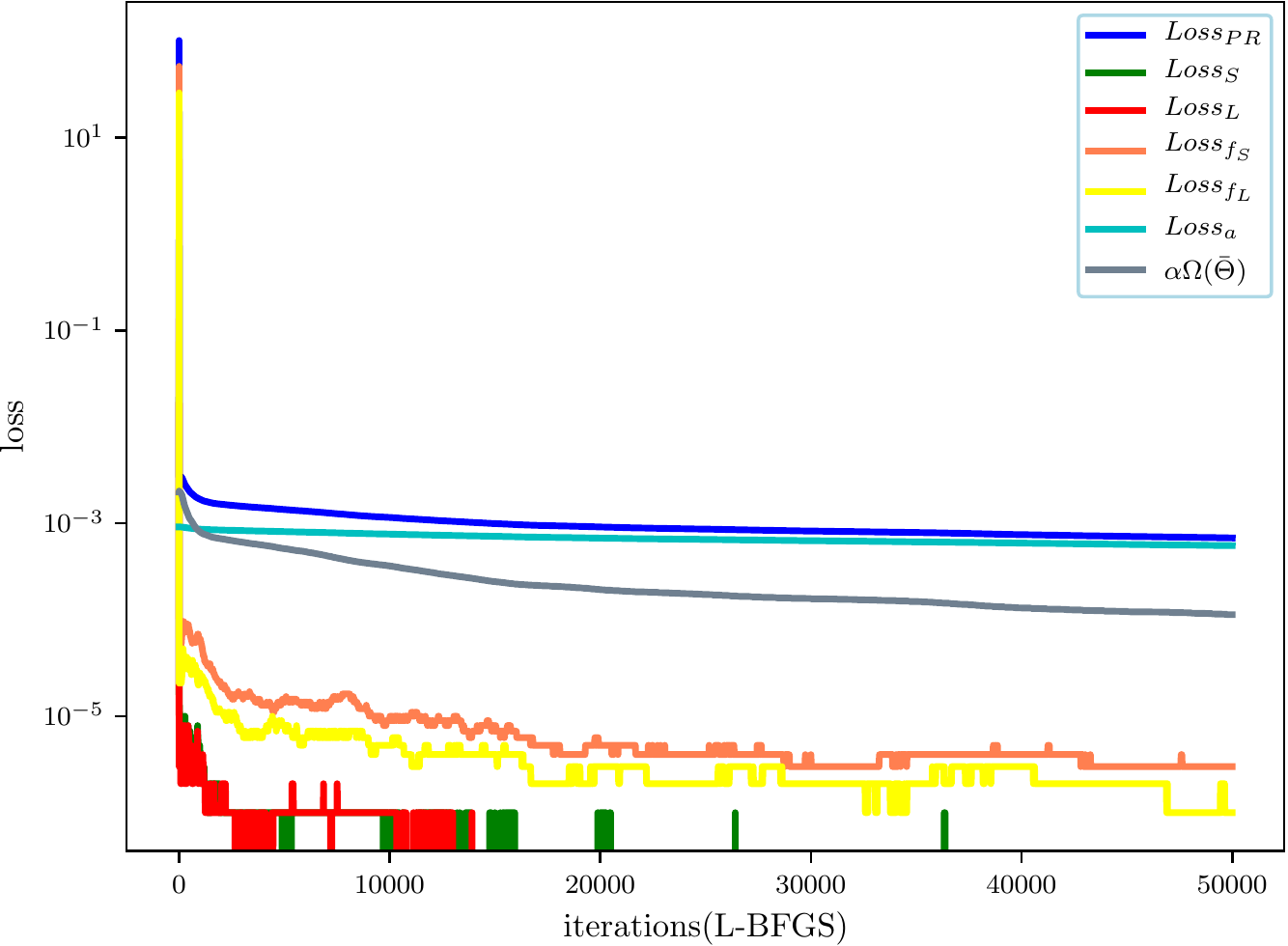}
\end{minipage}%
}%
\centering
\caption{(Color online) The loss function curve figures of the bright-bright RWs $S(x,t)$ and $L(x,t)$ arising from the improved PINN with the 20000 steps Adam and 50000 steps L-BFGS optimizations: (a) The loss function curve for the 20000 Adam optimization iterations; (b) The loss function curve for the 50000 L-BFGS optimization iterations.}
\label{F3}
\end{figure}

Figs. \ref{F2} - \ref{F3} present the deep learning results of the data-driven bright-bright RWs based on the improved PINN with parameter regularization for the YO system with the initial boundary value problem \eqref{E15} and \eqref{E16}. The left panels of Fig. \ref{F2} display the exact, learned and error dynamics with corresponding amplitude scale size on the right side for the bright RW $|S(x,t)|$ and bright RW $L(x,t)$, and exhibit the sectional drawings which contain the learned and explicit RWs at five different moments. From the density plots of learned dynamics and profiles which reveal amplitude and error of exact and prediction RWs in Fig. \ref{F2}, we observe that the amplitude of long-wave RW is much higher than that of short-wave RW. The right panels of Fig. \ref{F2} exhibit the 3D plots with contour map on three planes for the predicted bright-bright RWs. Fig. \ref{F3} showcases curve plots of the loss function after 20000 Adam optimization iterations and 50000 L-BFGS optimization iterations in improved PINN framework. In the left panel of Fig. \ref{F3}, one can see that the loss function curve $Loss_{PR}$ oscillatly descends in the process of optimizing the loss function by means of the Adam optimizer, and the gradient of the loss function descends very fast in the last about 4000 iterations. However, from the right panel of Fig. \ref{F3}, the loss function curve $Loss_{PR}$ linearly descends by means of the L-BFGS optimization algorithm. Furthermore, the loss function curves $Loss_S$ and $Loss_L$ of L-BFGS optimization are missing after a certain number of iterations, that is because the loss function value in this part is less than $\rm1e^{-6}$, which is beyond the statistical range of $``\%\mathrm{f}"$ in Python 3.7. During the process of optimizing the loss function $Loss_{PR}$ both in two optimization algorithms, $Loss_a$ and $\alpha\Omega(\bar{\Theta})$ descend linearly, which depend on their topological structure and mathematical form. In fact, $Loss_a$ ensures the better and faster convergence of the loss function by means of PINN algorithm with neuron-wise locally adaptive activation functions, while $\alpha\Omega(\bar{\Theta})$ guarantees the weight decay continuously by imposing the parameter regularization term.

\subsection{The data-driven intermediate-bright RWs}
In what follows, we will consider the initial-boundary value problem of the YO system for obtaining the data-driven intermediate-bright RWs by applying the multilayer improved PINN. Similarly, taking $[X_0,X_1]$ and $[T_0,T_1]$ in Eq. \eqref{E1} as $[-5.0,5.0]$ and $[-3.0,3.0]$ respectively, we derive the corresponding initial conditions $S^0(x)$ and $L^0(x)$, and Dirichlet boundary conditions as shown in the following formulas
\begin{align}\label{E17}
\begin{split}
&S^0(x)=S_{\rm irw}(x,-3.0),\,L^0(x)=L_{\rm brw2}(x,-3.0),\,x\in[-5.0,5.0],
\end{split}
\end{align}
and
\begin{align}\label{E18}
\begin{split}
&S^{\mathrm{lb}}(t)=S_{\rm irw}(-5.0,t),\,L^{\mathrm{lb}}(t)=L_{\rm brw2}(-5.0,t),\,t\in[-3.0,3.0],\\
&S^{\mathrm{ub}}(t)=S_{\rm irw}(5.0,t),\,L^{\mathrm{ub}}(t)=L_{\rm brw2}(5.0,t),\,t\in[-3.0,3.0].
\end{split}
\end{align}

By means of Matlab, we discretize the exact intermediate-bright RWs $S_{\rm irw}$ and $L_{\rm brw2}$ by utilizing the traditional finite difference scheme on even grids, and obtain the original training data set which contains initial data \eqref{E17} and boundary data \eqref{E18} by dividing the spatial region $[-5.0,5.0]$ into 2000 points and the temporal region $[-3.0,3.0]$ into 1000 points, as well as the remaining data will be used to obtain $\mathbb{L}_2$ norm error by comparing with predicted intermediate-bright RWs. Since the form of initial boundary value condition for intermediate-bright RWs is complex, it is necessary to increase the number of initial-boundary value training data points. Therefore, we generate a smaller training dataset where contains initial-boundary data by randomly extracting $N_q=2000$ initial-boundary value data points from original training dataset and $N_f=30000$ collocation points produced via LHS in the corresponding spatiotemporal region. Then, the intermediate-bright RWs $S(x,t)$ and $L(x,t)$ have been successfully learned by imposing a 9-hidden-layer improved PINN with 40 neurons per layer, and the related loss functions are optimized through 20000 Adam iterations and 50000 L-BFGS iterations. The relative $\mathbb{L}_2$ errors of the improved PINN model are 1.168852$\rm e$-03 for $S(x,t)$ and 6.766132$\rm e$-03 for $L(x,t)$, the total number of iterations is 70000.

\begin{figure}[htbp]
\centering
\subfigure[]{
\begin{minipage}[t]{0.48\textwidth}
\centering
\includegraphics[height=3.5cm,width=7cm]{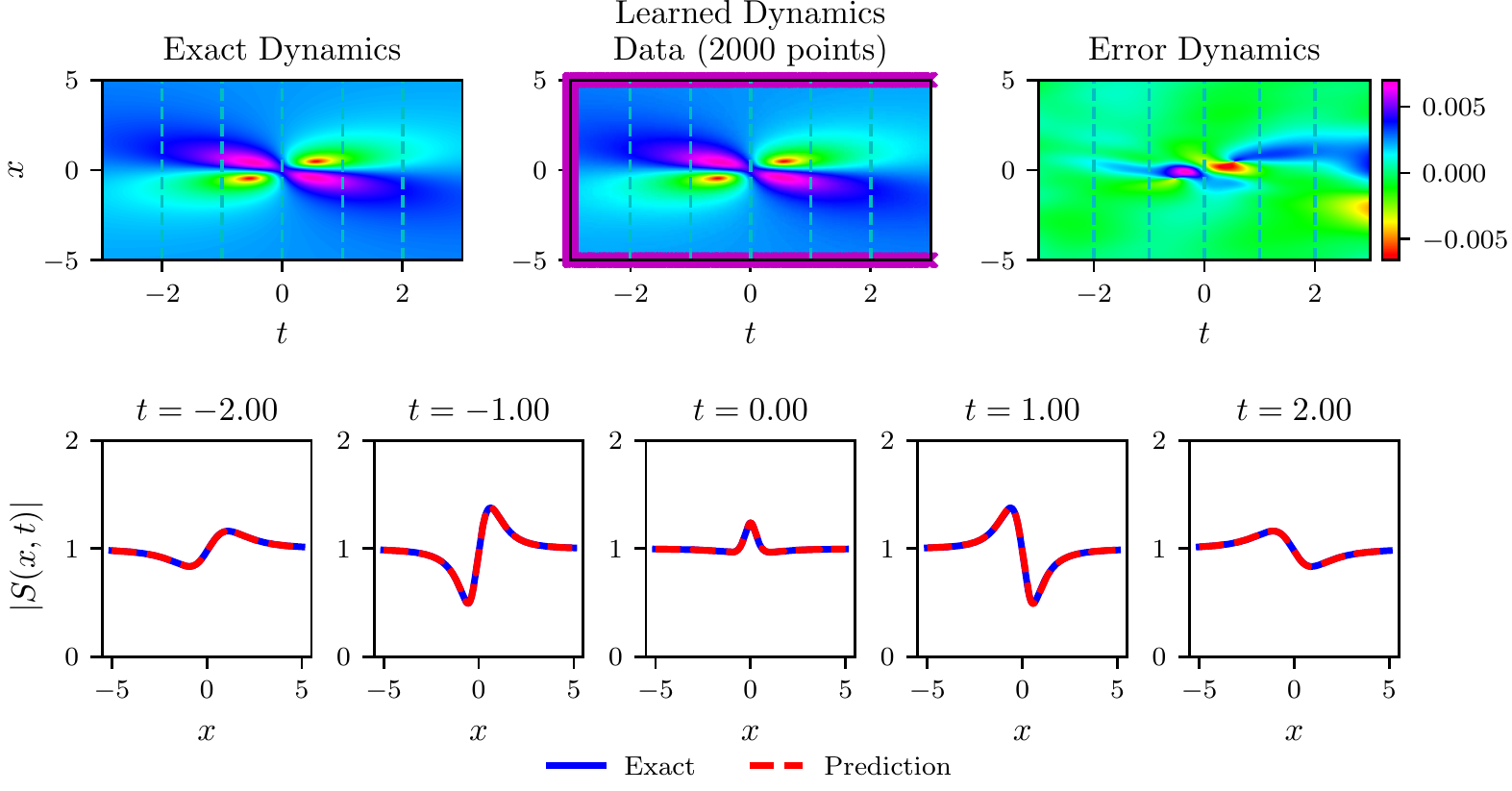}
\end{minipage}
}%
\subfigure[]{
\begin{minipage}[t]{0.48\textwidth}
\centering
\includegraphics[height=3.5cm,width=6.5cm]{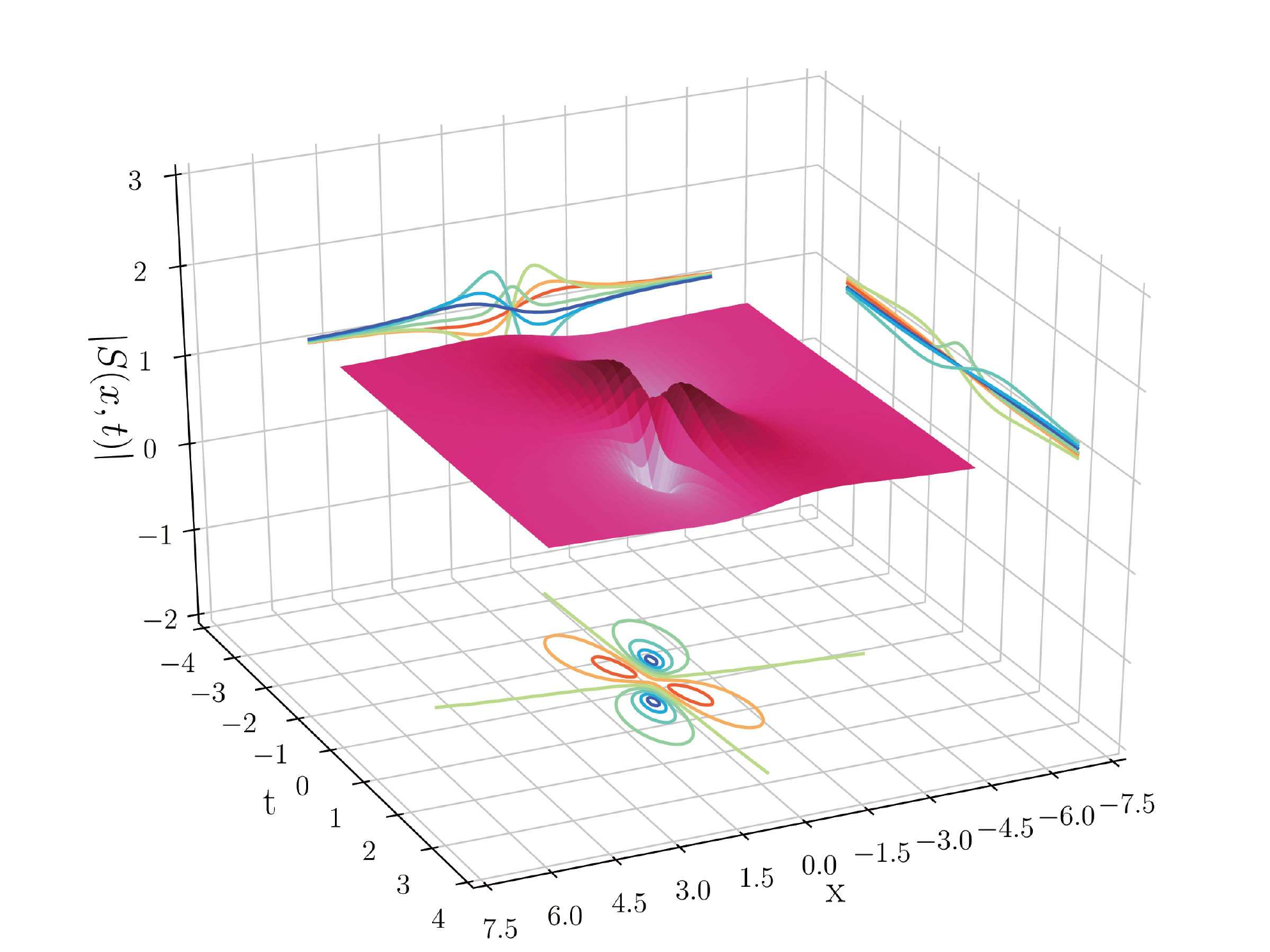}
\end{minipage}%
}\\%
\subfigure[]{
\begin{minipage}[t]{0.48\textwidth}
\centering
\includegraphics[height=3.5cm,width=7cm]{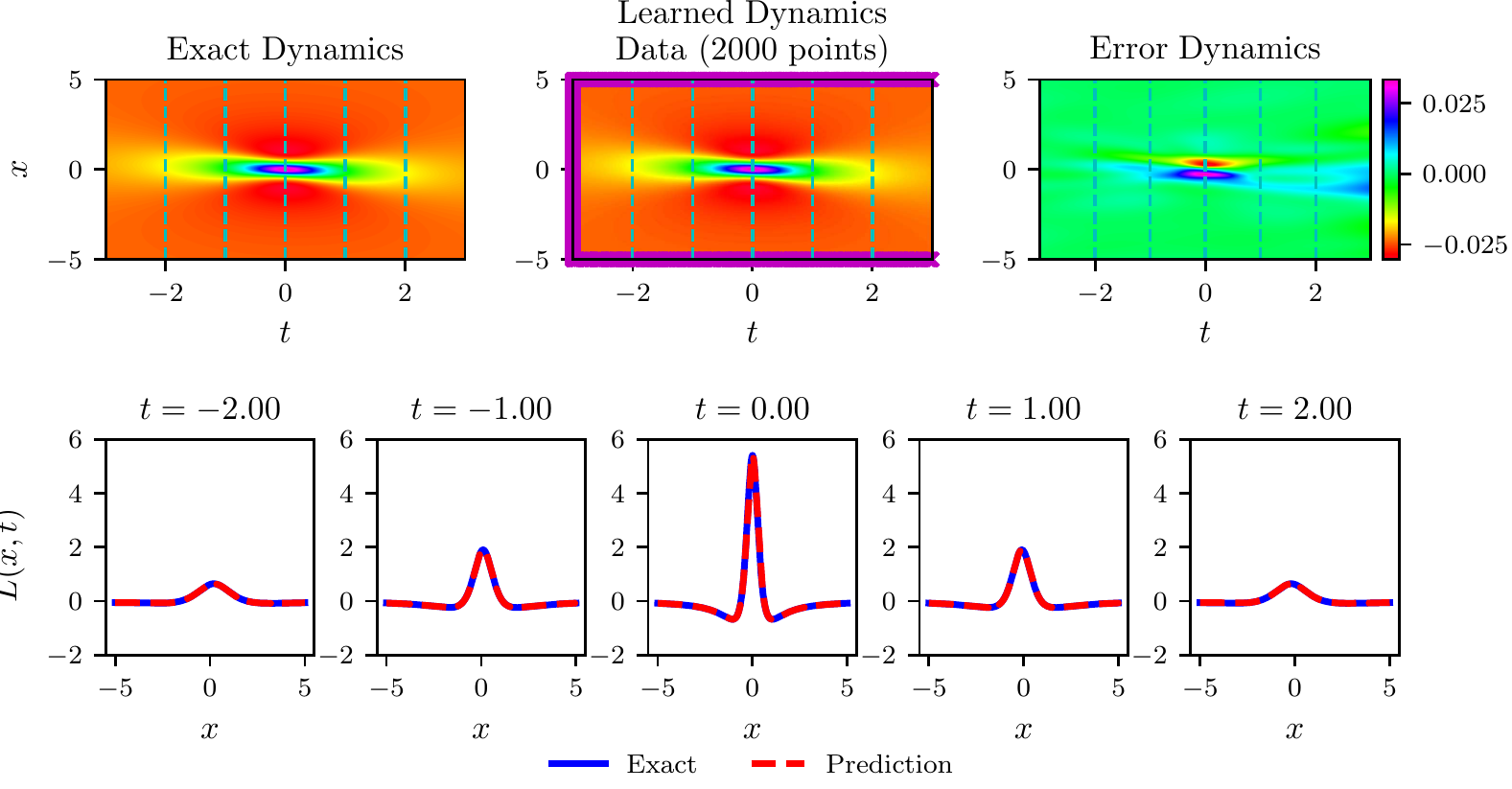}
\end{minipage}
}%
\subfigure[]{
\begin{minipage}[t]{0.48\textwidth}
\centering
\includegraphics[height=3.5cm,width=6.5cm]{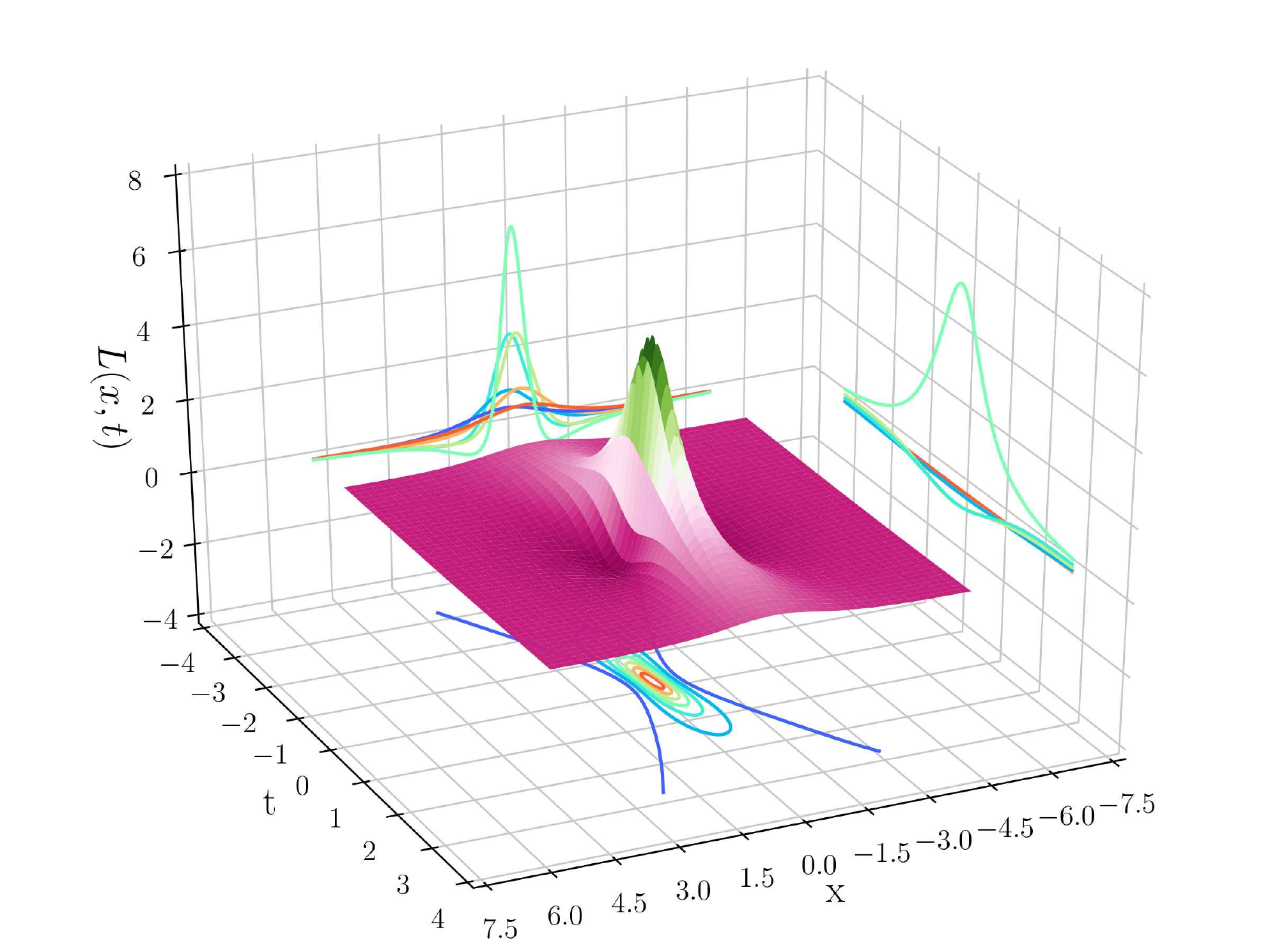}
\end{minipage}%
}%
\centering
\caption{(Color online) The data-driven intermediate-bright RWs $S(x,t)$ and $L(x,t)$ arising from the improved PINN with parameter regularization of hyper-parameter $\alpha=0.0001$ by randomly choosing $N_q=2000$ initial and boundary points which have been shown by using mediumorchid $``\times"$ in learned dynamics , as well as $N_f = 30000$ collocation points in the corresponding spatiotemporal region: (a) and (c) The exact, learned and error dynamics density plots with five distinct tested times $t=-2.00, -1.00, 0.00, 1.00$ and 2.00 (darkturquoise dashed lines), as well as sectional drawings which contain the learned and explicit intermediate-bright RWs at the aforementioned five distinct times; (b) and (d) The 3D plot with contour map for the data-driven intermediate-bright RWs.}
\label{F4}
\end{figure}

\begin{figure}[htbp]
\centering
\subfigure[]{
\begin{minipage}[t]{0.48\textwidth}
\centering
\includegraphics[height=4cm,width=6cm]{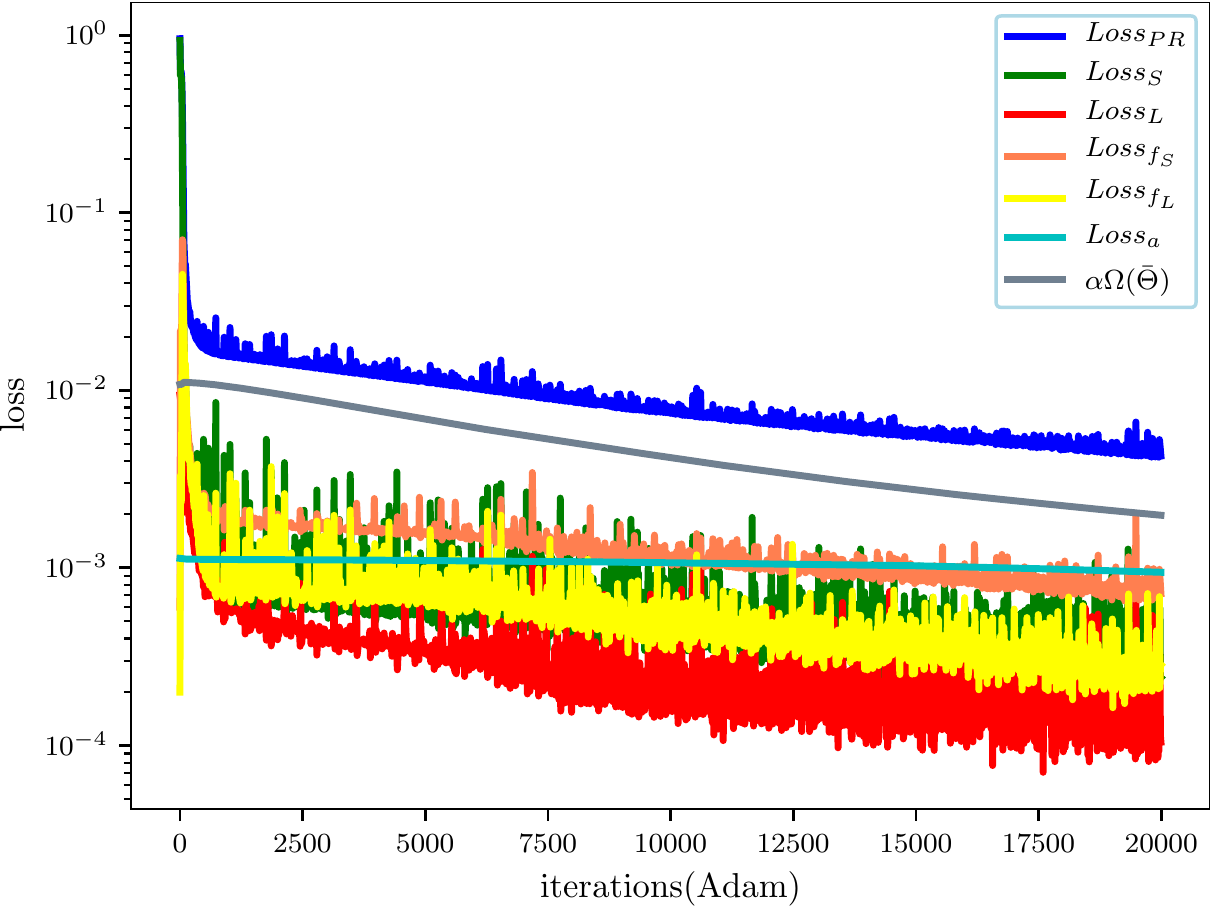}
\end{minipage}
}%
\subfigure[]{
\begin{minipage}[t]{0.48\textwidth}
\centering
\includegraphics[height=4cm,width=6cm]{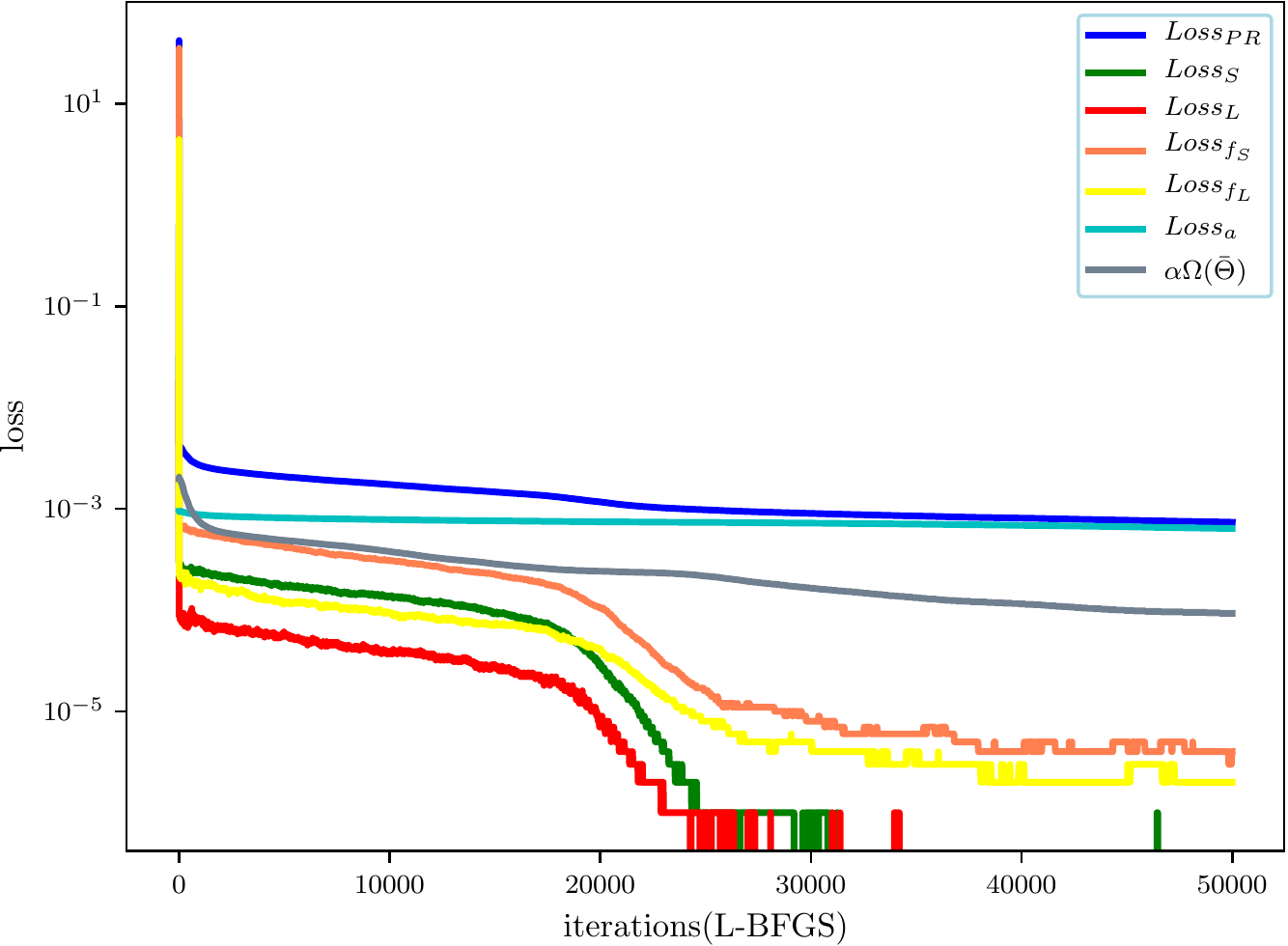}
\end{minipage}%
}%
\centering
\caption{(Color online) The loss function curve figures of the intermediate-bright RWs $S(x,t)$ and $L(x,t)$ arising from the improved PINN with the 20000 steps Adam and 50000 steps L-BFGS optimizations: (a) The loss function curve for the 20000 Adam optimization iterations; (b) The loss function curve for the 50000 L-BFGS optimization iterations.}
\label{F5}
\end{figure}

Figs. \ref{F4} - \ref{F5} display the training results of the data-driven intermediate-bright RWs $S(x,t)$ and $L(x,t)$ based on the improved PINN related to the initial boundary value problem \eqref{E17} and \eqref{E18} of the YO system \eqref{E1}. The left panels of Fig. \ref{F4} depict various dynamic density plots with corresponding amplitude scale size on the right side and sectional drawing at different moments, in which the panel (a) corresponds to short-wave intermediate RW $S(x,t)$ and panel (c) corresponds to the long-wave bright RW $L(x,t)$ for the YO system \eqref{E1}. As we can see from the right panels in Fig. \ref{F4}, the 3D plots with contour map on three planes for the intermediate-bright RWs are shown in Fig. (b) and Fig. (d) respectively. Apparently, from Fig. \ref{F4}, one can see that the maximum amplitudes of short wave RW and long wave RW are lower than those of the two RWs in Fig. \ref{F2}. Fig. \ref{F5} showcases curve plots of the loss function after 20000 Adam optimization iterations and 50000 L-BFGS optimization iterations in improved PINN framework. Different from the curve plots of the loss function in Figs. \ref{F3}, the loss function curves of Adam optimization for the intermediate-bright RWs descend steadily. However, in the process of optimizing the loss function using the L-BFGS optimizer, the loss function curve descends faster after about 20000 iterations, which is different from Fig. \ref{F3}.

\subsection{The data-driven dark-bright RWs}

As we all know, it is extremely difficult to observe and capture the generation and evolution of dark RW in physical experiments. Therefore, we consider utilize the computer-dependent deep learning method to recover the dark RW of short wave $|S|$ by means of the initial boundary value conditions. Similarly, considering the initial condition and Dirichlet boundary condition of the YO system to obtain the dark-bright RWs by using the 9-layer improved PINN with 40 neurons per layer, the $[X_0, X_1]$ and $[T_0, T_1]$ in Eq. \eqref{E1} are taken as $[-5.0, 5.0]$ and $[-3.0, 3.0]$, respectively. We immediately obtain the initial value conditions
\begin{align}\label{E19}
\begin{split}
&S^0(x)=S_{\rm drw}(x,-3.0),\,L^0(x)=L_{\rm brw3}(x,-3.0),\,x\in[-5.0,5.0],
\end{split}
\end{align}
and the Dirichlet boundary conditions
\begin{align}\label{E20}
\begin{split}
&S^{\mathrm{lb}}(t)=S_{\rm drw}(-5.0,t),\,L^{\mathrm{lb}}(t)=L_{\rm brw3}(-5.0,t),\,t\in[-3.0,3.0],\\
&S^{\mathrm{ub}}(t)=S_{\rm drw}(5.0,t),\,L^{\mathrm{ub}}(t)=L_{\rm brw3}(5.0,t),\,t\in[-3.0,3.0].
\end{split}
\end{align}

Similarly, discretizing exact dark-bright RWs $S_{\rm drw}$ and $L_{\rm brw3}$ with the aid of the traditional finite difference scheme on even grids in Matlab, then we obtain the original training data which contains initial data \eqref{E19} and boundary data \eqref{E20} by separately dividing the spatial region $[-5.0, 5.0]$ into 2000 points and the temporal region $[-3.0, 3.0]$ into 1000 points. Then, one can yield a smaller training dataset that contains partial initial-boundary data by randomly extracting $N_q = 2000$ initial-boundary value data points from original dataset and $N_f = 30000$ collocation points which are produced by the LHS. After that, the latent dark-bright RWs $S(x,t)$ and $L(x,t)$ have been successfully learned by tuning all learnable parameters of the improved PINN, and the network achieved relative $\mathbb{L}_2$ error of 1.964839$\rm e$-03 for the dark RW $S(x,t)$ and relative $\mathbb{L}_2$ error of 1.692152$\rm e$-02 for the bright RW $L(x,t)$, and the total number of iterations is 70000.

\begin{figure}[htbp]
\centering
\subfigure[]{
\begin{minipage}[t]{0.48\textwidth}
\centering
\includegraphics[height=3.5cm,width=7cm]{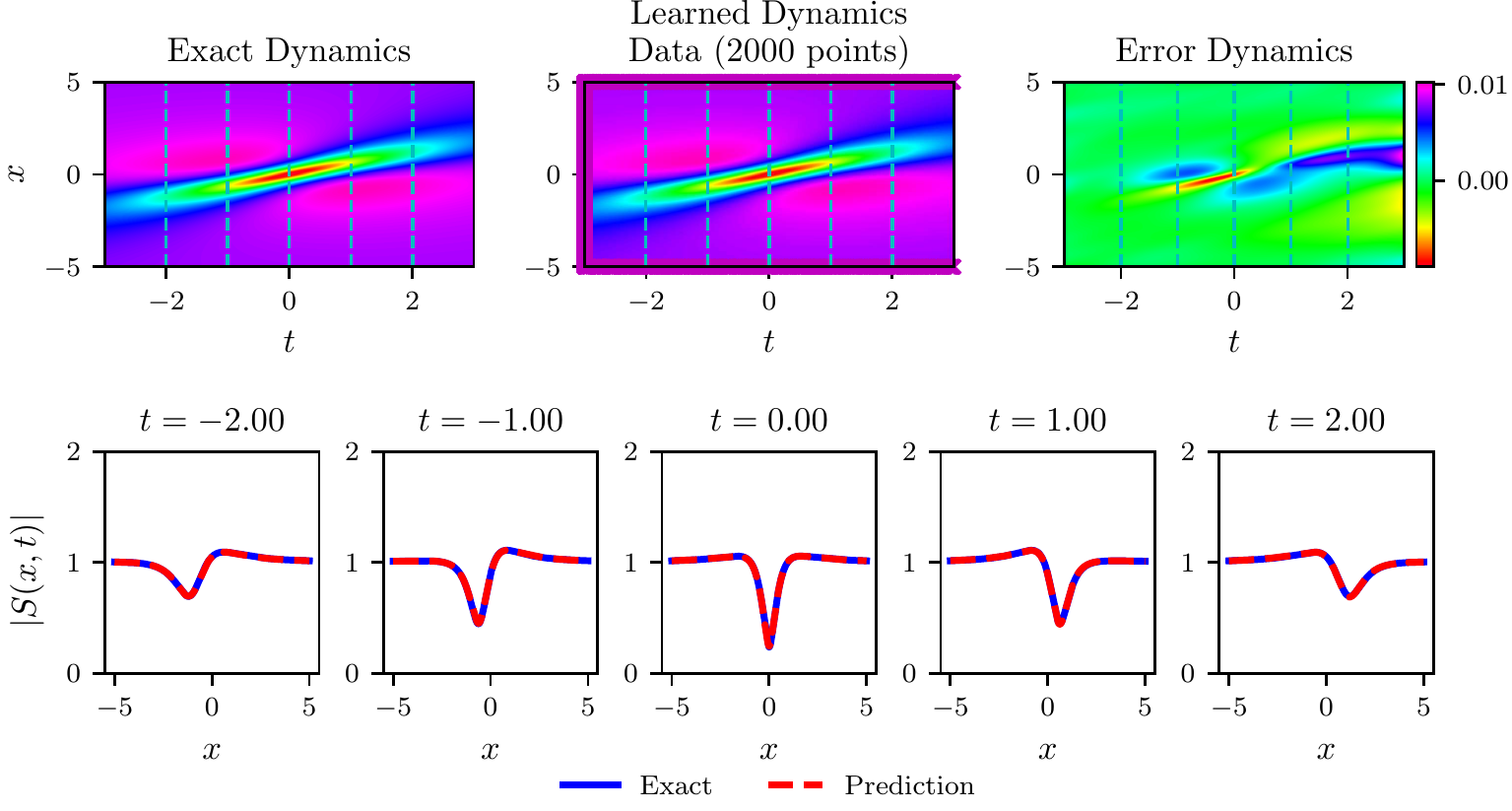}
\end{minipage}
}%
\subfigure[]{
\begin{minipage}[t]{0.48\textwidth}
\centering
\includegraphics[height=3.5cm,width=6.5cm]{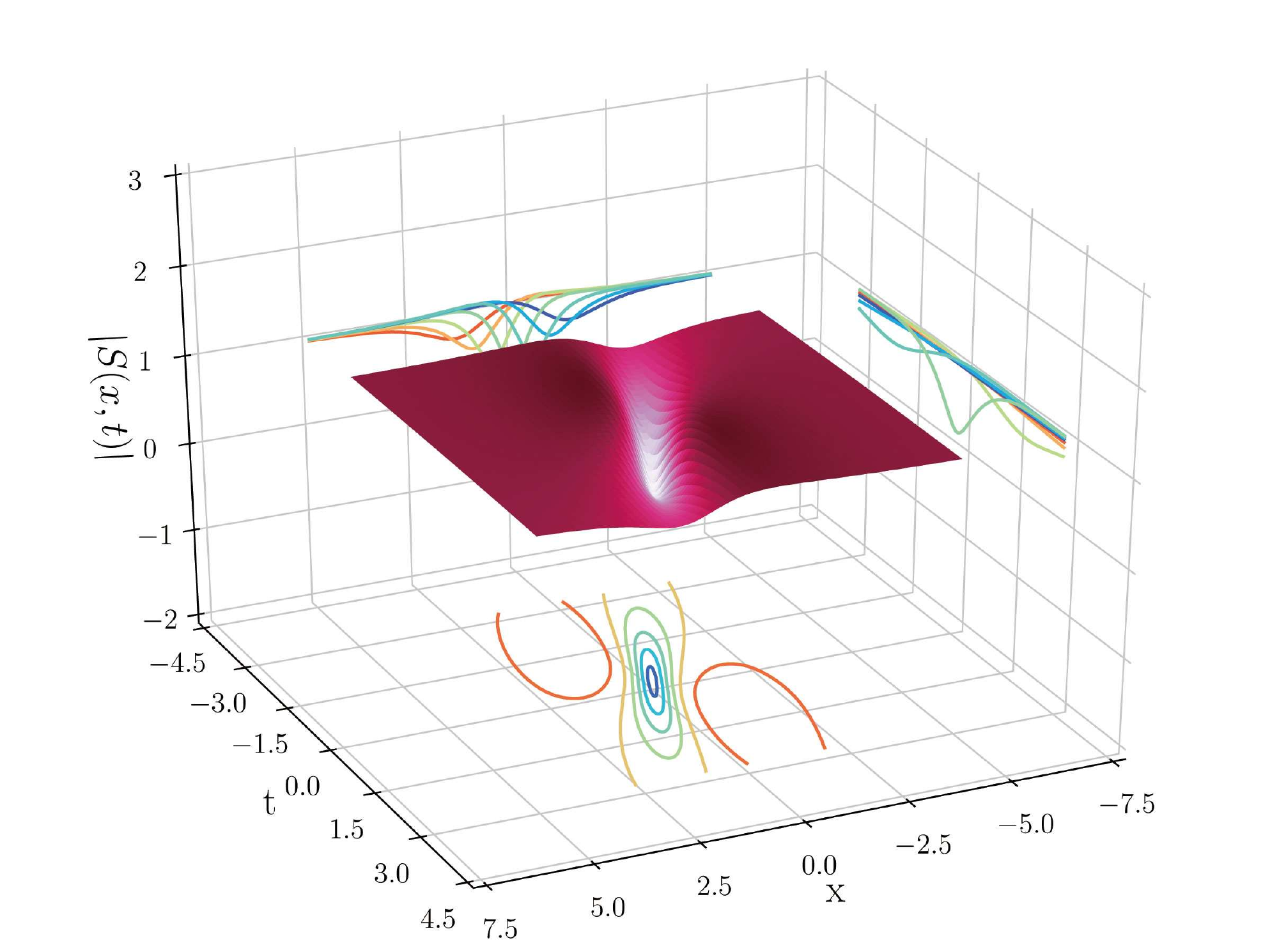}
\end{minipage}%
}\\%
\subfigure[]{
\begin{minipage}[t]{0.48\textwidth}
\centering
\includegraphics[height=3.5cm,width=7cm]{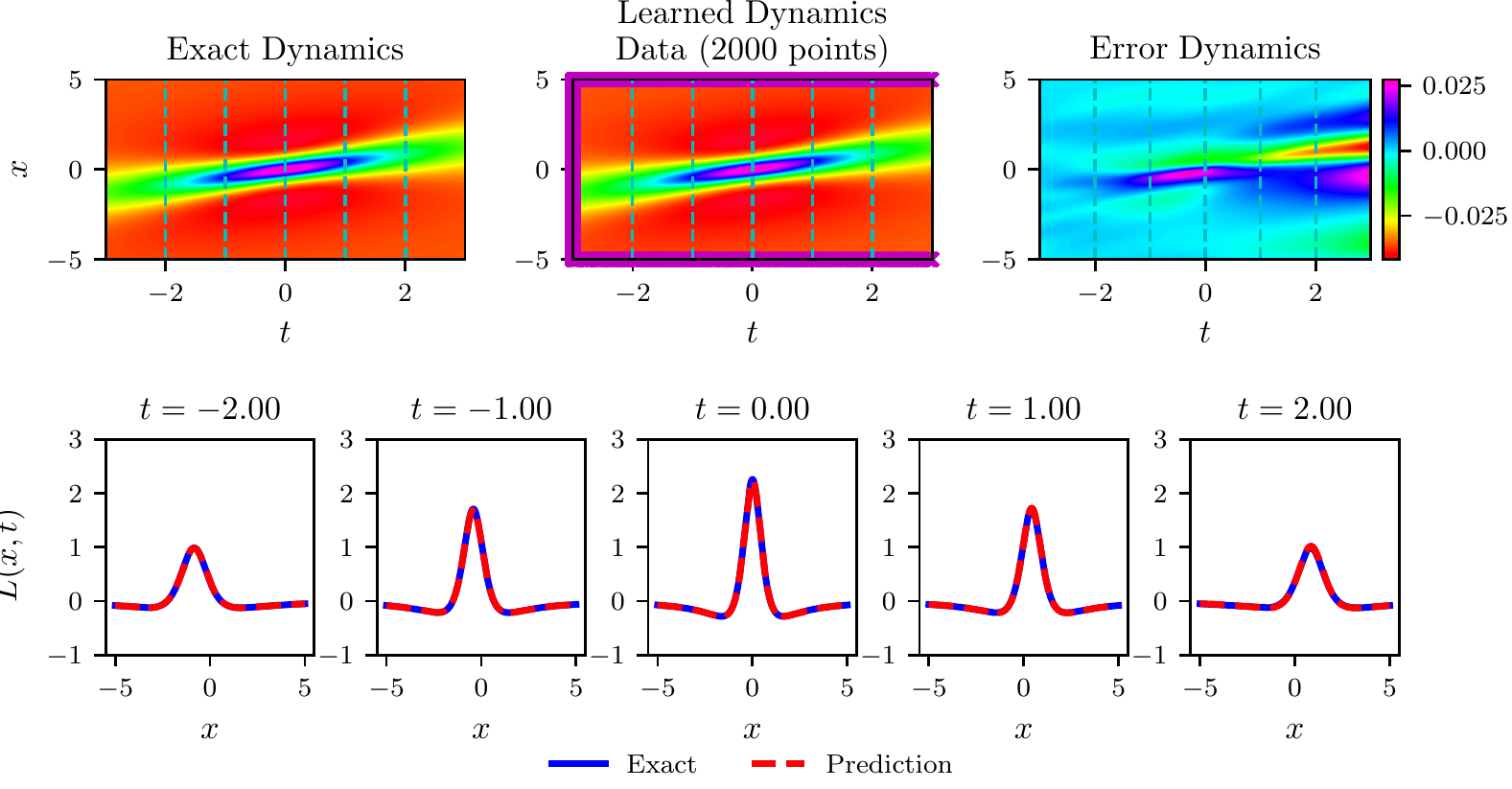}
\end{minipage}
}%
\subfigure[]{
\begin{minipage}[t]{0.48\textwidth}
\centering
\includegraphics[height=3.5cm,width=6.5cm]{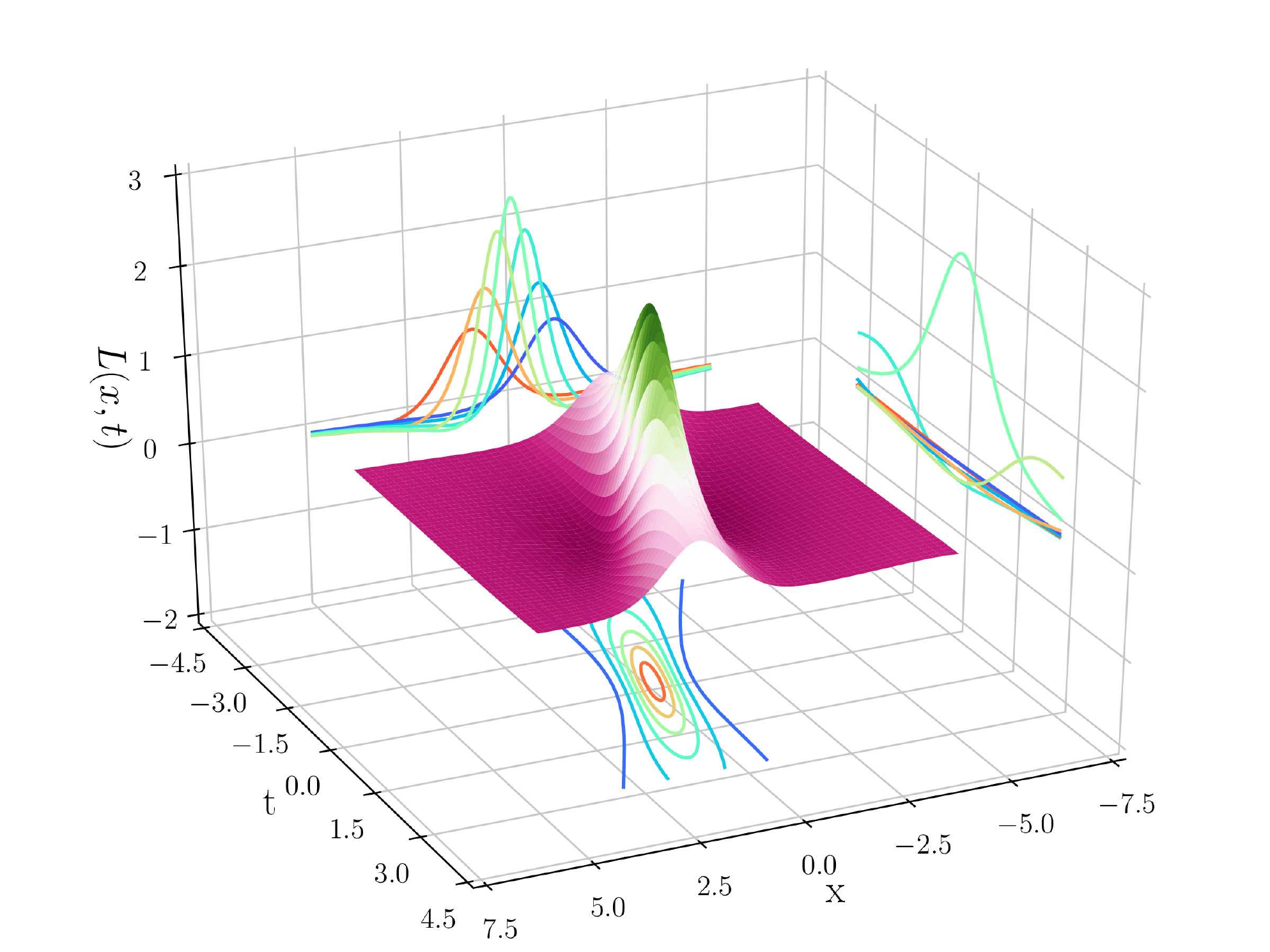}
\end{minipage}%
}%
\centering
\caption{(Color online) The data-driven dark-bright RWs $S(x,t)$ and $L(x,t)$ arising from the improved PINN with parameter regularization of hyper-parameter $\alpha=0.0001$ by randomly choosing $N_q=2000$ initial and boundary points which have been shown by using mediumorchid $``\times"$ in learned dynamics , as well as $N_f = 30000$ collocation points in the corresponding spatiotemporal region: (a) and (c) The exact, learned and error dynamics density plots with five distinct tested times $t=-2.00, -1.00, 0.00, 1.00$ and 2.00 (darkturquoise dashed lines), as well as sectional drawings which contain the learned and explicit dark-bright RWs at the aforementioned five distinct times; (b) and (d) The 3D plot with contour map for the data-driven dark-bright RWs.}
\label{F6}
\end{figure}

\begin{figure}[htbp]
\centering
\subfigure[]{
\begin{minipage}[t]{0.48\textwidth}
\centering
\includegraphics[height=4cm,width=6cm]{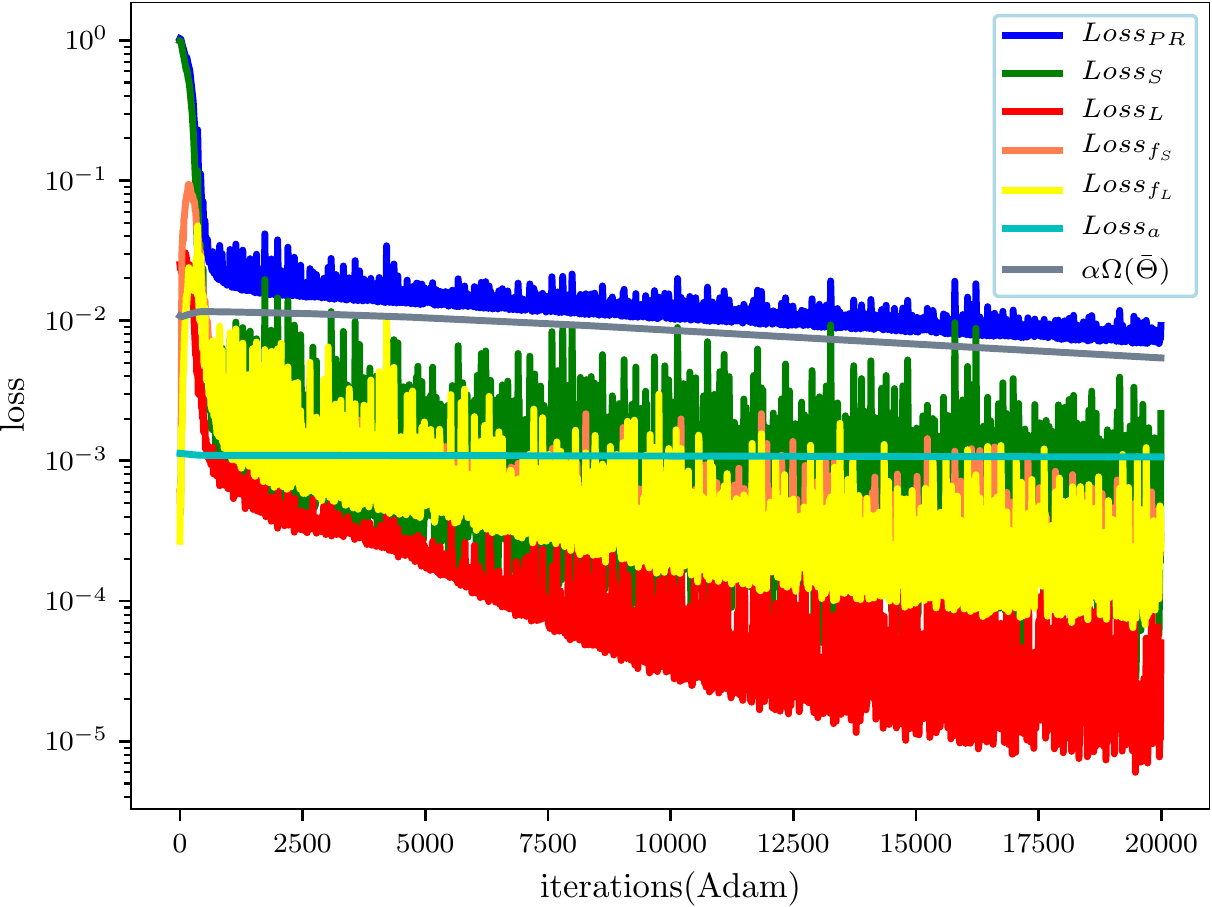}
\end{minipage}
}%
\subfigure[]{
\begin{minipage}[t]{0.48\textwidth}
\centering
\includegraphics[height=4cm,width=6cm]{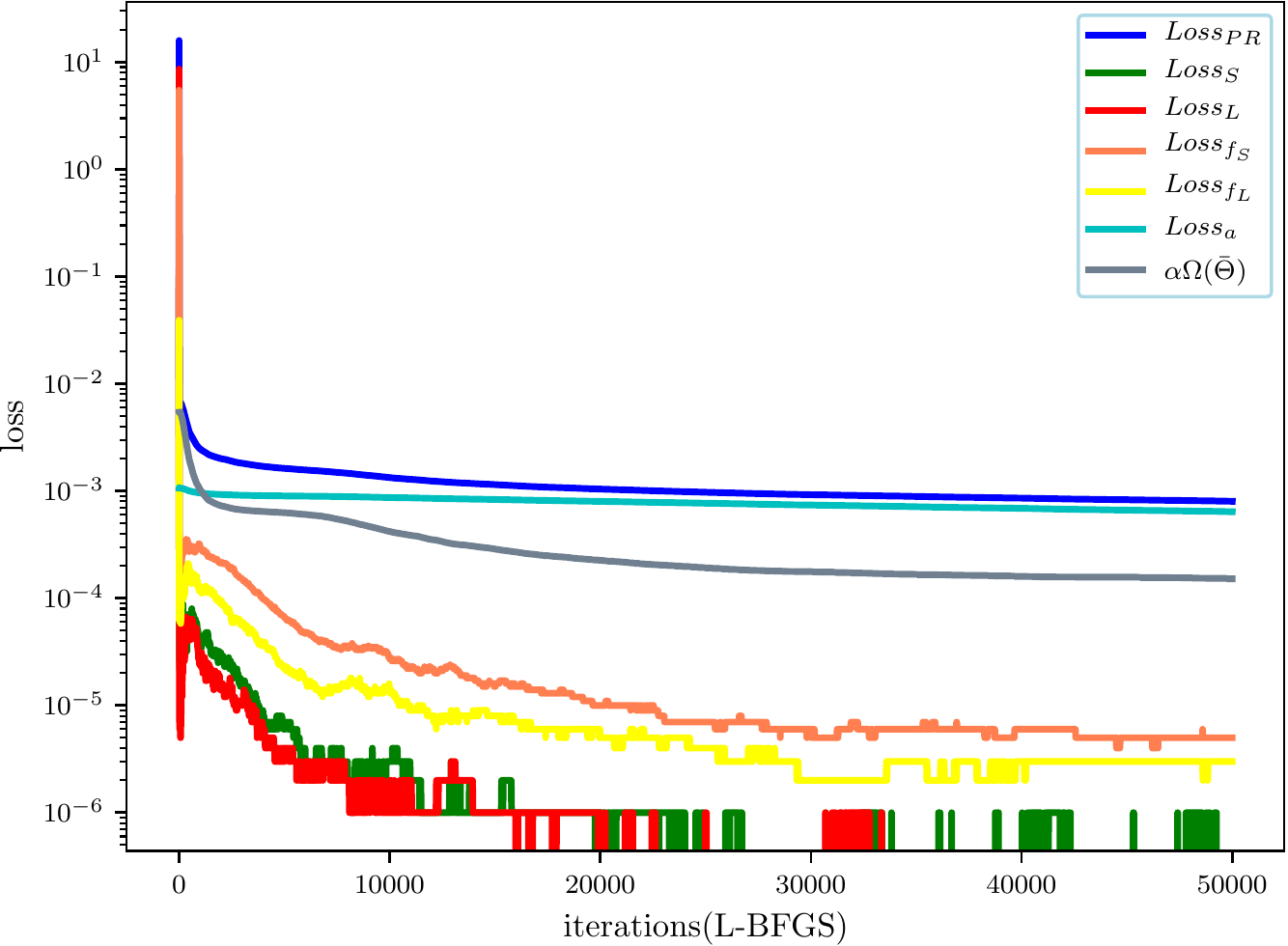}
\end{minipage}%
}%
\centering
\caption{(Color online) The loss function curve figures of the dark-bright RWs $S(x,t)$ and $L(x,t)$ arising from the improved PINN with the 20000 steps Adam and 50000 steps L-BFGS optimizations: (a) The loss function curve for the 20000 Adam optimization iterations; (b) The loss function curve for the 50000 L-BFGS optimization iterations.}
\label{F7}
\end{figure}

Figs. \ref{F6} - \ref{F7} provide the training results arising from the improved PINN for the data-driven dark-bright RWs $S(x,t)$ and $L(x,t)$ of the YO system with the initial boundary value problem \eqref{E19} and \eqref{E20}. In the left panels of Fig. \ref{F6}, the exact, learned and error dynamics density plots with corresponding amplitude scale size on the right side have been exhibited, it is worth mentioning that the $N_q=2000$ training data points involved in the initial-boundary condition are marked by mediumorchid symbol $``\times"$ in the learned density plots both in (a) and (c) of Fig. \ref{F6}. Meanwhile, the sectional drawings which include the learned and exact dark-bright RWs have been shown at the five distinct times pointed out in the exact, learned and error dynamics density plots by using darkturquoise dashed lines in the bottom panels of (a) and (c). The right panels of Fig. \ref{F6} display the three-dimensional plots with contour map on three planes for the predicted dark-bright RWs $S(x,t)$ and $L(x,t)$ based on the improved PINN. Fig. \ref{F7} exhibits the loss function curve figures of the dark-bright RWs $S(x,t)$ and $L(x,t)$ arising from the improved PINN with the 20000 steps Adam and 50000 steps L-BFGS optimizations on the loss function $\widetilde{\mathscr{L}}(\bar{\Theta})$. We are surprised to find that the generation and evolution process for dark RW of short wave in a certain spatiotemporal interval can be accurately recovered through the deep learning method, which will provide more opportunities for many physical experiment designs involving dark RW.

In addition, a large number of experimental data show that the training error of improved PINN with parameter regularization ($\alpha=0.0001$) is smaller than that of PINN without parameter regularization ($\alpha=0$). Tab. \ref{Tab1} gives the relative $\mathbb{L}_2$ norm errors of three different types of RWs with and without parameter regularization. From Tab. \ref{Tab1}, once the hyper-parameter $\alpha$ is small enough, one can see that the training error of the improved PINN model with parameter regularization is mostly lower than that of the PINN model without parameter regularization.

\begin{table}[htbp]
  \caption{Relative $\mathbb{L}_2$ errors of three different RW types in diverse PINN types}
  \label{Tab1}
  \centering
  \scalebox{0.8}{
  \begin{tabular}{l|c|c|c}
  \toprule
  \diagbox{\small{\textbf{PINN Types}}}{\small{\textbf{RW Types}}} & Bright-bright RWs & Intermediate-bright RWs & Dark-bright RWs\\
  \hline
  Hyper-parameter $\alpha=0$   & \makecell[c]{$S(x,t)$: 7.566667$\rm e$-04\\ $L(x,t)$: 2.414187$\rm e$-03} & \makecell[c]{$S(x,t)$: 1.975757$\rm e$-03\\ $L(x,t)$: 8.500360$\rm e$-03}& \makecell[c]{$S(x,t)$: 1.788549$\rm e$-03\\ $L(x,t)$: 1.286998$\rm e$-02}\\
  \hline
  Hyper-parameter $\alpha=0.0001$   & \makecell[c]{$S(x,t)$: 4.968430$\rm e$-04\\ $L(x,t)$: 1.763312$\rm e$-03} & \makecell[c]{$S(x,t)$: 1.168852$\rm e$-03\\ $L(x,t)$: 6.766132$\rm e$-03}& \makecell[c]{$S(x,t)$: 1.964839$\rm e$-03\\ $L(x,t)$: 1.692152$\rm e$-02}\\
  \bottomrule
  \end{tabular}}
\end{table}

\section{The inverse problem of the YO system}

In this section, we focus on the inverse problem of the YO system, that is parameters discovery problem for a data-driven YO system model \eqref{E1} by utilizing small training data set. In this situation, unknown parameters $\lambda_1$ and $\lambda_2$ are learned by introducing $L^2$ norm regularization into 9-layer improved PINN with 40 neurons per layer, and the physics-informed parts Eq. \eqref{E6} of the improved PINN for YO sysem \eqref{E1} become the following formula
\begin{align}\label{E21}
\begin{split}
&f_u:=-v_t+\lambda_1u_{xx}+uL,\quad f_v:=u_t+\lambda_1v_{xx}+vL,\quad f_L:=L_t-\lambda_2(2uu_x+2vv_x).
\end{split}
\end{align}

In order to learn the unknown parameters $\lambda_1$ and $\lambda_2$ in Eq. \eqref{E1} with the aid of the improved PINN with neuron-wise locally adaptive activation functions and parametric regularization term of different trade-off coefficients, and considering the initial conditions and Dirichlet boundary conditions of Eq. \eqref{E1} arising from the bright-bright RWs \eqref{E14}, here we set the spatiotemporal region $(x,t)\in[-5,5]\times[-0.5,0.5]$. Thus the corresponding initial conditions can be written as belows
\begin{align}\label{E22}
\begin{split}
&S^0(x)=S_{\rm brw}(x,-0.5),\,L^0(x)=L_{\rm brw1}(x,-0.5),\,x\in[-5.0,5.0],
\end{split}
\end{align}
and the Dirichlet boundary conditions
\begin{align}\label{E23}
\begin{split}
&S^{\mathrm{lb}}(t)=S_{\rm brw}(-5.0,t),\,L^{\mathrm{lb}}(t)=L_{\rm brw1}(-5.0,t),\,t\in[-0.5,0.5],\\
&S^{\mathrm{ub}}(t)=S_{\rm brw}(5.0,t),\,L^{\mathrm{ub}}(t)=L_{\rm brw1}(5.0,t),\,t\in[-0.5,0.5].
\end{split}
\end{align}

Here, we employ the same data discretization method in section 3, and produce the training data which consists of initial data \eqref{E22} and boundary data \eqref{E23} by dividing the spatial region $[-5.0,5.0]$ into 2000 points and temporal region $[-0.5,0.5]$ into 1000 points. Then we obtain a smaller training dataset that containing initial-boundary data by randomly extracting $N_q=2000$ initial boundary value data points from original dataset and $N_f=30000$ collocation points which are generated by the LHS method. After giving the training dataset, the latent data-driven unknown parameters $\lambda_1$ and $\lambda_2$ have been successfully predicted by tuning all learnable parameters of the improved PINN with 20000 Adam iterations and different number of L-BFGS iterations to regulate the loss function $\widetilde{\mathscr{L}}(\bar{\Theta})$. The unknown parameters $\lambda_1$ and $\lambda_2$ are initialized to $\lambda_1=\lambda_2=0$. the relative error of unknown parameters is defined as $RE=(|\hat{\lambda}_{\kappa}-\lambda_{\kappa}|/ \lambda_{\kappa})\times100\%$ $(\kappa=1, 2)$ with the predicted value $\hat{\lambda}_{\kappa}$ and true value $\lambda_{\kappa}$. All noise interference in this part is added to the randomly chosen small data set, the details are as shown below
\begin{align}\nonumber
\begin{split}
Data_{-}train = &Data_{-}train + noise*np.std(Data_{-}train)*np.random.randn\\
&(Data_{-}train.shape[0], Data_{-}train.shape[1]),
\end{split}
\end{align}
where $Data_{-}train$ and $noise$ represent a small randomly chosen training data set and the noise intensity, respectively. The $np.std(\cdot)$ returns the standard deviation of an array element, and $np.random.randn(\cdot,\cdot)$ returns a set of samples with a standard normal distribution.

\begin{figure}[htbp]
\centering
\subfigure[]{
\begin{minipage}[t]{0.48\textwidth}
\centering
\includegraphics[height=4.0cm,width=6.5cm]{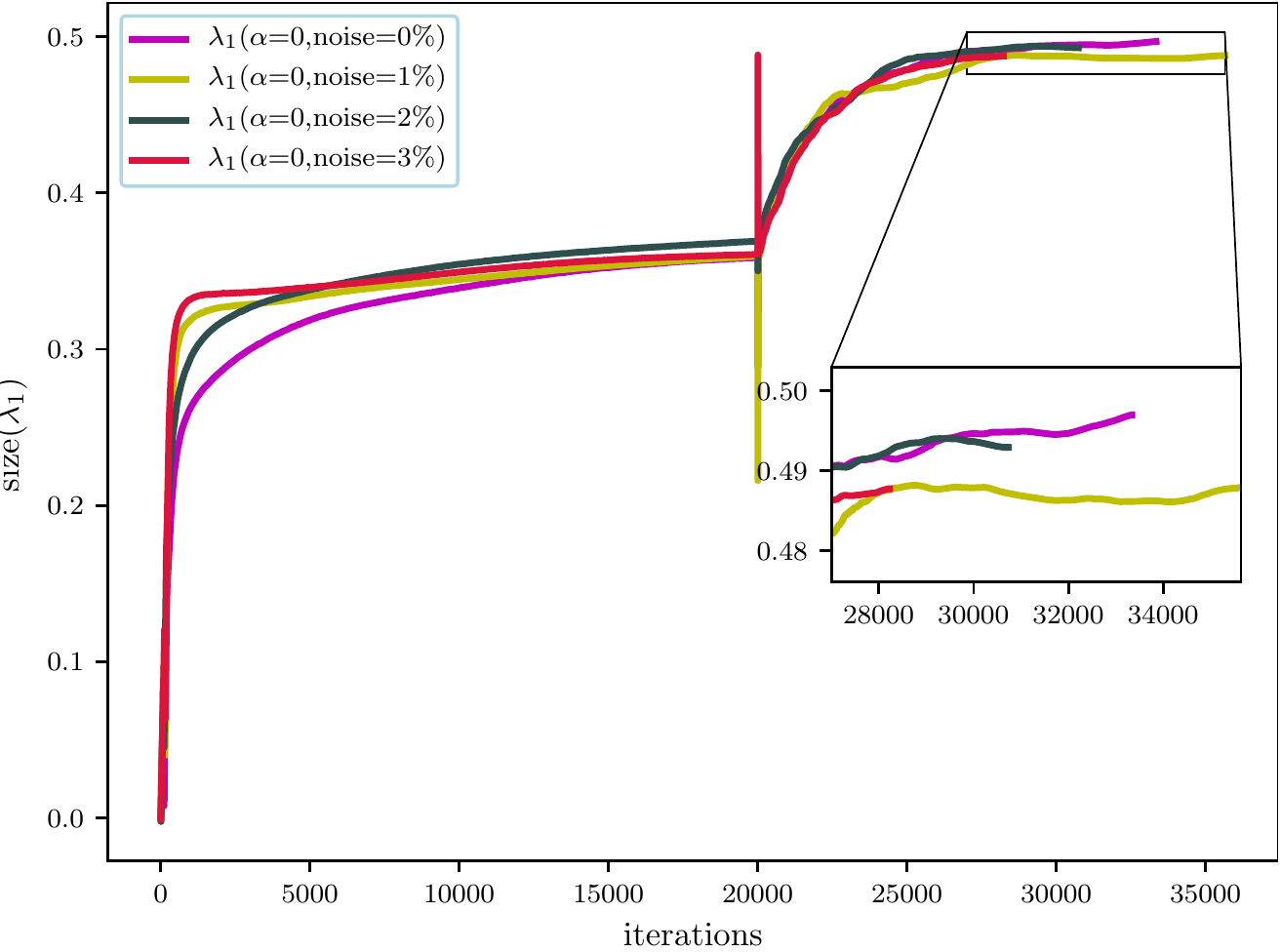}
\end{minipage}
}%
\subfigure[]{
\begin{minipage}[t]{0.48\textwidth}
\centering
\includegraphics[height=4.0cm,width=6.5cm]{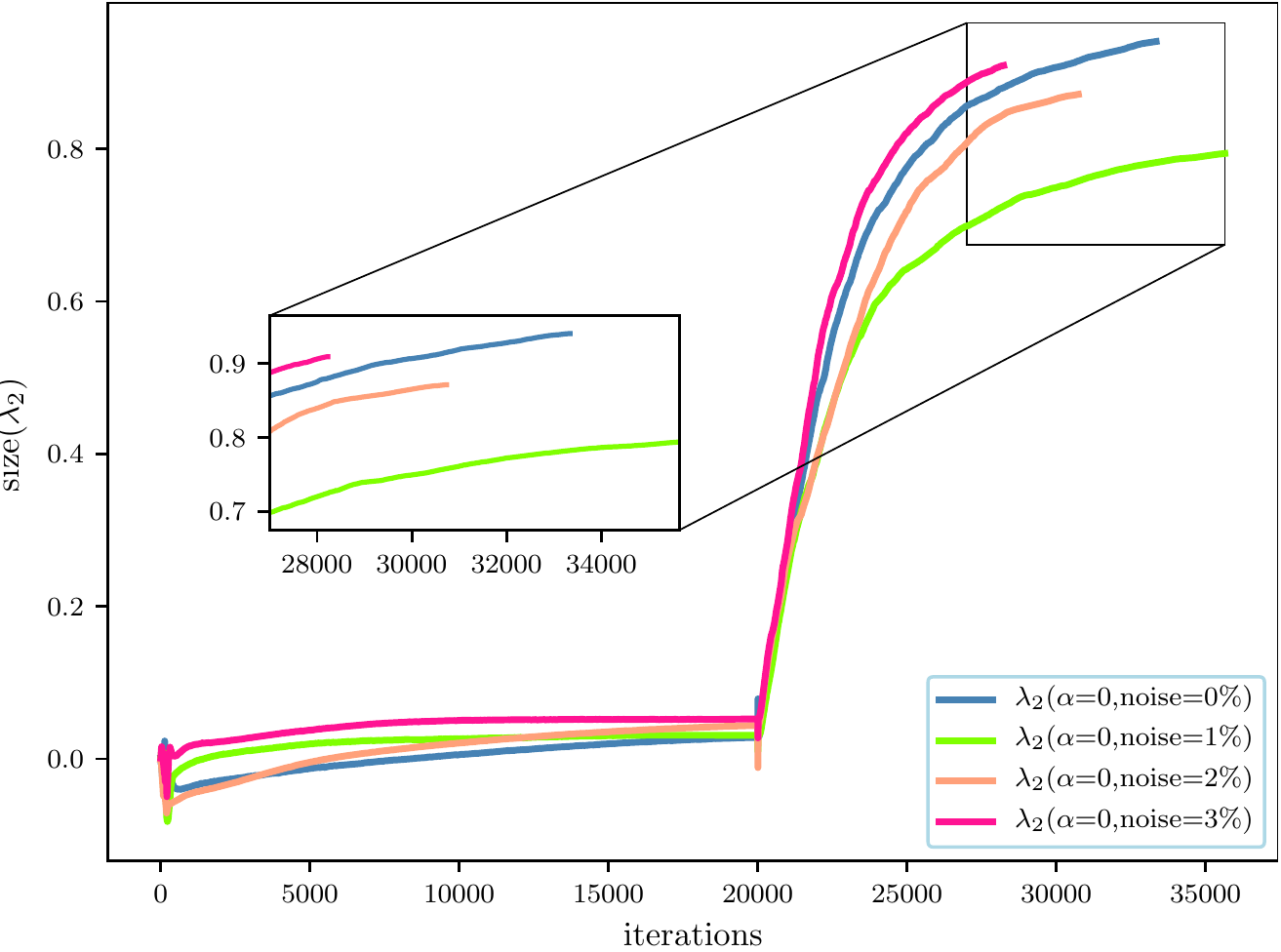}
\end{minipage}%
}\\%
\subfigure[]{
\begin{minipage}[t]{0.48\textwidth}
\centering
\includegraphics[height=4.0cm,width=6.5cm]{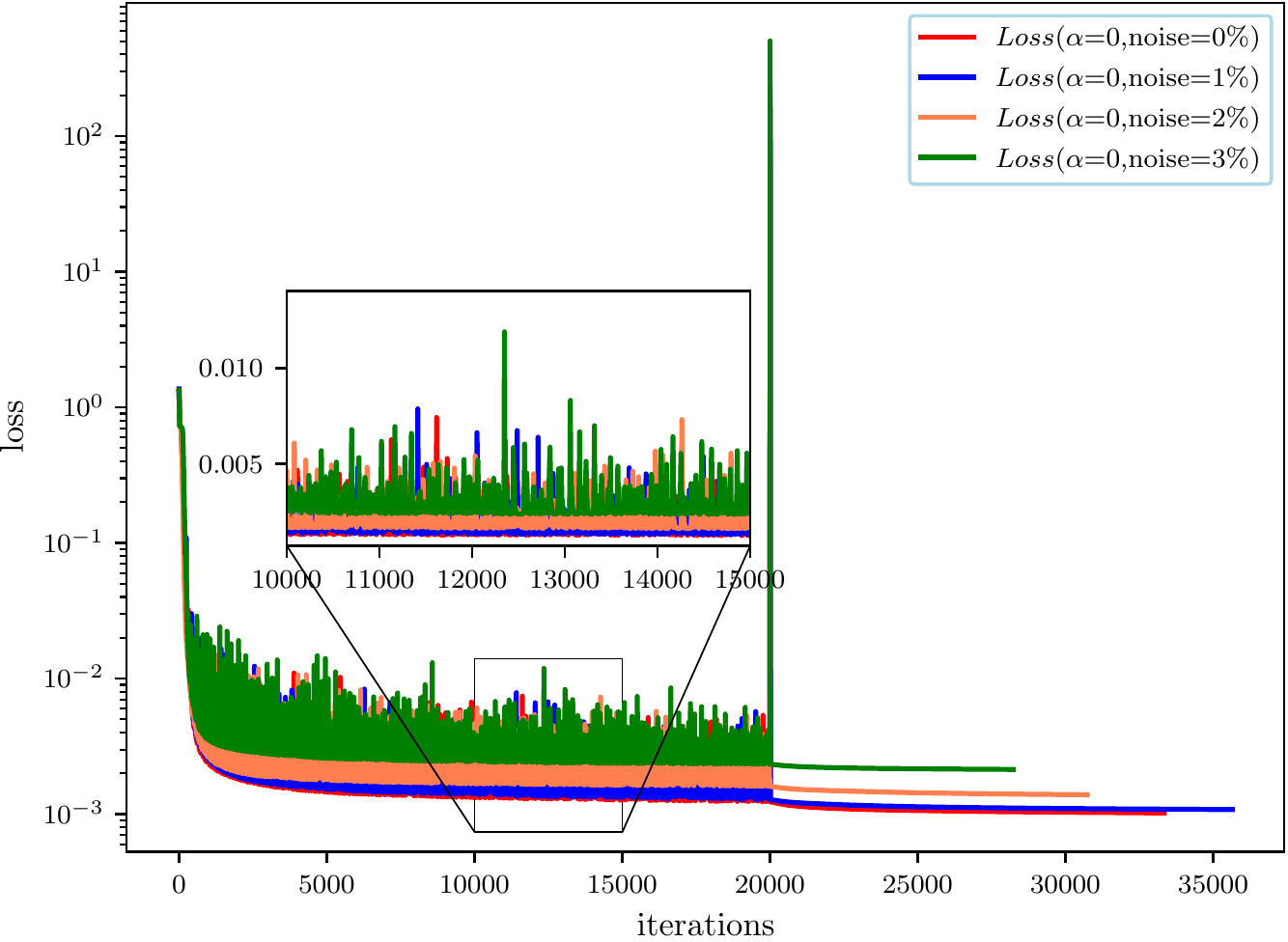}
\end{minipage}
}%
\subfigure[]{
\begin{minipage}[t]{0.48\textwidth}
\centering
\includegraphics[height=4.0cm,width=6.5cm]{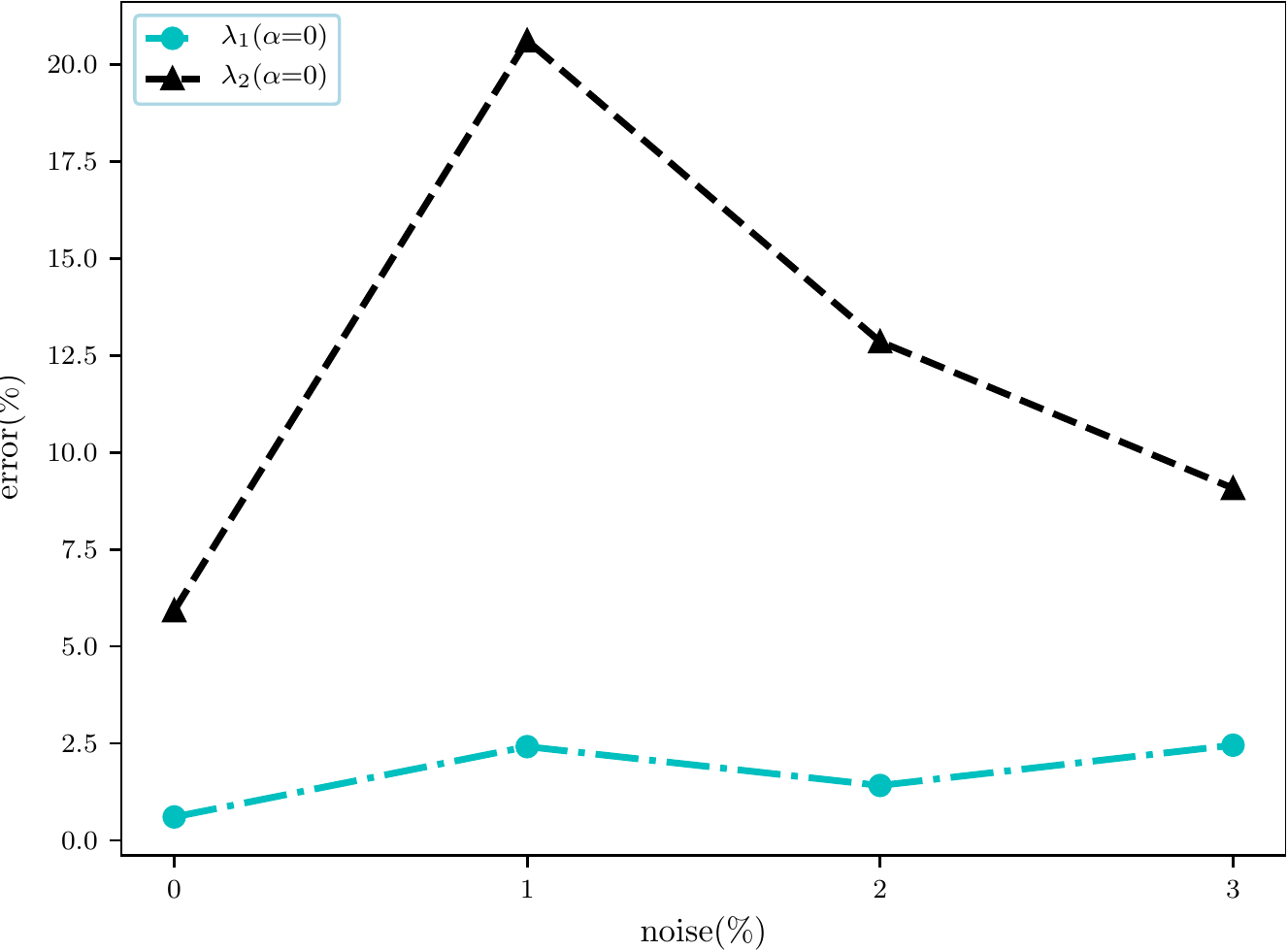}
\end{minipage}%
}%
\centering
\caption{(Color online) The parameter discover arising from the PINN without parameter regularization ($\alpha=0$): (a)-(b) the variation of unknown coefficients $\lambda_1$ and $\lambda_2$ with different noise intensity; (c) the variation of loss function  with different noise intensity; (d) unknown coefficients $\lambda_1$ and $\lambda_2$ error variation under different interference noise.}
\label{F8}
\end{figure}

Next we analyze the training result of the NN from different perspectives, such as the size of hyper-parameters $\alpha$, intensity of noise and relative error of prediction parameters. In order to more directly verify the effect of improved PINN with parameter regularization, we first showcase the training effect of PINN without parameter regularization (namely $\alpha=0$) in Fig. \ref{F8}. In the absence of a parametric regularization strategy, (a) and (b) of Fig. \ref{F8} describe the numerical variation curves of unknown parameters $\lambda_1$ and $\lambda_2$ during iteration, one can find that $\lambda_1$ increases from 0 to more than 0.3 during the previous 20000 Adam optimizations, while $\lambda_2$ hardly increases significantly in the aforementioned iterations. Instead, $\lambda_1$ increases slowly to about 0.5 in the later L-BFGS optimization process, while $\lambda_2$ increases sharply to about 1.0 in this iterative process. The panel (c) of Fig. \ref{F8} indicates that the greater the noise intensity, the more intense the fluctuation of the loss function curve in the first 20000 Adam optimization processes, and the larger the overall value of the loss function in the later L-BFGS optimization processes. Fig. \ref{F8} (d) shows that the relative error of $\lambda_2$ is more sensitive to the change of noise intensity than that of $\lambda_1$. Under four different noise intensity environments, one can observe that the overall relative error of $\lambda_2$ is greater than that of $\lambda_1$ from Fig. \ref{F8} (d), this is because there are two physical constraints related to $\lambda_1$, while there is only one physical constraint related to $\lambda_2$ in the physics-informed parts of PINN.

\begin{figure}[htbp]
\centering
\subfigure[]{
\begin{minipage}[t]{0.48\textwidth}
\centering
\includegraphics[height=4.0cm,width=6.5cm]{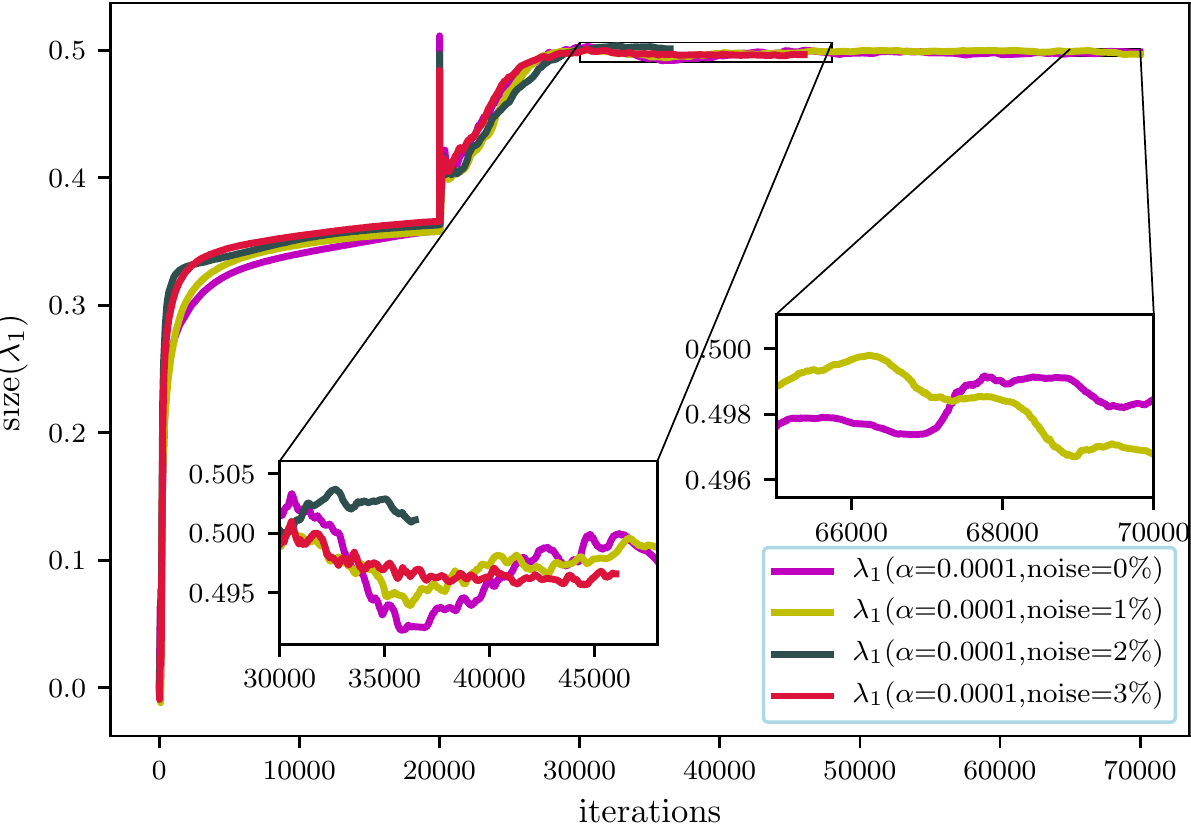}
\end{minipage}
}%
\subfigure[]{
\begin{minipage}[t]{0.48\textwidth}
\centering
\includegraphics[height=4.0cm,width=6.5cm]{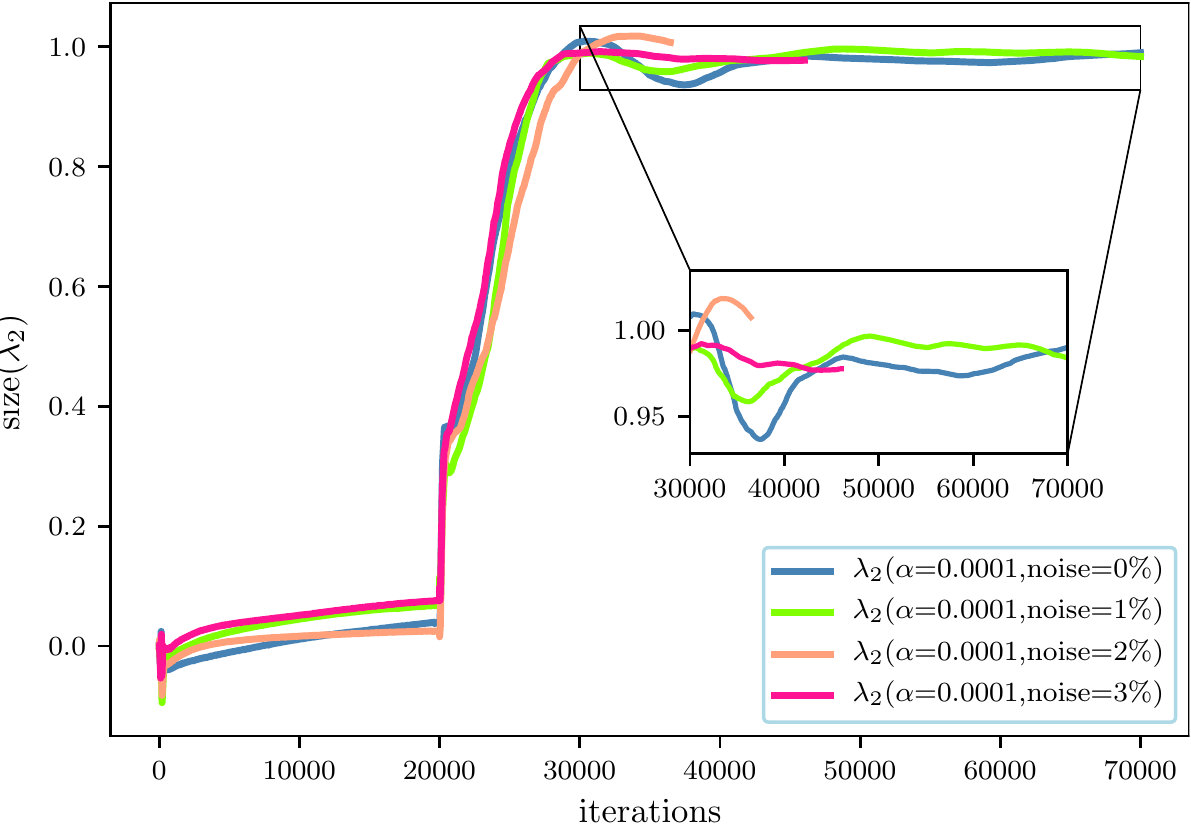}
\end{minipage}%
}\\%
\subfigure[]{
\begin{minipage}[t]{0.48\textwidth}
\centering
\includegraphics[height=4.0cm,width=6.5cm]{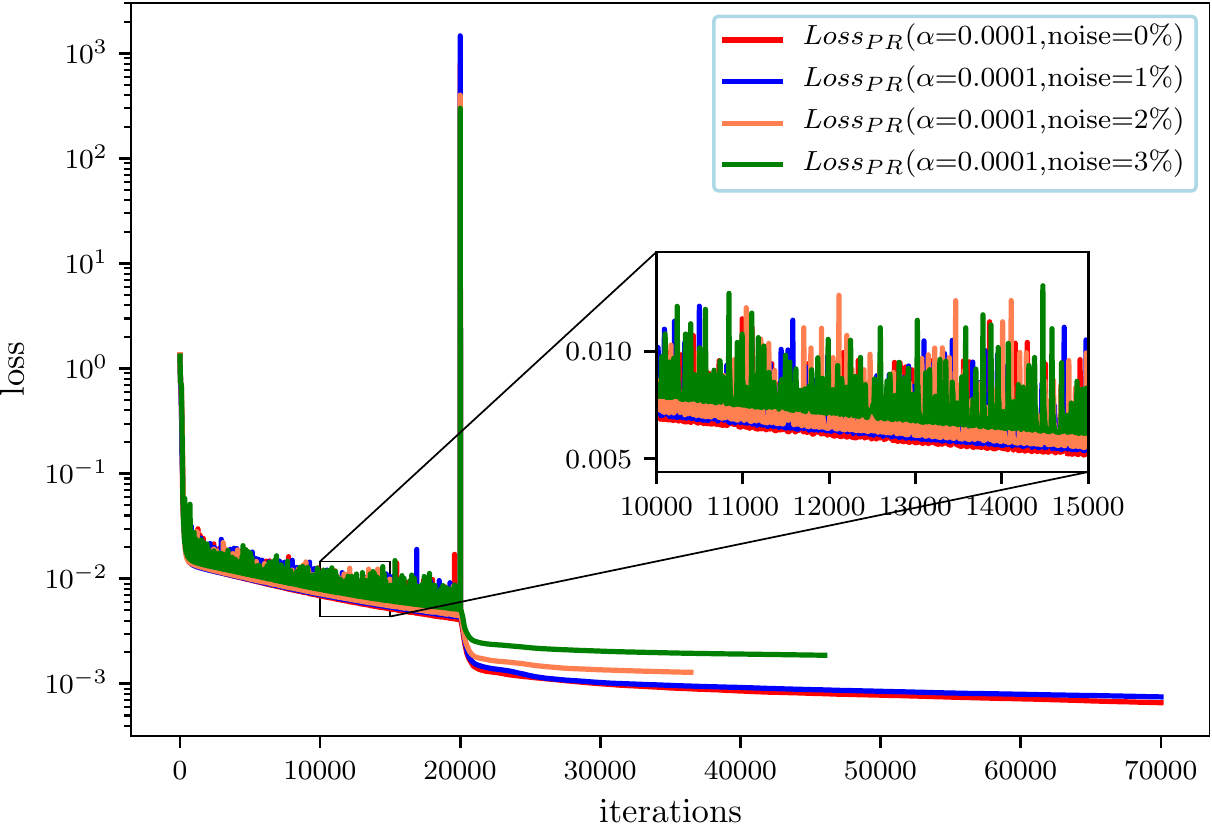}
\end{minipage}
}%
\subfigure[]{
\begin{minipage}[t]{0.48\textwidth}
\centering
\includegraphics[height=4.0cm,width=6.5cm]{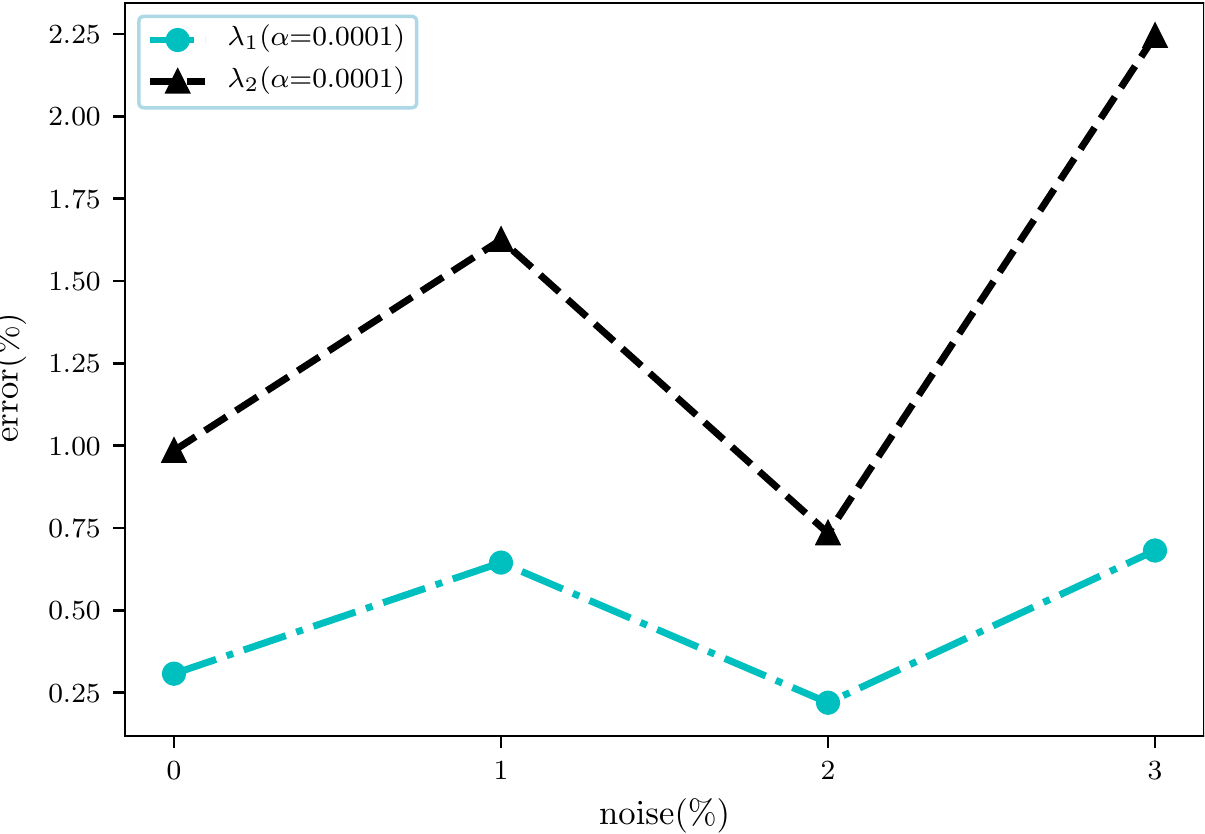}
\end{minipage}%
}%
\centering
\caption{(Color online) The parameter discover arising from the improved PINN with parameter regularization ($\alpha=0.0001$): (a)-(b) the variation of unknown coefficients $\lambda_1$ and $\lambda_2$ with different noise intensity; (c) the variation of loss function with different noise intensity; (d) unknown coefficients $\lambda_1$ and $\lambda_2$ error variation under different interference noise.}
\label{F9}
\end{figure}

Next, we impose the parametric regularization strategy with penalty coefficient $\alpha=0.0001$ to the improved PINN, the corresponding training results are shown in Fig. \ref{F9}. Due to the introduction of parameter regularization, the load of loss function will increase, thus the number of iterations will be greater than that in the PINN without parameter regularization strategy in some cases, but the maximum number of iterations is artificially set to 70000 (including 2000 Adam optimizer iterations and 50000 L-BFGS optimizer iterations). Fig. \ref{F9} (a)-(b) show the variation curves of unknown parameters $\lambda_1$ and $\lambda_2$ in the iterative process under different noise intensity conditions, in which both figures indicate that there is little difference between the final learning results of parameters $\lambda_1$ and $\lambda_2$ to be learned in the case of noise and no noise. Interestingly, panel (b) of Fig. \ref{F9} demonstrates that the value of $\lambda_2$ learned with 2$\%$ noise intensity is closer to the real value than that learned without noise case. Different from Fig. \ref{F8} (c), the Fig. \ref{F9} (c) shows that the relationship between the fluctuation degree of the loss function curve and the noise intensity is not obvious. Fig. \ref{F9} (d) indicates that in improved PINN with parameter regularization, the relative error of parameter training results is the smallest when the noise intensity is 2$\%$, which further reveals that improved PINN with parameter regularization strategy has the ability to suppress data noise interference.

\begin{figure}[htbp]
\centering
\subfigure[]{
\begin{minipage}[t]{0.48\textwidth}
\centering
\includegraphics[height=4.0cm,width=6.5cm]{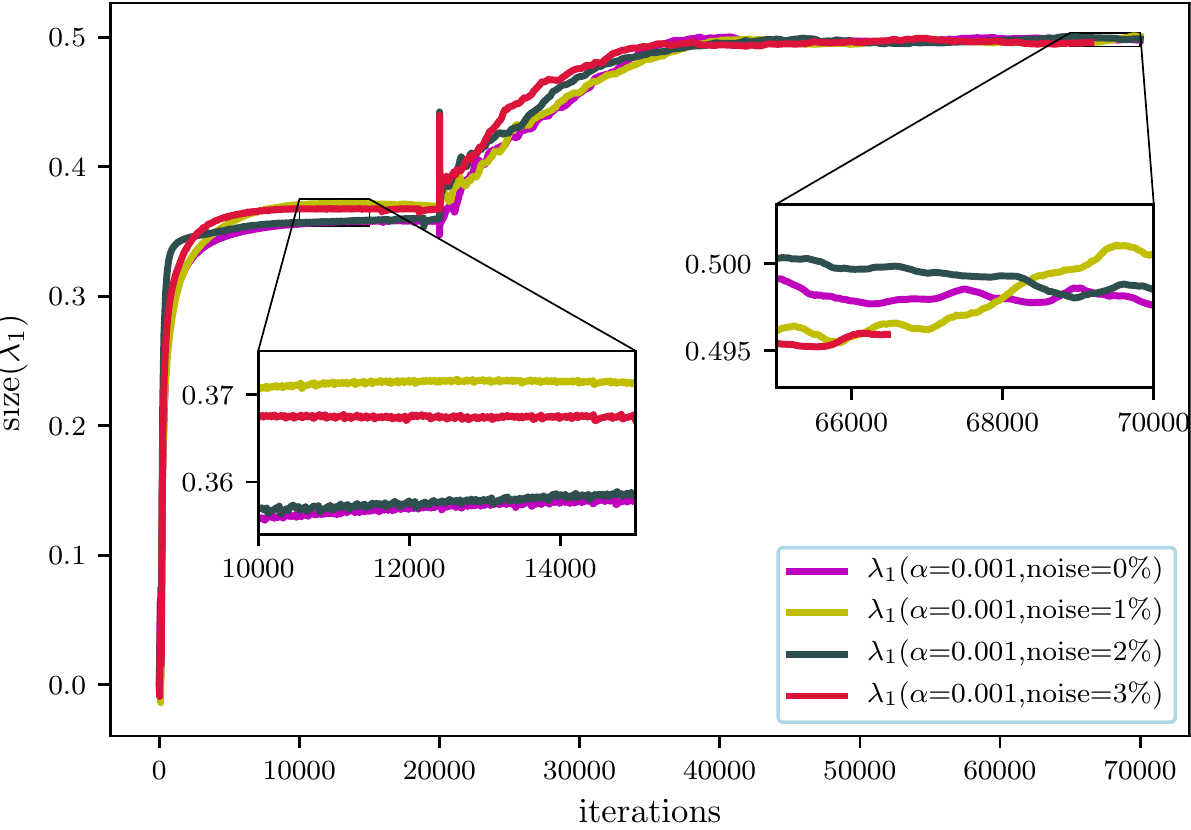}
\end{minipage}
}%
\subfigure[]{
\begin{minipage}[t]{0.48\textwidth}
\centering
\includegraphics[height=4.0cm,width=6.5cm]{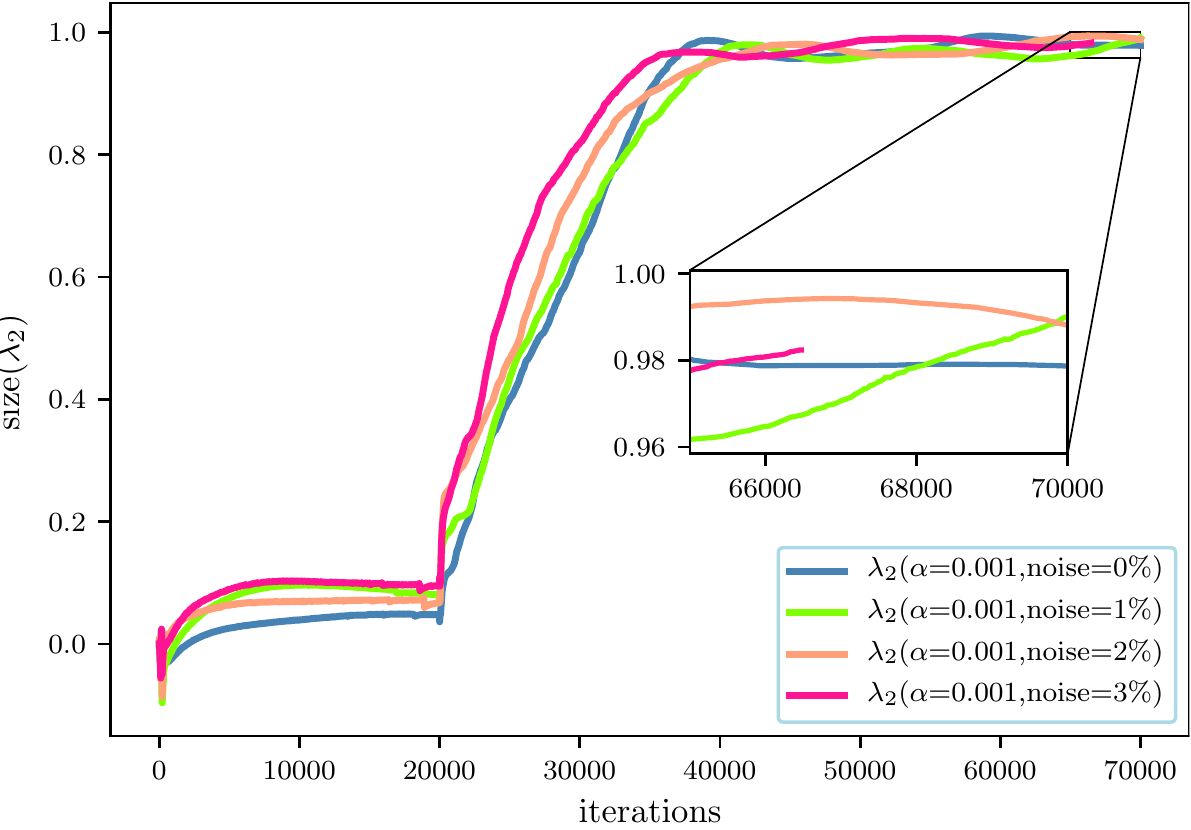}
\end{minipage}%
}\\%
\subfigure[]{
\begin{minipage}[t]{0.48\textwidth}
\centering
\includegraphics[height=4.0cm,width=6.5cm]{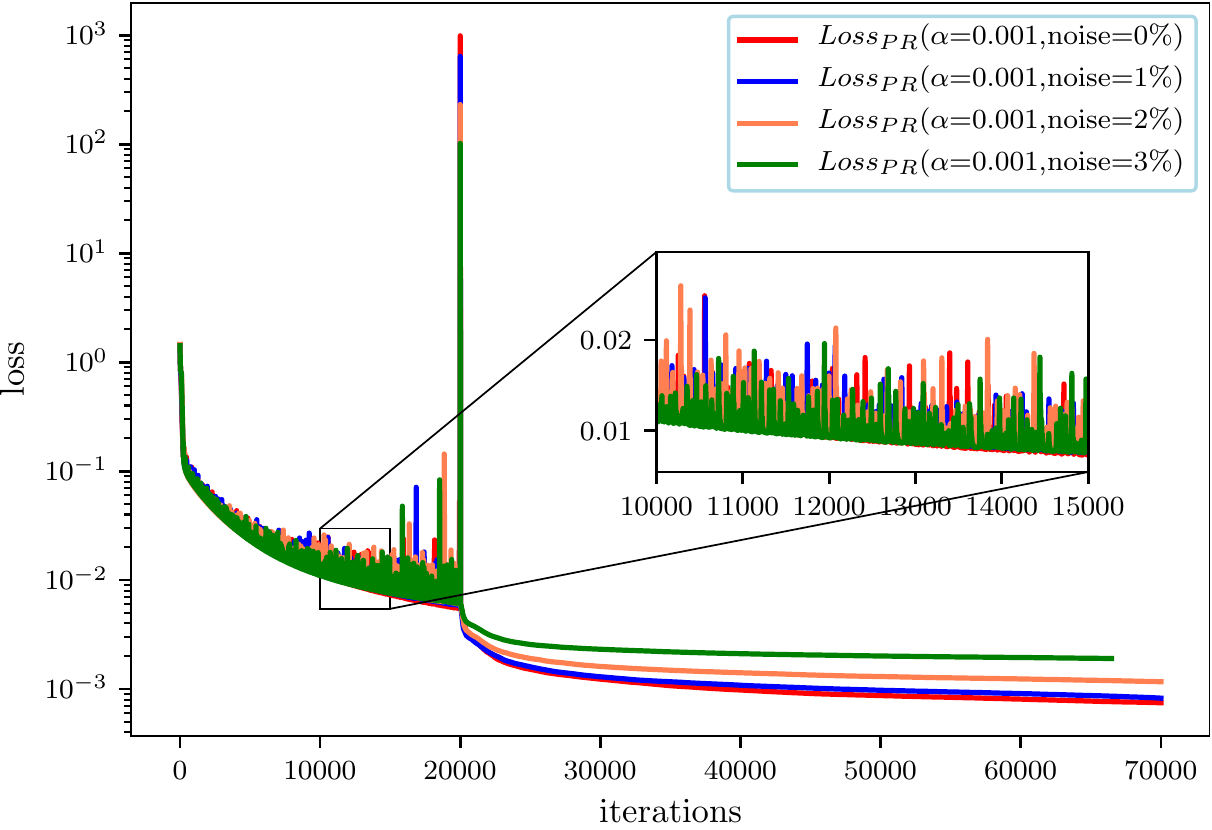}
\end{minipage}
}%
\subfigure[]{
\begin{minipage}[t]{0.48\textwidth}
\centering
\includegraphics[height=4.0cm,width=6.5cm]{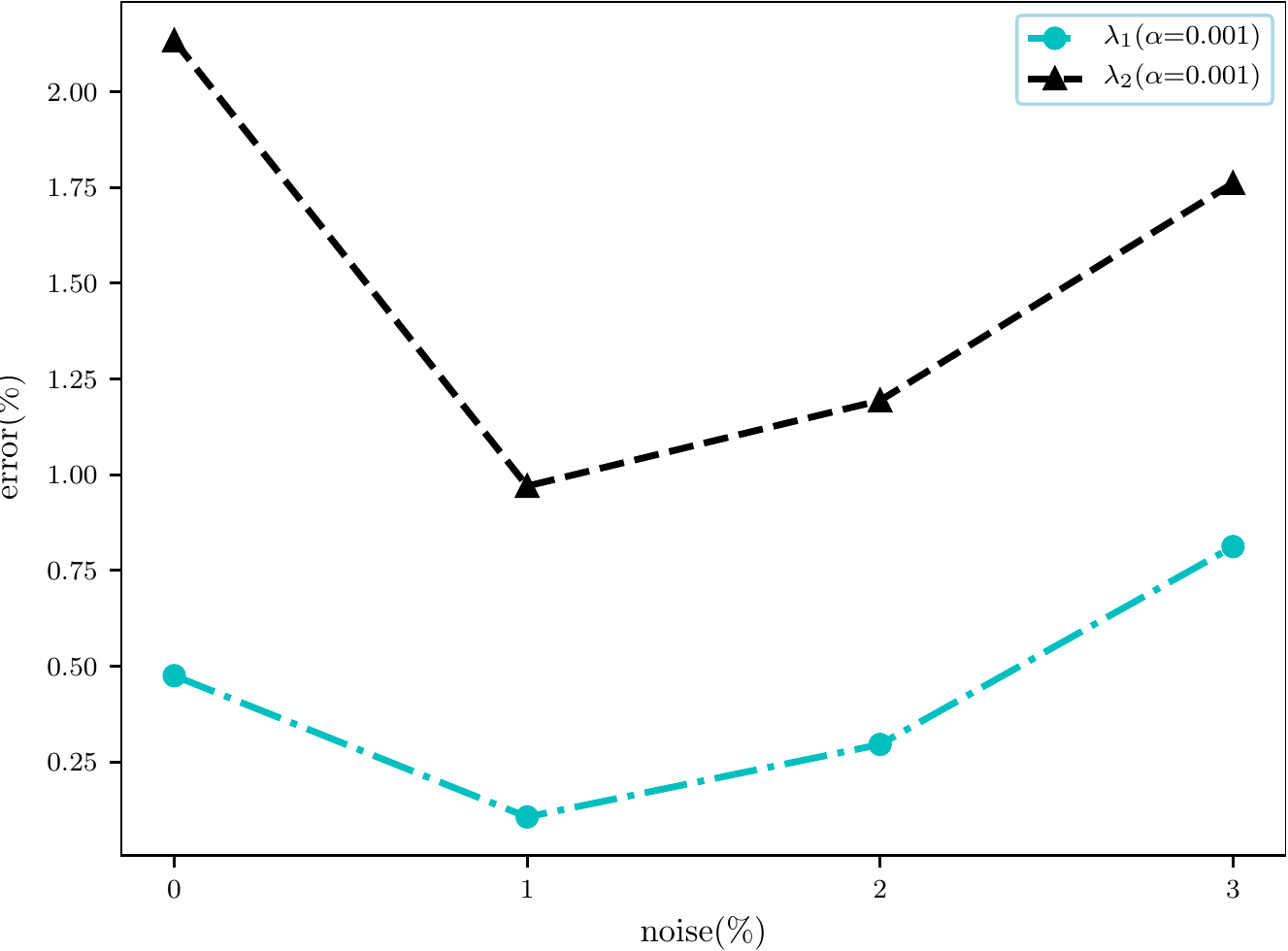}
\end{minipage}%
}%
\centering
\caption{(Color online) The parameter discover arising from the improved PINN with parameter regularization ($\alpha=0.001$): (a)-(b) the variation of unknown coefficients $\lambda_1$ and $\lambda_2$ with different noise intensity; (c) the variation of loss function with different noise intensity; (d) unknown coefficients $\lambda_1$ and $\lambda_2$ error variation under different interference noise.}
\label{F10}
\end{figure}

Furthermore, we expand the weight coefficient of parameter regularization by one order of magnitude, namely $\alpha=0.001$. Fig. \ref{F10} displays the dynamic behavior similar to that the case of $\alpha=0.0001$ in Fig. \ref{F9}, except that when the noise intensity is 1$\%$, the relative training error of parameters to be learned is smaller than that in the other three noise intensity cases. As can be seen from Fig. \ref{F10} (a) - (b), when the noise intensity is $3\%$, the NN completes the training fastest, but the number of iterations is also greater than 66000. However, when the noise intensity is $2\%$ in Fig. \ref{F9} (a) - (b), the NN completes the training fastest, and the number of iterations is less than 40000. This means that whether from the perspective of iteration times or parameter error, the NN performs best when the regularization weight coefficient is $\alpha=0.0001$, which is why we utilize the improved PINN with $\alpha=0.0001$ parametric regularization in the section 3.

\begin{figure}[htbp]
\centering
\subfigure[]{
\begin{minipage}[t]{0.48\textwidth}
\centering
\includegraphics[height=4.0cm,width=6.5cm]{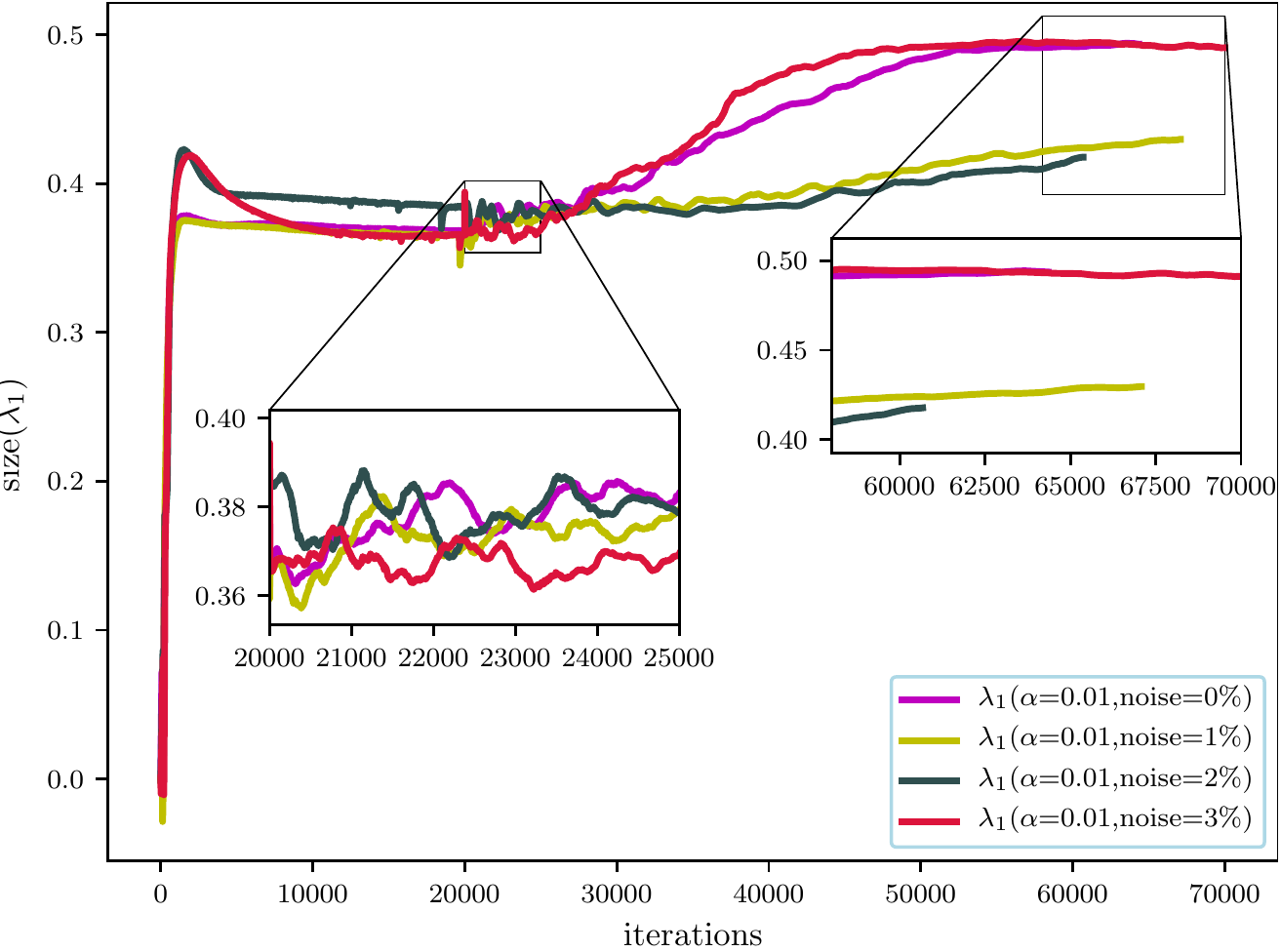}
\end{minipage}
}%
\subfigure[]{
\begin{minipage}[t]{0.48\textwidth}
\centering
\includegraphics[height=4.0cm,width=6.5cm]{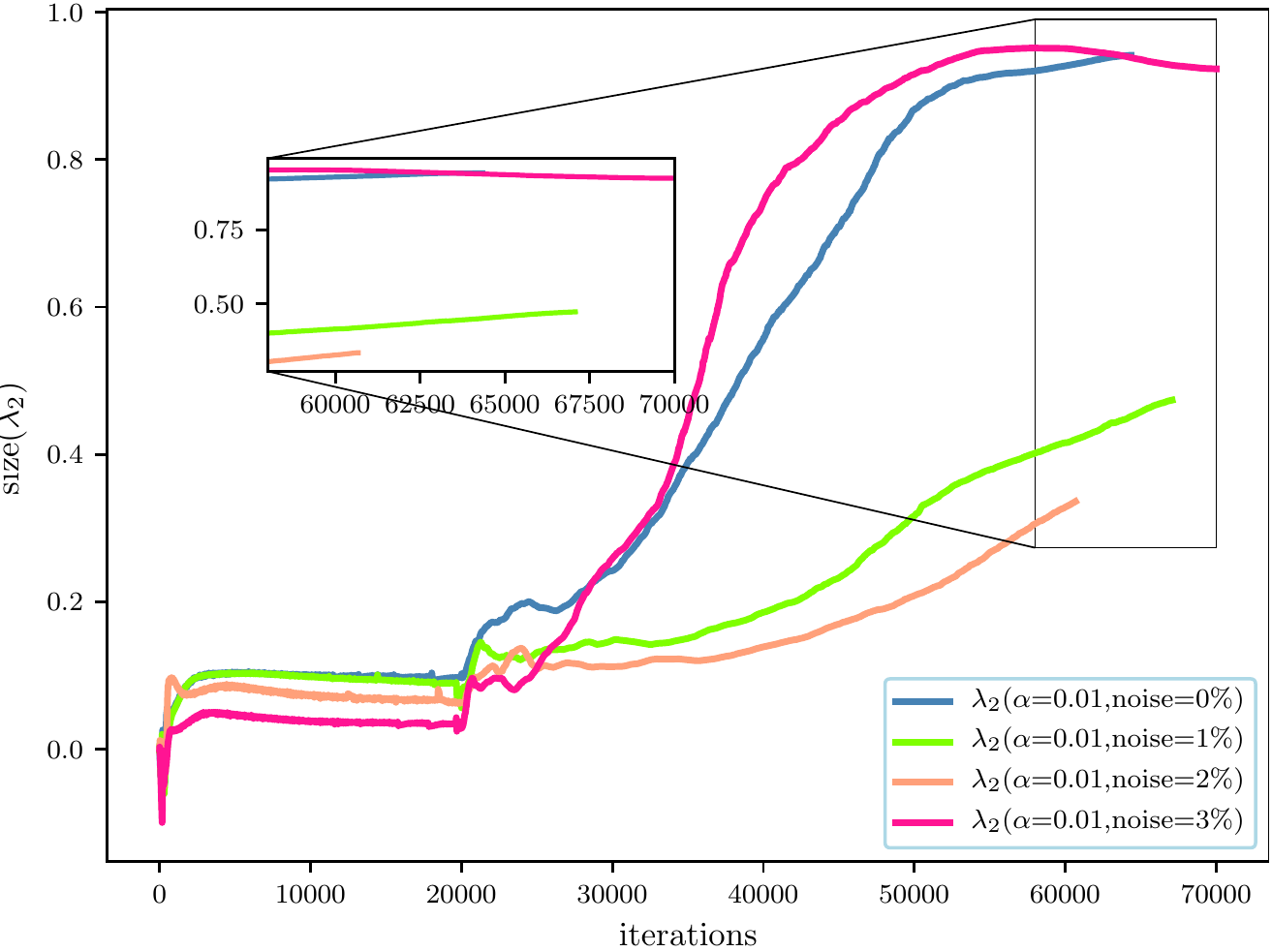}
\end{minipage}%
}\\%
\subfigure[]{
\begin{minipage}[t]{0.48\textwidth}
\centering
\includegraphics[height=4.0cm,width=6.5cm]{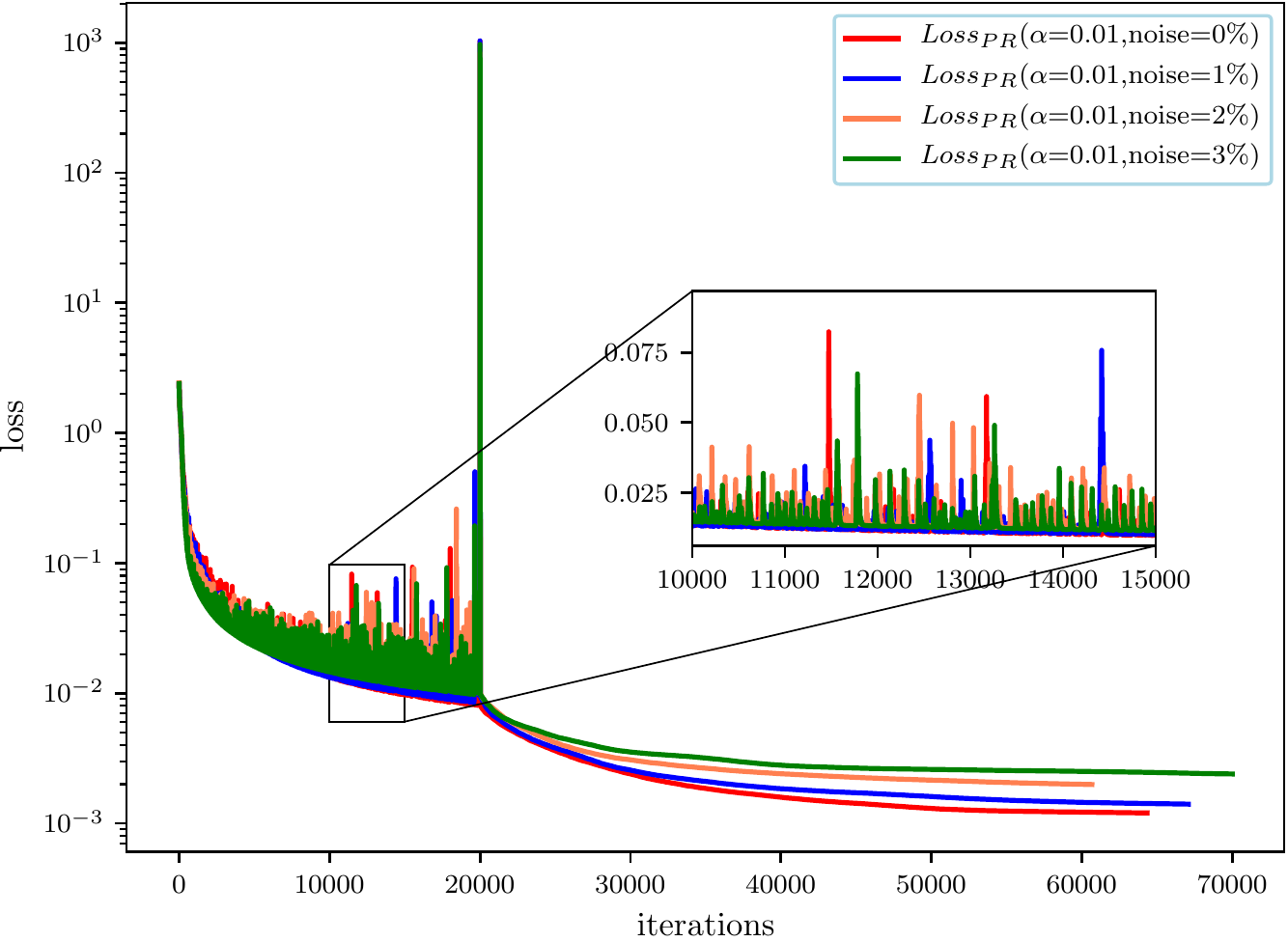}
\end{minipage}
}%
\subfigure[]{
\begin{minipage}[t]{0.48\textwidth}
\centering
\includegraphics[height=4.0cm,width=6.5cm]{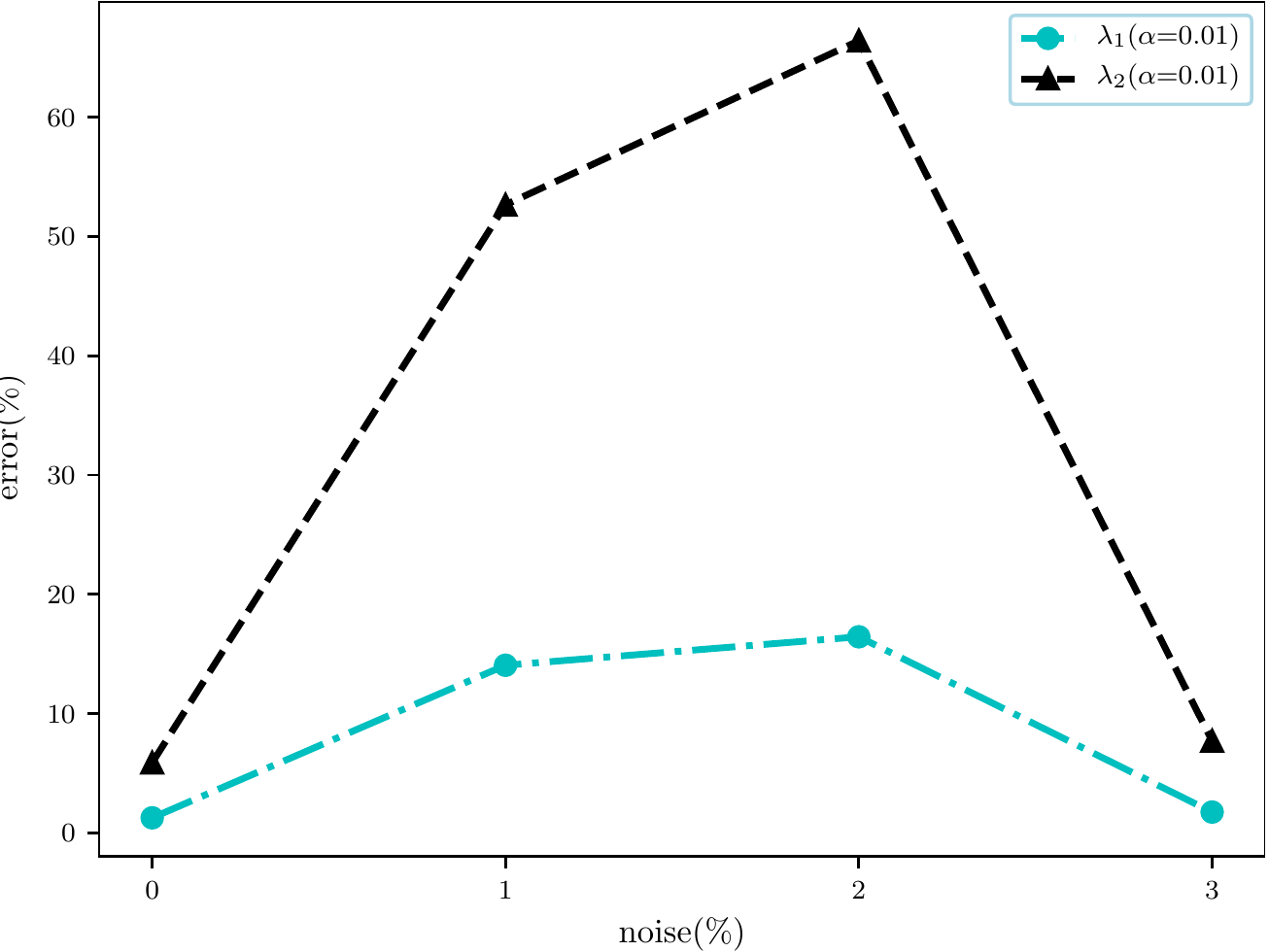}
\end{minipage}%
}%
\centering
\caption{(Color online) The parameter discover arising from the improved PINN with parameter regularization ($\alpha=0.01$): (a)-(b) the variation of unknown coefficients $\lambda_1$ and $\lambda_2$ with different noise intensity; (c) the variation of loss function with different noise intensity; (d) unknown coefficients $\lambda_1$ and $\lambda_2$ error variation under different interference noise.}
\label{F11}
\end{figure}

In order to further understand the influence of parameter regularization with larger weight ratio on improved PINN, we again expand the weight coefficient by one order of magnitude, now $\alpha=0.01$, then obtain the corresponding inverse problem training results of YO system in Fig. \ref{F11}. Fig. \ref{F11} (a)-(b) showcase the parameters discovery curves of $\lambda_1$ and $\lambda_2$, where panel (b) indicates that there is a large gap between the training value and the actual value of $\lambda_2$ when the noise intensities are 1$\%$ and 2$\%$, and combined with panel (d) of Fig. \ref{F11}, one can observe that the relative error of $\lambda_2$ is very large at this time. Although the prediction values with $3\%$ noise for the unknown parameters are close to the training results without noise, it is quite different from the real parameter values, we also notice the relative error of $\lambda_2$ has exceeded $5\%$ from Fig. \ref{F11} (d). This also means that the larger the trade-off coefficient is not always better. As shown in Fig. \ref{F11}, once it is expanded to a certain extent, the training effect is not ideal. Therefore, the selection of the value for the trade-off coefficient is very important for the parameter regularization strategy.

So to summarise, this section mainly describes the results and corresponding analysis when the improved PINN model use parameter regularization technology with different weight ratio to study the inverse problem of YO system. One can also find that when using the parameter regularization strategy with appropriate weight coefficients $\alpha$ (generally, $\alpha$ should not be too large), the training effect of parameter discovery is much better than that without parameter regularization strategy, especially after adding noise interference to the training data. In fact, these specific training results in this section further confirm the detailed analysis of parameter regularization weight coefficient $\alpha$ in section 2.4. This shows that improved PINN with $L^2$ norm parameter regularization can not only reduce the $\mathbb{L}_2$ norm error of the prediction solutions in the forward problem, but also accurately learn the unknown parameters by using the training data with noise interference in the inverse problem. Finally, we provide a training results of Figs \ref{F8} - \ref{F11} in following Tab. \ref{Tab2}.

\begin{table}[htbp]
  \caption{Comparison of correct YO system and identified YO system obtained by means of the PINN with different noise intensities and diverse weight hyper-parameters $\alpha$.}
  \label{Tab2}
  \centering
  \scalebox{0.55}{
  \begin{tabular}{l|c|c|c|c}
  \toprule
  \diagbox{\scriptsize{\textbf{YO system}}}{\scriptsize{\textbf{hyper-parameters}}} & $\alpha=0$ & $\alpha=0.0001$ & $\alpha=0.001$ & $\alpha=0.01$\\
  \hline
  Correct YO system   & \makecell[c]{$\mathrm{i}S_{t}+0.5S_{xx}+SL=0$ \\ $\mathrm{i}L_{t}-(|S|^2)_x=0$ \\ $\lambda_1$ error: 0$\%$ \\ $\lambda_2$ error: 0$\%$} & \makecell[c]{$\mathrm{i}S_{t}+0.5S_{xx}+SL=0$ \\ $\mathrm{i}L_{t}-(|S|^2)_x=0$ \\ $\lambda_1$ error: 0$\%$ \\ $\lambda_2$ error: 0$\%$} & \makecell[c]{$\mathrm{i}S_{t}+0.5S_{xx}+SL=0$ \\ $\mathrm{i}L_{t}-(|S|^2)_x=0$ \\ $\lambda_1$ error: 0$\%$ \\ $\lambda_2$ error: 0$\%$} & \makecell[c]{$\mathrm{i}S_{t}+0.5S_{xx}+SL=0$ \\ $\mathrm{i}L_{t}-(|S|^2)_x=0$ \\ $\lambda_1$ error: 0$\%$ \\ $\lambda_2$ error: 0$\%$}\\
  \hline
  Identified YO system (clean data)   & \makecell[c]{$\mathrm{i}S_{t}+0.496981S_{xx}+SL=0$ \\ $\mathrm{i}L_{t}-0.940780(|S|^2)_x=0$ \\ $\lambda_1$ error: 0.603867$\%$ \\ $\lambda_2$ error: 5.922013$\%$} & \makecell[c]{$\mathrm{i}S_{t}+0.498461S_{xx}+SL=0$ \\ $\mathrm{i}L_{t}-0.990150(|S|^2)_x=0$ \\ $\lambda_1$ error: 0.307775$\%$ \\ $\lambda_2$ error: 0.984997$\%$} & \makecell[c]{$\mathrm{i}S_{t}+0.497624S_{xx}+SL=0$ \\ $\mathrm{i}L_{t}-0.978666(|S|^2)_x=0$ \\ $\lambda_1$ error: 0.475115$\%$ \\ $\lambda_2$ error: 2.133435$\%$} & \makecell[c]{$\mathrm{i}S_{t}+0.493709S_{xx}+SL=0$ \\ $\mathrm{i}L_{t}-0.941199(|S|^2)_x=0$ \\ $\lambda_1$ error: 1.258254$\%$ \\ $\lambda_2$ error: 5.880136$\%$}\\
  \hline
  Identified YO system (1$\%$ noise)   & \makecell[c]{$\mathrm{i}S_{t}+0.487918S_{xx}+SL=0$ \\ $\mathrm{i}L_{t}-0.793893(|S|^2)_x=0$ \\ $\lambda_1$ error: 2.416390$\%$ \\ $\lambda_2$ error: 20.610661$\%$} & \makecell[c]{$\mathrm{i}S_{t}+0.496775S_{xx}+SL=0$ \\ $\mathrm{i}L_{t}-0.983750(|S|^2)_x=0$ \\ $\lambda_1$ error: 0.645000$\%$ \\ $\lambda_2$ error: 1.625025$\%$} & \makecell[c]{$\mathrm{i}S_{t}+0.500531S_{xx}+SL=0$ \\ $\mathrm{i}L_{t}-0.990297(|S|^2)_x=0$ \\ $\lambda_1$ error: 0.106263$\%$ \\ $\lambda_2$ error: 0.970280$\%$} & \makecell[c]{$\mathrm{i}S_{t}+0.429766S_{xx}+SL=0$ \\ $\mathrm{i}L_{t}-0.473574(|S|^2)_x=0$ \\ $\lambda_1$ error: 14.046884$\%$ \\ $\lambda_2$ error: 52.642570$\%$}\\
  \hline
  Identified YO system (2$\%$ noise)   & \makecell[c]{$\mathrm{i}S_{t}+0.492921S_{xx}+SL=0$ \\ $\mathrm{i}L_{t}-0.871463(|S|^2)_x=0$ \\ $\lambda_1$ error: 1.415741$\%$ \\ $\lambda_2$ error: 12.853670$\%$} & \makecell[c]{$\mathrm{i}S_{t}+0.501098S_{xx}+SL=0$ \\ $\mathrm{i}L_{t}-1.007344(|S|^2)_x=0$ \\ $\lambda_1$ error: 0.219584$\%$ \\ $\lambda_2$ error: 0.734389$\%$} & \makecell[c]{$\mathrm{i}S_{t}+0.498522S_{xx}+SL=0$ \\ $\mathrm{i}L_{t}-0.988062(|S|^2)_x=0$ \\ $\lambda_1$ error: 0.295609$\%$ \\ $\lambda_2$ error: 1.193810$\%$} & \makecell[c]{$\mathrm{i}S_{t}+0.417833S_{xx}+SL=0$ \\ $\mathrm{i}L_{t}-0.335997(|S|^2)_x=0$ \\ $\lambda_1$ error: 16.433418$\%$ \\ $\lambda_2$ error: 66.400314$\%$}\\
  \hline
  Identified YO system (3$\%$ noise)   & \makecell[c]{$\mathrm{i}S_{t}+0.487733S_{xx}+SL=0$ \\ $\mathrm{i}L_{t}-0.909239(|S|^2)_x=0$ \\ $\lambda_1$ error: 2.453375$\%$ \\ $\lambda_2$ error: 9.076095$\%$} & \makecell[c]{$\mathrm{i}S_{t}+0.496592S_{xx}+SL=0$ \\ $\mathrm{i}L_{t}-0.977557(|S|^2)_x=0$ \\ $\lambda_1$ error: 0.681674$\%$ \\ $\lambda_2$ error: 2.244323$\%$} & \makecell[c]{$\mathrm{i}S_{t}+0.495939S_{xx}+SL=0$ \\ $\mathrm{i}L_{t}-0.982388(|S|^2)_x=0$ \\ $\lambda_1$ error: 0.812221$\%$ \\ $\lambda_2$ error: 1.761216$\%$} & \makecell[c]{$\mathrm{i}S_{t}+0.491294S_{xx}+SL=0$ \\ $\mathrm{i}L_{t}-0.922920(|S|^2)_x=0$ \\ $\lambda_1$ error: 1.741284$\%$ \\ $\lambda_2$ error: 7.707989$\%$}\\
  \bottomrule
  \end{tabular}}
\end{table}

\section{Conclusions}
From many of our previous work \cite{PuND2021,PuWM2021,PuJC2021}, we realize that by introducing scalable hyper-parameters into the activation function, the PINN with locally adaptive activation functions employed to simulate localized waves of nonlinear integrable systems not only speed up the convergence of loss function, but also improves better accuracy and network performance. However, during the process of studying the inverse problem of YO system, we find that the PINN method with only neuron-wise locally adaptive activation functions has good training effect by means of clean data, but the network displays poor parameters learning capacity in the case of data training with noise interference. Therefore, considering that $L^2$ norm parameter regularization can drive the weights closer to origin and reduce the influence of large weight on the network, it can enhance the performance of NN when learning unknown parameters via training data with noisy, so we further improve the deep PINN algorithm by introducing $L^2$ norm regularization into the PINN method with neuron-wise locally adaptive activation functions. In this paper, the data-driven forward-inverse problems of YO system are showcased by means of the improved PINN with locally adaptive activation functions and parameter regularization strategy. For the data-driven forward problems of YO system, we find that improved PINN can well reveal three different forms of RWs via three distinct initial-boundary value data, including bright RW, intermediate RW and dark RW from the perspective of short wave. Surprisingly, the relative $\mathbb{L}_2$ norm errors of RWs for YO system simulated arising from improved PINN model are smaller after introducing $L^2$ norm parameter regularization technique into PINN, although the effect is not obvious in some training results, as shown in Tab. \ref{Tab1}. However, during the further research of the inverse problem for YO system, we find that the data-driven unknown parameters trained by improved PINN with $L^2$ norm parameter regularization are more accurate than those trained by means of the PINN without parameter regularization. Compared with the general PINN model with only neuron-wise locally adaptive activation functions, the improved PINN with parameter regularization shows excellent noise immunity in the inverse problem of YO system.

Since the parameter regularization strategy is a very direct and efficient approach to improve the performance of NNs, and the introduction of $L^2$ norm regularization strategy in deep learning is a very mature and commonly used means. Naturally, we successfully applied this technology to PINN algorithm and found that it has a positive effect in the forward problem of YO system, and achieves an amazing effect in the parameter discovery process of inverse problem in this paper. Except for $L^2$ norm parameter regularization, there are other parameter regularization strategies, such as $L^1$ norm parameter regularization. In comparison to $L^2$ regularization, $L^1$ regularization results in a solution that is more sparse, here sparsity in this context refers to the fact that some parameters have an optimal value of zero. The $L^2$ regularization does not cause the parameters to become sparse, while $L^1$ regularization may do so for large enough $\alpha$. Indeed, we also introduce $L^1$ norm parameter regularization strategy into PINN model to investigate forward-inverse problems of YO system, but we find that the simulation effect is not ideal, so how to better combine $L^1$ parameter regularization technology with PINN to deal with the various data-driven problems of nonlinear systems needs further research. Of course, there are other parameter regularization methods, but how to properly introduce these regularization methods into deep learning methods to solve the problem at hand is an eternal topic. We anticipate that our work will impact significantly on the design of deep learning experiments for various nonlinear systems and may open new research opportunities for integrable deep learning.

\section*{Declaration of competing interest}
The authors declare that they have no known competing financial interests or personal relationships that could have appeared to influence the work reported in this paper.

\section*{Acknowledgements}
\hspace{0.3cm}
The authors gratefully acknowledge the support of the National Natural Science Foundation of China (No. 12175069) and Science and Technology Commission of Shanghai Municipality (No. 21JC1402500 and No. 18dz2271000).


\begin{thebibliography}{99}

\bibitem{Hinton2006}
G.E. Hinton, S. Osindero, Y. Teh, A fast learning algorithm for deep belief nets, Neural Computation, 18(17) (2006) 1527-1554.

\bibitem{LeCun2015}
Y. LeCun, Y. Bengio, G. E. Hinton, Deep learning, Nature 521 (2015) 436-444.

\bibitem{Krizhevsky2012}
A. Krizhevsky, I. Sutskever, G. E. Hinton, ImageNet classification with deep convolutional neural networks, In Advances in Neural Information Processing Systems 25 (NIPS’2012) (2012).

\bibitem{Cao2020}
X.Y. Cao, J. Yao, Z.B. Xu, D.Y. Meng, Hyperspectral image classification with convolutional neural network and active learning, IEEE Trans. Geosci. Remote Sens. 58(7) (2020) 4604-4616.

\bibitem{Hinton2012}
G.E. Hinton, L. Deng, D. Yu, G.E. Dahl, A. Mohamed, N. Jaitly, A. Senior, V. Vanhoucke, P. Nguyen, T.N. Sainath, B. Kingsbury, Deep neural networks for acoustic modeling in speech recognition: The shared views of four research groups, IEEE Signal Process. Mag. 29(6) (2012) 82-97.

\bibitem{LiGF2020}
G.F. Li, Y.F. Yang, X.D. Qu, Deep learning approaches on pedestrian detection in hazy weather, IEEE Trans. Ind. Electron. 67(10) (2020) 8889-8899.

\bibitem{ZengNY2021}
N.Y. Zeng, H. Li, Z.D. Wang, W.B. Liu, S.M. Liu, F.E. Alsaadi, X.H. Liu, Deep-reinforcement-learning-based images segmentation for quantitative analysis of gold immunochromatographic strip, Neurocomputing 425(15) (2021) 173-180.

\bibitem{Zhang2020}
J.M. Zhang, W. Wang, C.Q. Lu, J. Wang, A.K. Sangaiah, Lightweight deep network for traffic sign classification, Ann. Telecommun. 75 (2020) 369-379.

\bibitem{Collobert2011}
R. Collobert, J. Weston, L. Bottou, M. Karlen, K. Kavukcuoglu, P. Kuksa, Natural language processing (almost) from scratch, J. Mach. Learn. Res. 12 (2011) 2493-2537.

\bibitem{Alipanahi2015}
B. Alipanahi, A. Delong, M.T. Weirauch, B.J. Frey, Predicting the sequence specificities of DMA- and RNA-binding proteins by deep learning, Nat. Biotechnol. 33 (2015) 831-838.

\bibitem{Sun2018}
X.D. Sun, P.C. Wu, S.C.H. Hoi, Face detection using deep learning: an improved faster RCNN approach, Neurocomputing 299 (2018) 42-50.

\bibitem{Shao2020}
Z.F. Shao, L.G. Wang, Z.Y. Wang, W. Du, W.J. Wu, Saliency-aware convolution neural network for ship detection in surveillance video, IEEE Trans. Circuits Syst. Video Technol. 30(3) (2020) 781-794.

\bibitem{Raissi2019}
M. Raissi, P. Perdikaris, G.E. Karniadakis, Physics-informed neural networks: a deep learning framework for solving forward and inverse problems involving nonlinear partial differential equations, J. Comput. Phys. 378 (2019) 686-707.

\bibitem{Haykin2008}
S. Haykin, Neural Networks and Learning Machines, Prentice-Hall, New York, 2008.

\bibitem{Hornik1989}
K. Hornik, M. Stinchcombe, H. White, Multilayer feedforward networks are universal approximators, Neural Networks 2(5) (1989) 359-366.

\bibitem{Rumelhart1986}
D.E. Rumelhart, G.E. Hinton, R.J. Williams, Learning representations by back-propagating errors, Nature, 323 (1986) 533-536.

\bibitem{Ruder2017}
S. Ruder, An overview of gradient descent optimization algorithms, arXiv:1609 .04747v2, 2017.

\bibitem{Kingma2014}
D.P. Kingma, J.L. Ba, Adam: a method for stochastic optimization, arXiv: 1412.6980v9, 2014.

\bibitem{Liu1989}
D.C. Liu, J. Nocedal, On the limited memory BFGS method for large scale optimization, Math. Program. 45 (1989) 503-528.

\bibitem{Baydin2018}
A.G. Baydin, B.A. Pearlmutter, A.A. Radul, J.M. Siskind, Automatic differentiation in machine learning: a survey, Journal of Machine Learning Research 18 (2018) 1-43.

\bibitem{Lagaris1998}
I.E. Lagaris, A. Likas, D.I. Fotiadis, Artificialneural network for solving ordinary and partial differential equations, IEEE Trans. Neural Netw. 9(5) (1998) 987-1000.

\bibitem{KarniadakisGE2021}
G.E. Karniadakis, I.G. Kevrekidis, L. Lu, P. Perdikaris, S.F. Wang, L. Yang, Physics-informed machine learning, Nat. Rev. Phys. 3 (2021) 422-440.

\bibitem{RaissiM2019}
M. Raissi, H. Babaee, P. Givi, Deep learning of turbulent scalar mixing, Phys. Rev. Fluids 4 (2019) 124501.

\bibitem{PangG2019}
G.F. Pang, L. Lu, G.E. Karniadakis, FPINNs: fractional physics-informed neural networks, SIAM J. Sci. Comput. 41(4) (2019) A2603-A2626.

\bibitem{Mao2020}
Z.P. Mao, A.D. Jagtap, G.E. Karniadakis, Physics-informed neural networks for high-speed flows, Comput. Methods Appl. Mech. Engrg. 360 (2020) 112789.

\bibitem{CaiS2021}
S.Z. Cai, Z.C. Wang, S.F. Wang, P. Perdikaris, G.E. Karniadakis, Physics-informed neural networks for heat transfer problems. J. Heat Transfer, 143(6) (2021) 060801.

\bibitem{YangL2020}
L. Yang, D.K. Zhang, G.E. Karniadakis, Physics-informed generative adversarial networks for stochastic differential equations, SIAM J. Sci. Comput. 42(1) (2020) A292-A317.

\bibitem{Goodfellow-book2016}
I.J. Goodfellow, Y. Bengio, A. Courville, Deep Learning, MIT Press, Cambridge, 2016.

\bibitem{Jagtap2020}
A.D. Jagtap, K. Kawaguchi, G.E. Karniadakis, Adaptive activation functions accelerate convergence in deep and physics-informed neural networks, J. Comput. Phys. 404 (2020) 109136.

\bibitem{JagtapA2020}
A.D. Jagtap, K. Kawaguchi, G.E. Karniadakis, Locally adaptive activation functions with slope recovery for deep and physics-informed neural networks, Proc. R. Soc. A 476 (2020) 20200334.

\bibitem{PuND2021}
J.C. Pu, J. Li, Y. Chen, Solving localized wave solutions of the derivative nonlinear Schr\"odinger equation using an improved PINN method, Nonlinear Dyn. 105 (2021) 1723-1739.

\bibitem{PuWM2021}
J.C. Pu, W.Q. Peng, Y. Chen, The data-driven localized wave solutions of the derivative nonlinear Schr\"{o}dinger equation by using improved PINN approach, Wave Motion 107 (2021) 102823.

\bibitem{Poggio1990}
T. Poggio, F. Girosi, Regularization algorithms for learning that are equivalent to multilayer networks, Science, 247(4945) (1990) 978-982.

\bibitem{ZhangR2021}
R. Zhang, X.L. Li, T. Xu, Y. Zhao, Data clustering via uncorrelated ridge regression, IEEE Trans. Neural Netw. Learn. Syst. 32(1) (2021) 450-456.

\bibitem{Chada2020}
N.K. Chada, A.M. Stuart, X.T. Tong, Tikhonov regularization within ensemble kalman inversion, SIAM J. Numer. Anal. 58(2) (2020) 1263-1294.

\bibitem{PuCPB2021}
J.C. Pu, J. Li, Y. Chen, Soliton, breather and rogue wave solutions for solving the nonlinear Schr\"odinger equation using a deep learning method with physical constraints, Chin. Phys. B 30 (2021) 060202.

\bibitem{Peng2021}
W.Q. Peng, J.C. Pu, Y. Chen, PINN deep learning for the Chen-Lee-Liu equation: rogue wave on the periodic background, Commun. Nonlinear Sci. Numer. Simul. 105 (2022) 106067.

\bibitem{Manakov1974}
S.V. Manakov, On the theory of two-dimensional stationary self-focusing of electromagnetic waves, Sov. Phys. JETP 38(2) (1974) 248-253.

\bibitem{Benney1977}
D.J. Benney, A general theory for interactions between short and long waves, Stud. Appl. Math. 56(1) (1977) 81-94.

\bibitem{Yajima1976}
N. Yajima, M. Oikawa, Formation and interaction of Sonic-Langmuir solitons: inverse scattering method, Prog. Theor. Phys. 56(6) (1976) 1719-1739.

\bibitem{Djordjevic1977}
V.D. Djordjevic, L.G. Redekopp, On the two-dimensional packets of capillary-gravitywaves, J. Fluid Mech. 79 (1977) 703-714.

\bibitem{Zakharov1972}
V.E. Zakharov, Collapse of Langmuir waves, Sov. Phys. JETP 35(5) (1972) 908-914.

\bibitem{Chowdhury2008}
A. Chowdhury, J.A. Tataronis, Long wave-short wave resonance in nonlinear negative refractive index media, Phys. Rev. Lett. 100 (2008) 153905.

\bibitem{Grimshaw1977}
R.H.J. Grimshaw, The modulation of an internal gravity-wave packet, and the resonance with the mean motion, Stud. Appl. Math. 56(3) (1977) 241-266.

\bibitem{Funakoshi1983}
M. Funakoshi, M. Oikawa, The resonant interaction between a long internal gravity wave and a surface gravity wave packet, J. Phys. Soc. Jpn. 52 (1983) 1982-1995.

\bibitem{PuJC2021}
J.C. Pu, Y. Chen, The data-driven vector localized waves of Manakov system using improved PINN approach, arXiv: 2109.09266, 2021.

\bibitem{Chow2013}
K.W. Chow, H.N. Chan, D.J. Kedziora, R.H.J. Grimshaw, Rogue wave modes for the long wave-short wave resonance model, J. Phys. Soc. Jpn. 82 (2013) 074001.

\bibitem{ChenJC2018}
J.C. Chen, Y. Chen, B.F. Feng, K.I. Maruno, Y. Ohta, General high-order rogue waves of the (1+1)-dimensional Yajima-Oikawa System, J. Phys. Soc. Jpn. 87 (2018) 094007.

\bibitem{Stein1987}
M. Stein, Large sample properties of simulations using Latin hypercube sampling, Technometrics 29(2) (1987) 143-151.

\bibitem{ChenSH2014}
S.H. Chen, Darboux transformation and dark rogue wave states arising from two-wave resonance interaction, Phys. Lett. A 378 (2014) 1095-1098.

\bibitem{Peregrine1983}
D.H. Peregrine, Water waves, nonlinear Schr\"{o}dinger equations and their solutions, J. Austral. Math. Soc. Ser. B Appl. Math. 25(1) (1983) 16-43.




\end{thebibliography}
\end{document}